\documentclass[12pt]{amsart}
\usepackage{amsmath}
\usepackage{amsthm}
\usepackage{amsfonts}
\usepackage{amssymb}
\usepackage{mathrsfs}
\usepackage{graphicx}
\usepackage{stmaryrd}
\usepackage{dsfont}
\usepackage[all]{xy}
\usepackage[shortlabels]{enumitem}
\usepackage[usenames, dvipsnames]{color}
\usepackage[margin=1in]{geometry} 
\usepackage[bookmarks, bookmarksdepth=2, colorlinks=true, linkcolor=blue, citecolor=blue, urlcolor=blue]{hyperref}
\usepackage{eucal}
\usepackage{contour}
\usepackage[normalem]{ulem}
\usepackage{tikz}

\makeatletter
\@addtoreset{equation}{section}
\makeatother

\numberwithin{equation}{section}
\newtheorem{theorem}[equation]{Theorem}

\newtheorem{proposition}[equation]{Proposition}
\newtheorem{lemma}[equation]{Lemma}
\newtheorem{corollary}[equation]{Corollary}

\theoremstyle{definition}
\newtheorem{rmk}[equation]{Remark}
\newenvironment{remark}[1][]{\begin{rmk}[#1] \pushQED{\qed}}{\popQED \end{rmk}}
\newtheorem{eg}[equation]{Example}
\newenvironment{example}[1][]{\begin{eg}[#1] \pushQED{\qed}}{\popQED \end{eg}}
\newtheorem{defnaux}[equation]{Definition}
\newenvironment{definition}[1][]{\begin{defnaux}[#1]\pushQED{\qed}}{\popQED \end{defnaux}}

\newcommand{\arxiv}[1]{\href{http://arxiv.org/abs/#1}{{\tiny\tt arXiv:#1}}}
\newcommand{\DOI}[1]{\href{http://doi.org/#1}{\color{purple}{\tiny\tt DOI:#1}}}

\newcommand{\bC}{\mathbf{C}}
\newcommand{\cC}{\mathcal{C}}

\newcommand{\cD}{\mathcal{D}}

\newcommand{\bF}{\mathbf{F}}

\newcommand{\rH}{\mathrm{H}}

\newcommand{\bI}{\mathbf{I}}

\newcommand{\rI}{\mathrm{I}}

\newcommand{\bJ}{\mathbf{J}}

\newcommand{\bK}{\mathbf{K}}

\newcommand{\rK}{\mathrm{K}}

\newcommand{\bN}{\mathbf{N}}

\newcommand{\bO}{\mathbf{O}}

\newcommand{\bQ}{\mathbf{Q}}

\newcommand{\bR}{\mathbf{R}}

\newcommand{\bS}{\mathbf{S}}

\newcommand{\fS}{\mathfrak{S}}

\newcommand{\bZ}{\mathbf{Z}}

\newcommand{\rf}{\mathrm{f}}

\let\ol\overline
\let\ul\underline

\newcommand{\defn}[1]{\emph{#1}}
\renewcommand{\phi}{\varphi}
\renewcommand{\emptyset}{\varnothing}
\newcommand{\lw}{{\textstyle \bigwedge}}
\DeclareMathOperator{\sgn}{sgn} 
\DeclareMathOperator{\im}{im}

\DeclareMathOperator{\End}{End}

\DeclareMathOperator{\Aut}{Aut}

\DeclareMathOperator{\Mod}{Mod}
\DeclareMathOperator{\Ind}{Ind}
\DeclareMathOperator{\ind}{ind}
\DeclareMathOperator{\res}{res}
\DeclareMathOperator{\Hom}{Hom}

\DeclareMathOperator{\gr}{gr}
\DeclareMathOperator{\Rep}{Rep}
\DeclareMathOperator{\vol}{vol}

\newcommand{\id}{\mathrm{id}}
\newcommand{\op}{\mathrm{op}}

\DeclareMathOperator{\Res}{Res}
\newcommand{\GL}{\mathbf{GL}}
\newcommand{\Mat}{\mathrm{Mat}}
\newcommand{\Vect}{\mathrm{Vec}}
\newcommand{\ev}{\mathrm{ev}}

\newcommand{\bone}{\mathbf{1}}
\newcommand{\bbone}{\mathds{1}}
\newcommand{\uotimes}{\mathbin{\ul{\otimes}}}
\newcommand{\OI}{\mathbf{OI}}

\newcommand{\wa}{{\bullet}}
\newcommand{\wb}{{\circ}}

\contourlength{1pt}

\contourlength{0.8pt}
\newcommand{\myuline}[1]{%
  \uline{\phantom{#1}}%
  \llap{\contour{white}{#1}}%
}
\DeclareMathOperator{\uRep}{\text{\myuline{\rm Rep}}}
\DeclareMathOperator{\uPerm}{\ul{Perm}}

\title{The Delannoy category}
\author{Nate Harman}
\author{Andrew Snowden}
\author{Noah Snyder}

\date{January 6, 2023}

\begin{document}

\begin{abstract}
Let $G$ be the group of all order-preserving self-maps of the real line. In previous work, the first two authors constructed a pre-Tannakian category $\uRep(G)$ associated to $G$. The present paper is a detailed study of this category, which we name the \emph{Delannoy category}. We classify the simple objects, determine branching rules to open subgroups, and give a combinatorial rule for tensor products. The Delannoy category has some remarkable features: it is semi-simple in all characteristics; all simples have categorical dimension $\pm 1$; and the Adams operations on its Grothendieck group are trivial. We also give a combinatorial model for $\uRep(G)$ based on Delannoy paths.
\end{abstract}

\maketitle
\tableofcontents

\section{Introduction} \label{s:intro}

Pre-Tannakian tensor categories provide a natural generalization of representation categories of algebraic groups or supergroups. A topic of recent interest is to understand to what extent pre-Tannakian categories go beyond those coming from algebraic supergroups. An important theorem of Deligne \cite{Deligne2} states that any pre-Tannakian category of moderate growth in characteristic zero comes from an algebraic supergroup, and there has been much work devoted to generalizing this result to positive characteristic \cite{BE, BEO, Cou1, AbEnv, CEO, CEOP, EtingofHarmanOstrik, EOf, EOV, Ostrik}. Deligne \cite{Deligne3} also introduced the first examples of pre-Tannakian categories not associated to algebraic supergroups: the interpolation categories $\uRep(\fS_t)$, $\uRep(\GL_t)$, $\uRep(\bO_t)$; see  \cite{ComesOstrik1, ComesOstrik, EntovaAizenbudHeidersdorf, Harman, Harman2, Knop, Knop2} for subsequent work.

In recent work \cite{HarmanSnowden}, the first two authors introduced a general construction of pre-Tannakian categories associated to oligomorphic groups with an appropriate measure. In the case of the infinite symmetric group, this construction recovers Deligne's category $\uRep(\fS_t)$. The construction also applies to the oligomorphic group $G=\Aut(\bR, <)$ of order-preserving self-bijections of the real line. The resulting category $\uRep(G)$ is a fundamentally new example of a tensor category: for example, in positive characteristic, it is the first known example of a semi-simple pre-Tannakian category that is not of moderate growth. We note that \cite{HarmanSnowden} is mostly concerned with developing the general theory of pre-Tannkian categories associated to oligomorphic groups, and proved little about $\uRep(G)$ in particular.


This paper is a detailed study of the category $\uRep(G)$. We christen this the \defn{Delannoy category}, for reasons that will be shortly apparent. Our results show that it is a remarkable object in possession of a number of features distinguishing it from other known pre-Tannakian categories. In the rest of the introduction, we briefly describe what $\uRep(G)$ is, and then explain a number of our results about it.

\subsection{The Delannoy category}

From here on, we fix a field $k$ and write $G$ for $\Aut(\bR,<)$. We now give a brief discussion of what $\uRep(G)$ is. We hope this provides the reader with enough of a sense of the category to appreciate our main results. A more detailed discussion of the construction of $\uRep(G)$ is provided in \S \ref{s:oligo} and \S \ref{s:AutR}.

Let $\bR^{(n)}$ be the subset of $\bR^n$ consisting of tuples $(x_1, \ldots, x_n)$ with $x_1<\cdots<x_n$. The group $G$ acts transitively on $\bR^{(n)}$. We say that a function $\phi \colon \bR^{(n)} \to k$ is a \defn{Schwartz function} if its level sets can be defined by first-order formulas using $<$, $=$, and finitely many real constants; equivalently, this means that the subgroup of $G$ stabilizing $\phi$ is open, with respect to the natural topology on $G$. We define the \defn{Schwartz space} $\cC(\bR^{(n)})$ to be the $k$-vector space of all Schwartz functions.

A concept of fundamental importance to this paper is integration with respect to Euler characteristic, as developed by Schapira and Viro (see, e.g., \cite{Viro}). Precisely, if $\phi$ is a Schwartz function on $\bR^{(n)}$ that takes the values $a_1, \ldots, a_r$ then we define
\begin{displaymath}
\int_{\bR^{(n)}} \phi(x) dx = \sum_{i=1}^r a_i \cdot \chi_c(\phi^{-1}(a_i)),
\end{displaymath}
where $\chi_c$ is the Euler characteristic of compactly supported cohomology. For example, if $\phi \in \cC(\bR)$ is the characteristic function of an open interval, its integral is $-1$.

Suppose that $\Psi$ is a Schwartz function on $\bR^{(n)} \times \bR^{(m)}$ (defined in the obvious manner). We can then define an integral operator
\begin{displaymath}
A_{\Psi} \colon \cC(\bR^{(m)}) \to \cC(\bR^{(n)}), \qquad (A_{\Psi} \phi)(x) = \int_{\bR^{(m)}} \Psi(x,y) \phi(y) dy.
\end{displaymath}
The composition of two such integral operators is an integral operator of the same kind.

We can use the above constructions to define a category $\uPerm(G)$ of ``permutation modules.'' The objects are formal direct sums of $\cC(\bR^{(n)})$'s. The morphisms of basic objects are $G$-invariant integral operators. This category carries a symmetric tensor structure $\uotimes$ defined on the basic objects by
\begin{displaymath}
\cC(\bR^{(n)}) \uotimes \cC(\bR^{(m)}) = \cC(\bR^{(n)} \times \bR^{(m)}).
\end{displaymath}
The space $\bR^{(n)} \times \bR^{(m)}$ decomposes into $G$-orbits isomorphic to $\bR^{(s)}$ (for various $s$), and the right side above is defined to be the corresponding sum of Schwartz spaces. The tensor category $\uPerm(G)$ is rigid: in fact, every object is self-dual.

The category $\uRep^{\rf}(G)$ of finite length objects in $\uRep(G)$ is equivalent to the Karoubian envelope of $\uPerm(G)$, and $\uRep(G)$ is equivalent to the ind-completion of $\uRep^{\rf}(G)$; these statements rely on non-trivial results of \cite{HarmanSnowden}, but for the present purposes one can take them as definitions of $\uRep^{\rf}(G)$ and $\uRep(G)$. As we have already indicated, $\uRep^{\rf}(G)$ is a semi-simple pre-Tannakian tensor category; this is again a non-trivial result of \cite{HarmanSnowden}. We emphasize that semi-simplicity holding in all characteristics is rather remarkable.

\subsection{The path model} \label{ss:intro-path}

Before stating our results, we give another description of the Delannoy category that highlights the relevant combinatorics. This description is new to this paper, and developed in more detail in \S \ref{s:delannoy}.

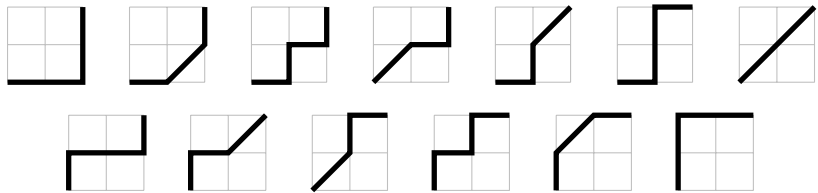
\begin{figure}
\def\delannoy#1{\begin{tikzpicture}[scale=0.5]
\draw[step=1, color=gray!50] (0, 0) grid (2,2);
\draw[line width=2pt] #1;
\end{tikzpicture}}
\begin{center}
\delannoy{(0,0)--(1,0)--(2,0)--(2,1)--(2,2)} \quad
\delannoy{(0,0)--(1,0)--(2,1)--(2,2)} \quad
\delannoy{(0,0)--(1,0)--(1,1)--(2,1)--(2,2)} \quad
\delannoy{(0,0)--(1,1)--(2,1)--(2,2)} \quad
\delannoy{(0,0)--(1,0)--(1,1)--(2,2)} \quad
\delannoy{(0,0)--(1,0)--(1,1)--(1,2)--(2,2)} \quad
\delannoy{(0,0)--(1,1)--(2,2)}
\end{center}
\vskip.5\baselineskip
\begin{center}
\delannoy{(0,0)--(0,1)--(1,1)--(2,1)--(2,2)} \quad
\delannoy{(0,0)--(0,1)--(1,1)--(2,2)} \quad
\delannoy{(0,0)--(1,1)--(1,2)--(2,2)} \quad
\delannoy{(0,0)--(0,1)--(1,1)--(1,2)--(2,2)} \quad
\delannoy{(0,0)--(0,1)--(1,2)--(2,2)} \quad
\delannoy{(0,0)--(0,1)--(0,2)--(1,2)--(2,2)}
\end{center}
\caption{The thirteen $(2,2)$-Delannoy paths.}
\label{fig:delannoy}
\end{figure}

An \defn{$(n,m)$-Delannoy path} is a path in the plane from $(0,0)$ to $(n,m)$ composed of steps of the form $(1,0)$, $(0,1)$, and $(1,1)$. The \defn{Delannoy number} $D(n,m)$ is the number of $(n,m)$-Delannoy paths, and the \defn{central Delannoy number} $D(n)$ is $D(n,n)$. For example, $D(2)=13$; see Figure~\ref{fig:delannoy}. The Delannoy numbers are well-known in the literature; see, e.g., \cite{Banderier}.

We introduce an algebra $\cD(n)$, called the \defn{Delannoy algebra}, which is modeled on Delannoy paths. The elements of $\cD(n)$ are $k$-linear combinations of $(n,n)$-Delannoy paths; we write $[p]$ for the basis element of $\cD(n)$ corresponding to the path $p$. Suppose that $p_1$, $p_2$, and $p_3$ are three $(n,n)$-Delannoy paths. It turns out that there is at most one three-dimensional $(n,n,n)$-Delannoy path $q$ whose three projections are $p_1$, $p_2$, and $p_3$. We define
\begin{displaymath}
\epsilon(p_1,p_2,p_3) = \begin{cases}
(-1)^{\ell(q)+\ell(p_3)} & \text{if $q$ exists} \\
0 & \text{otherwise} \end{cases}
\end{displaymath}
where $\ell(q)$ is the length of $q$. Multiplication in $\cD(n)$ is then defined by
\begin{displaymath}
[p_1] \cdot [p_2] = \sum_{p_3} \epsilon(p_1,p_2,p_3) [p_3].
\end{displaymath}
This product is associative and unital. We show that $\cD(n)$ is isomorphic to the endomorphism algebra of $\cC(\bR^{(n)})$ in $\uRep(G)$; in particular, the dimension of this endomorphism algebra is the central Delannoy number $D(n)$. The key point here is that $G$-orbits on $\bR^{(n)} \times \bR^{(n)}$ are naturally paremetrized by $(n,n)$-Delannoy paths; see Proposition~\ref{prop:orbit}.

More generally, define $\cD(n,m)$ to be the space spanned by $(n,m)$-Delannoy paths. A similar combinatorial rule to the above defines a composition
\begin{displaymath}
\cD(n,m) \times \cD(m,\ell) \to \cD(n,\ell).
\end{displaymath}
We define a category with objects indexed by natural numbers and where the $\Hom$ sets are $\cD(n,m)$. We show that the additive envelope of this category is $\uPerm(G)$. We also explain how to construct the tensor product from the combinatorial point of view. See \S \ref{s:delannoy} for details.

\subsection{Main results} \label{ss:results}

Having given some sense of $\uRep(G)$, we now explain our results.

\textit{(a) Classification of simples.} A \defn{weight} is a word in the alphabet $\{\wa,\wb\}$. Given a weight $\lambda$ of length $n$, we define a submodule $L_{\lambda}$ of $\cC(\bR^{(n)})$ by writing down an explicit family of Schwartz functions and taking the submodule they generate. We show that the $L_{\lambda}$'s are exactly the simple objects of $\uRep(G)$ (Theorem~\ref{thm:simple} and Corollary~\ref{cor:all-simples}).

\textit{(b) Decomposition of Schwartz space.} We determine the simple decomposition of $\cC(\bR^{(n)})$ (Theorem~\ref{thm:decomp}). In particular, we find that it has length $3^n$. This decomposition gives a representation-theoretic meaning to a classical formula for the central Delannoy number $D(n)$ (Remark~\ref{rmk:selberg}).

\textit{(c) Projectors and dimension.} Let $\lambda$ be a weight of length $n$. From (b), it follows that $L_{\lambda}$ has multiplicity one in $\cC(\bR^{(n)})$. We determine the projection operator $\cC(\bR^{(n)}) \to L_{\lambda}$ explicitly (Proposition~\ref{prop:projector}). As a corollary, we compute the categorical dimension of $L_{\lambda}$: it is $(-1)^n$ (Corollary~\ref{cor:cat-dim}). It is a remarkable feature of $\uRep(G)$ that every simple object has dimension $\pm 1$. We do not know of any other ``interesting'' examples of pre-Tannakian category with this property.

\textit{(d) Branching rules.} Let $G(0) \subset G$ be the stabilizer of $0 \in \bR$, which is an open subgroup of $G$. There are induction and restriction operations between $G$ and $G(0)$. A thematic problem in this situation is to determine the branching rules, i.e., the simple decompositions of the induction or restriction of simple objects. We solve this problem completely. A convenient feature of the present situation is that $G(0)$ is isomorphic to $G \times G$, and so we know exactly what its irreducible representations are.

To give a little more detail, we consider the induction
\begin{displaymath}
\Ind_{G \times G}^G(L_{\lambda} \boxtimes L_{\mu}),
\end{displaymath}
and completely describe its decomposition into simples (Theorem~\ref{thm:ind}); an interesting phenomenon here is that this induction always has length three. Similarly, we consider the restriction
\begin{displaymath}
\Res^G_{G \times G}(L_{\lambda}),
\end{displaymath}
and completely determine its decomposition into simples (Theorem~\ref{thm:res}). We note that one can deduce the branching rules between $G$ and any open subgroup from the $G(0)$ case. The computation of restriction plays an especially important role in the remainder of the paper; see Example~\ref{ex:res-dim} for a demonstration of how it can be a powerful tool.

\textit{(e) Tensor products.} Another thematic problem throughout representation theory is to describe the tensor product of irreducible representations. We solve this problem completely for $\uRep(G)$. Precisely, given weights $\lambda$ and $\mu$, we consider the tensor product
\begin{displaymath}
L_{\lambda} \uotimes L_{\mu}
\end{displaymath}
and give an explicit combinatorial rule for its simple decomposition (Theorem~\ref{thm:tensor}). Roughly speaking, the simples that appear are indexed by weights obtained by shuffling the words $\lambda$ and $\mu$; however, in these shuffles we allow letters to collide, and there are some rules for what happens in this case.

\textit{(f) The Grothendieck group.} Let $\rK$ be the Grothendieck group of the category $\uRep^{\rf}(G)$. We let $a_{\lambda} \in \rK$ be the class of the simple $L_{\lambda}$; these elements form a $\bZ$-basis for $\rK$. The tensor product on $\uRep(G)$ endows $\rK$ with the structure of a commutative ring. The restriction functor
\begin{displaymath}
\Res \colon \uRep(G) \to \uRep(G(0)) \cong \uRep(G \times G)
\end{displaymath}
defines a (non-co-commutative) co-multiplication on $\rK$. We show that this gives $\rK$ the structure of a Hopf algebra (see \S \ref{ss:hopf}). The Grothendieck group $\rK$ carries a natural ascending filtration, which is compatible with the Hopf algebra structure. We show that the associated graded is isomorphic (as a Hopf algebra) to the shuffle algebra (see \S \ref{ss:K-ring}). In particular, we find that $\rK \otimes \bQ$ is a polynomial ring on the $a_{\lambda}$ where $\lambda$ is a Lyndon word (Corollary~\ref{cor:K-poly}).

\textit{(f) Adams operations.} The Grothendieck group $\rK$ naturally carries the structure of a $\lambda$-ring, and thus admits Adams operations $\psi^n$. We prove that these operations are trivial, i.e., $\psi^n$ is the identity for all $n$ (Theorem~ \ref{thm:adams}); this implies that $\rK$ is a binomial ring (Corollary~\ref{cor:binom}). As a consequence, it is particularly easy to compute the action of Schur functors on $\uRep(G)$ (see \S \ref{ss:schur}). The triviality of all Adams operations is another remarkable feature of $\uRep(G)$: the only semi-simple pre-Tannakian categories we know with this property are $\uRep(G^n)$ with $n \in \bN$, and Czenky has recently shown that $\psi^2$ is non-trivial in any non-trivial symmetric fusion category of characteristic $\ne 2$ \cite[Theorem~4.8]{Czenky}.

\textit{(g) The path model.} As stated above, we show that $\uRep(G)$ is equivalent to a category defined using the combinatorics of Delannoy paths (Theorem~\ref{thm:equiv}).

\subsection{Connection to other work} \label{ss:other-work}

This paper relates to a few other topics.
\begin{itemize}
\item Representations of other oligomorphic groups, such as the infinite symmetric group, coincide with Deligne's interpolation categories. There is a substantial literature on them, such as \cite{ComesOstrik1, ComesOstrik, Deligne1, Deligne2, Deligne3, EntovaAizenbudHeidersdorf, Harman, Harman2, Knop, Knop2}.
\item Let $\OI$ be the category of finite totally ordered sets and monotonic injections. There is a natural functor $\OI^{\op} \to \uPerm(G)$ given on objects by $[n] \mapsto \cC(\bR^{(n)})$. It follows that $\OI$-modules, i.e., functors $\OI \to \Vect$, are closely related to $\uRep(G)$. These have been studied in \cite{GanLi, increp, unipotent, catgb}.
\item This paper studies a modified notion of representation for the group $G$. Ordinary representations of $G$ have been studied in \cite{DLLX, Nekrasov, Tsankov}.
\item The category $\uRep(G)$ is closely connected to the combinatorics of Delannoy paths. These have been studied extensively in the literature, e.g., \cite{Banderier, Covington, CDNS, OEIS, Sulanke1, Sulanke2, Tarnauceanu}.
\end{itemize}

\subsection{Outline}

We give an outline of the paper:
\begin{itemize}
\item In \S \ref{s:oligo}, we briefly review the general theory developed in \cite{HarmanSnowden}.
\item In \S \ref{s:AutR}, we specialize the discussion from \S \ref{s:oligo} to the case of $\Aut(\bR,<)$.
\item In \S \ref{s:simple}, we classify the simple objects of $\uRep(G)$, and determine the simple decomposition of the Schwartz space $\cC(\bR^{(n)})$.
\item In \S \ref{s:inv}, we analyze the invariant spaces of simple objects under open subgroups, determine the projection operator from $\cC(\bR^{(n)})$ to $L_{\lambda}$, and compute the categorical dimension of $L_{\lambda}$.
\item In \S \ref{s:ind}, we study the induction and restriction operations between $G$ and its open subgroup $G(0) \cong G \times G$. We also show that $\rK$ is a Hopf algebra.
\item In \S \ref{s:tensor}, we give a combinatorial rule for the tensor product in $\uRep(G)$ and show that $\gr(K)$ is isomorphic to the shuffle algebra.
\item In \S \ref{s:adams}, we prove that the Adams operations are trivial on the Grothendieck group.
\item In \S \ref{s:delannoy}, we discuss the path model for $\uRep(G)$.
\end{itemize}

\subsection{Notation}

We list some of the important notation here:
\begin{description}[align=right,labelwidth=2.5cm,leftmargin=!]
\item[ $k$ ] the coefficient ring (usually a field)
\item[ $\bbone$ ] unit object of a tensor category (e.g., the trivial representation)
\item[ $\bone$ ] the one-point set
\item[ $D(n,m)$ ] the $(n,m)$-Delannoy number (see \S \ref{ss:intro-path})
\item[ $\cC(X)$ ] the Schwartz space of $X$ (see \S \ref{ss:oligo-int})
\item[ $\uPerm(G)$ ] the category of permutation modules (see \S \ref{ss:oligo-perm})
\item[ $A(G)$ ] the completed group algebra of $G$ (see \S \ref{ss:oligo-rep})
\item[ $\uRep(G)$ ] the category of smooth $A(G)$-modules (see \S \ref{ss:oligo-rep})
\item[ $G$ ] the group $\Aut(\bR,<)$ (from \S \ref{s:AutR} onwards)
\item[ $G(A)$ ] the subgroup of $G$ fixing the set $A$ pointwise
\item[ $\bR^{(n)}$ ] the set of $x \in \bR^n$ with $x_1<\cdots<x_n$ (see \S \ref{ss:AutR})
\item[ $\vol(X)$ ] the volume of $X$ under the principal measure (see \S \ref{ss:AutR-meas})
\item[ $\ul{G}$ ] another notation for $A(G)$ (more or less; see \S \ref{ss:repcat})
\item[ $\rK$ ] the Grothendieck group of $\uRep^{\rf}(G)$ (see \S \ref{ss:groth})
\item[ $\res$ ] the restriction map on $\rK$ (see \S \ref{ss:groth})
\item[ $\ind$ ] the induction map on $\rK$ (see \S \ref{ss:groth})
\item[ $\lambda$ ] a weight, i.e., a word in the alphabet $\{ \wa, \wb \}$ (see \S \ref{ss:simples})
\item[ $L_{\lambda}$ ] the simple of weight $\lambda$ (see Definition~\ref{defn:simple})
\item[ $a_{\lambda}$ ] the class in $\rK$ of the simple $L_{\lambda}$ (see Definition~\ref{defn:simple})
\item[ $\odot$ ] the concatenation product on $\rK$ (see \S \ref{ss:concat})
\item[ {$\lambda[i,j]$} ] a substring of $\lambda$ (see \S \ref{ss:res})
\item[ $\pi(n)$ ] the length $n$ weight with all letters $\wa$ (see \S \ref{ss:simple-tensor})
\end{description}

\subsection*{Acknowledgments}

We thank Jordan Ellenberg, Victor Ostrik, and David Treumann for helpful discussions. We also thank Jeremy Miller, Peter Patzt, and Andrew Putman for organizing the conference ``Stability in Topology, Arithmetic, and Representation Theory'' at Purdue in March 2022, at which this project was conceived. NS was supported in part by NSF grant DMS-2000093. 

\section{Generalities on oligomorphic groups} \label{s:oligo}

In this section, we review the construction and basic properties of the category $\uRep(G)$ for a general oligomorphic group $G$. This discussion is a summary of the main points of \cite{HarmanSnowden}, except for \S \ref{ss:induction}, which is new.

\subsection{Oligomorphic groups}

An \defn{oligomorphic group} is a pair $(G, \Omega)$ consisting of a group $G$ and a set $\Omega$ equipped with a faithful action of $G$ such that $G$ has finitely many orbits on $\Omega^n$ for all $n$. For a finite subset $A$ of $\Omega$, let $G(A)$ be the subgroup of $G$ fixing each element of $A$. The $G(A)$'s form a neighborhood basis for a topology on $G$. This topology has the following three properties \cite[Proposition~2.4]{HarmanSnowden}:
\begin{enumerate}
\item It is Hausdorff.
\item It is non-archimedean: open subgroups form a neighborhood basis of the identity.
\item It is Roelcke precompact: if $U$ and $V$ are open subgroups then $U \backslash G/V$ is finite.
\end{enumerate}
An \defn{admissible group} is a topological group satisfying (a), (b), and (c). Thus every oligomorphic group comes with an admissible topology. Although we are most interested in oligomorphic groups, our constructions depend only on the topology and not the particular set $\Omega$. For this reason, it is most natural to work with admissible groups.

Let $G$ be an admissible group. We say that an action of $G$ on a set is \defn{smooth} if every stabilizer is open. We use the term \defn{$G$-set} for a set equipped with a smooth action of $G$. We say that a $G$-set is \defn{finitary} if it has finitely many orbits. See \cite[\S 2.3]{HarmanSnowden} for basic properties of $G$-sets.

We define a $\hat{G}$-set to be a set $X$ equipped with a smooth action of some open subgroup of $G$, which we call the \defn{group of definition} of $X$. Shrinking the group of definition does not change the $\hat{G}$-set. The notion of finitary is well-defined for $\hat{G}$-sets. The symbol $\hat{G}$ does not have any rigorous meaning on its own, but intuitively we think of $\hat{G}$ as an infinitesimal neighborhood of the identity in $G$. See \cite[\S 2.5]{HarmanSnowden} for additional details.

\begin{example}
Any finite group is admissible, under the discrete topoogy.
\end{example}

\begin{example}
Let $\fS$ be the infinite symmetric group, i.e., the group of all permutations of $\Omega=\{1,2,\ldots\}$. This is oligomorphic. Let $\fS(n)$ be the subgroup of $\fS$ fixing each of $1, \ldots, n$. The open subgroups of $\fS$ are those that contain some $\fS(n)$; in fact, every open subgroup is conjugate to one of the form $H \times \fS(n)$, where $H$ is a subgroup of the finite symmetric group $\fS_n$ \cite[Proposition~15.1]{HarmanSnowden}. The $\hat{\fS}$-subsets of $\Omega$ are the finite and cofinite subsets. More generally, the $\hat{\fS}$-subsets of $\Omega^n$ are those that can be defined by a first-order formula using using equality and finitely many constants.
\end{example}

\begin{example}
Many oligomorphic groups can be constructed using Fra\"iss\'e limits; for instance, the automorphism group of the Rado graph is oligomorphic. See \cite[\S 6.2]{HarmanSnowden} for a quick summary, and \cite{Cameron} for a detailed treatment.
\end{example}

\subsection{Measures and integration} \label{ss:oligo-int}

Fix an admissible group $G$ and a commutative ring $k$. We now come to a key concept introduced in \cite[\S 3.1]{HarmanSnowden}:

\begin{definition} \label{def:measure}
A \defn{measure} for $G$ valued in $k$ is a rule $\mu$ that assigns to every finitary $\hat{G}$-set a quantity $\mu(X)$ in $k$ such that the following axioms hold (in which $X$ and $Y$ denote finitary $\hat{G}$-sets):
\begin{enumerate}
\item Isomorphism invariance: $\mu(X)=\mu(Y)$ if $X \cong Y$.
\item Normalization: $\mu(\bone)=1$, where $\bone$ is the one-point $\hat{G}$-set.
\item Conjugation invariance: $\mu(X^g)=\mu(X)$, where $X^g$ is the conjugate of the $\hat{G}$-set $X$ by $g \in G$ (see \cite[\S 2.5]{HarmanSnowden} for the definition of $X^g$).
\item Additivity: $\mu(X \amalg Y)=\mu(X)+\mu(Y)$.
\item Multiplicativity in fibrations: if $X \to Y$ is a map of transitive $U$-sets, for some open subgroup $U$, with fiber $F$ (over some point) then $\mu(X)=\mu(F) \cdot \mu(Y)$. \qedhere
\end{enumerate}
\end{definition}

We note that a measure automatically satisfies $\mu(X \times Y)=\mu(X) \cdot \mu(Y)$. See \cite[\S 3.1]{HarmanSnowden} for additional details. There is a universal measure valued in a ring $\Theta(G)$. To construct $\Theta(G)$, start with the polynomial ring in symbols $[X]$, where $X$ is a finitary $G$-set (up to isomorphism), and impose relations corresponding to (b)--(e) above. See \cite[\S 4]{HarmanSnowden} for details. Computing $\Theta(G)$ is an important (and often difficult) problem.

Suppose now that we have a measure $\mu$ for $G$ valued in $k$. We then obtain a theory of integration, as follows. Let $X$ be a $\hat{G}$-set. We say that a function $\phi \colon X \to k$ is \defn{smooth} if it is invariant under some open subgroup of $G$ (contained in the group of definition for $X$) and we say $\phi$ has \defn{finitary support} if its support is a finitary $\hat{G}$-set. A \defn{Schwartz function} on $X$ is smooth $k$-valued functions on $X$ with finitary support, and the space of all Schwartz functions is called \defn{Schwartz space} and denoted $\cC(X)$. Given a Schwartz function $\phi$, let $c_1, \ldots, c_n$ be the distinct non-zero values it attains, and let $A_i=\phi^{-1}(c_i)$. We define the \defn{integral} of $\phi$ by
\begin{displaymath}
\int_X \phi(x) dx = \sum_{i=1}^n c_i \mu(A_i).
\end{displaymath}
Integration defines a $k$-linear map $\cC(X) \to k$. See \cite[\S 3.3]{HarmanSnowden} for more details.

If $f \colon X \to Y$ is a smooth map of $\hat{G}$-sets then there is an associated push-forward map $f_* \colon \cC(X) \to \cC(Y)$ given by integrating over fibers. This satisfies the expected properties, e.g., if $g \colon Y \to Z$ is a second map then $(gf)_*=g_*f_*$. If $X$ and $Y$ are themselves finitary, then there is also a pull-back map $f^* \colon \cC(Y) \to \cC(X)$. See \cite[\S 3.4]{HarmanSnowden} for more on push-forwards.

\begin{example}
A finite group $G$ admits a unique $\bZ$-valued measure via $\mu(X)=\# X$. This identifies $\Theta(G)$ with $\bZ$.
\end{example}

\begin{example} \label{ex:mut}
For each complex number $t$, there is a unique $\bC$-valued measure $\mu_t$ for $\fS$ with the property that $\mu_t(\Omega)=t$. This measure satisfies $\mu_t(\Omega^{(n)})=\binom{t}{n}$, where $\Omega^{(n)}$ denotes the set of $n$-element subsets of $\Omega$. A Schwartz function $\phi \colon \Omega \to \bC$ is simply one that is eventually constant. If $\phi(n)=a$ for all $n>N$ then
\begin{displaymath}
\int_{\Omega} \phi(x) dx = (t-N) a  + \sum_{i=1}^N \phi(i).
\end{displaymath}
The ring $\Theta(G)$ is identified with the ring of integer-valued polynomials; in particular, the $\mu_t$ account for all complex measures. See \cite[\S 15]{HarmanSnowden} for more details on this case.
\end{example}

\subsection{Permutation representations} \label{ss:oligo-perm}

Fix an admissible group $G$ and a $k$-valued measure $\mu$. If $X$ and $Y$ are finitary $\hat{G}$-sets, a \defn{$Y \times X$-matrix} $A$ is a Schwartz function $A \colon Y \times X \rightarrow k$. We write $\Mat_{Y,X}$ for the space of all $Y \times X$ matrices. We also write $\Mat_X$ for $\Mat_{X,X}$, and $\Vect_X$ for $\Mat_{X,\bone}$. If $A$ is a $Z \times Y$ matrix and $B$ is a $Y \times X$ matrix, we define the \defn{product matrix} $AB$ to be the $Z \times X$ matrix given by
\begin{displaymath}
AB(z,x) = \int_Y A(z,y) B(y,x) dy.
\end{displaymath}
Matrix multiplication gives $\Mat_X$ the structure of an associative and unital ring; the identity element is the \defn{identity matrix} $I_X$, which is simply the characteristic function of the diagonal. Much of ordinary linear algebra can be generalized to this setting; see \cite[\S 7]{HarmanSnowden} for details.

We now introduce the category of ``permutation modules'' $\uPerm(G)$; we refer to \cite[\S 8]{HarmanSnowden} for details. The category is defined as follows:
\begin{itemize}
\item The objects are $\Vect_X$, for $X$ a finitary $G$-set.
\item A morphism $\Vect_X \to \Vect_Y$ is a $G$-invariant $Y \times X$ matrix.
\item Composition is given by matrix multiplication.
\end{itemize}
The category $\uPerm(G)$ depends on the measure $\mu$, even though this is absent from the notation. The category $\uPerm(G)$ is naturally $k$-linear. It is also additive: we have
\begin{displaymath}
\Vect_X \oplus \Vect_Y = \Vect_{X \amalg Y}.
\end{displaymath}
The structure maps here are given by the usual projection matrices.

The category $\uPerm(G)$ carries a natural monoidal structure $\uotimes$. On objects, it is given by
\begin{displaymath}
\Vect_X \uotimes \Vect_Y = \Vect_{X \times Y}.
\end{displaymath}
On morphisms, it is given by the usual Kronecker product of matrices. The monoidal unit is  $\bbone = \Vect_\bone$. This monoidal structure is naturally symmetric; the symmetric structure is induced by the symmetric structure on the category of $G$-sets.

The tensor category $\uPerm(G)$ is rigid. Every object is self-dual. The evaluation and co-evaluation morphisms for $\Vect_X$ are both given by the characteristic function of the diagonal. The categorical trace of a morphism $A \in \Mat_{X,X}^G$ is the ``usual'' matrix trace
\begin{displaymath}
\int_X A(x,x) dx.
\end{displaymath}
In particular, the categorical dimension of $\Vect_X$ is $\mu(X)$. (See \cite[\S 4.7]{Etingof} for the definition of categorical trace and dimension.) 

\begin{example}
If $G$ is a finite group then $\uPerm(G)$ is equivalent the usual category of permutation $G$-modules.
\end{example}

\subsection{General representations} \label{ss:oligo-rep}

Let $G$ and $\mu$ be as above. The category $\uPerm(G)$ is typically not abelian. We now explain how we can produce an abelian representation category $\uRep(G)$, which is well-behaved in certain circumstances.

The \defn{completed group algebra} of $G$, denoted $A(G)$, is the inverse limit of the Schwartz spaces $\cC(G/U)$ over open subgroups $U$. Convolution endows $A(G)$ with the structure of an associative and unital algebra; see \cite[\S 10.3]{HarmanSnowden}. The algebra structure on $A(G)$ depends on the measure $\mu$, even though it is absent from the notation. An $A(G)$-module $M$ is called \defn{smooth} if for every element $x$ there is an open subgroup $U$ such that the action of $A(G)$ on $x$ factors through $\cC(G/U)$. We define $\uRep(G)$ to be the full subcategory of $\Mod_{A(G)}$ spanned by smooth modules. The category $\uRep(G)$ is always a Grothendieck abelian category.

In order for $\uRep(G)$ to be well-behaved, we must make some assumptions on the measure $\mu$. We say that $\mu$ \defn{regular} if $\mu(X)$ is a unit in $k$ for every transitive $G$-set $X$, and \defn{quasi-regular} if there is an open subgroup $U$ of $G$ such that $\mu \vert_U$ is regular. The measure $\mu_t$ in Example~\ref{ex:mut} is quasi-regular for all $t$, and regular if and only if $t \not\in \bN$. There is another condition on measures used in \cite{HarmanSnowden} called \defn{normal}, which is weaker than quasi-regular, but we will not need it here. We assume in the following discussion that $\mu$ is quasi-regular.

The most important objects of $\uRep(G)$ are the Schwartz spaces $\cC(X)$, where $X$ is a $G$-set. The $A(G)$-module structure on $\cC(X)$ is defined via convolution; see \cite[Proposition~10.11]{HarmanSnowden}. The Schwartz spaces $\cC(X)$ have a mapping property \cite[Corollary~11.19]{HarmanSnowden}, and every object of $\uRep(G)$ is a quotient of some $\cC(X)$ \cite[Proposition~11.14]{HarmanSnowden}. There is a functor
\begin{displaymath}
\Phi \colon \uPerm(G) \to \uRep(G), \qquad \Vect_X \mapsto \cC(X)
\end{displaymath}
that is fully faithful \cite[Proposition~11.12]{HarmanSnowden}. Note that $\Vect_X$ and $\cC(X)$ are the same vector space; we use different notation to emphasize that $\cC(X)$ is thought of in $\uRep(G)$.

Although we will not precisely define the completed group algebra or its action on $\cC(X)$, we can give a direct description of submodules of $\cC(X)$. Let $\phi$ belong to $\cC(X)$, and suppose that $\phi$ is invariant under the open subgroup $U \subset G$. Given $\psi \in \cC(G/U)$, we define $\psi \ast \phi$ to be the function on $X$ given by
\begin{displaymath}
(\psi \ast \phi)(x) = \int_{G/U} \psi(g) \phi(g^{-1} x) dg.
\end{displaymath}
We then have the following criterion: a $k$-subspace $V$ of $\cC(X)$ is an $A(G)$-submodule if and only if it is closed under the above convolutions, i.e., given $\phi \in V$ that is $U$-invariant and $\psi \in \cC(G/U)$ the element $\psi \ast \phi$ belongs to $V$. This also describes how to form the $A(G)$-submodule of $\cC(X)$ generated by a set of elements (take all such convolutions).

Since the Schwartz spaces generate $\uRep(G)$, there is at most one symmetric monoidal structure on $\uRep(G)$ that is compatible with $\Phi$ and co-continuous in each variable. Such a symmetric monoidal structure is constructed in \cite[Theorem~12.9]{HarmanSnowden}, and we denote it by $\uotimes$. It satisfies $\cC(X) \uotimes \cC(Y)=\cC(X \times Y)$ for arbitrary $G$-sets $X$ and $Y$.

In fact, $\uotimes$ satisfies a mapping property: giving a map from $M \uotimes N$ is equivalent to giving a $G$-equivariant bilinear map from $M \times N$ that commutes with certain integrals; see \cite[\S 12.2]{HarmanSnowden}. This permits us to give a more concrete description of the self-duality of $\cC(X)$: the pairing $\cC(X) \uotimes \cC(X) \to \bbone$ is induced from the bilinear form
\begin{displaymath}
\langle, \rangle \colon \cC(X) \times \cC(X) \to k, \qquad
\langle \phi, \psi \rangle = \int_X \phi(x) \psi(x) dx.
\end{displaymath}
This pairing will be used several times throughout this paper.

\begin{remark}
We warn the reader that the forgetful functor from $\uRep(G)$ to $k$-modules is not monoidal. For example, $\cC(X) \uotimes \cC(X) = \cC(X \times X)$, which is typically larger than $\cC(X) \otimes_k \cC(X)$ since one usually cannot express the characteristic function of the diagonal as a finite sum of pure tensors. One should think of $\cC(X) \uotimes \cC(Y)$ as a kind of completion of the algebraic tensor product.
\end{remark}

We recall the following general definition:

\begin{definition} \label{defn:pre-tannakian}
Suppose $k$ is a field. We say that a $k$-linear symmetric monoidal category  $\cC$ is \defn{pre-Tannakian} if it is abelian, all objects have finite length, all Hom spaces are finite dimensional, all objects have duals, and $\End_{\cC}(\bbone)=k$, where $\bbone$ is the unit object.
\end{definition}

The following is the most important general theorem on the structure of $\uRep(G)$. It uses the somewhat technical condition \defn{Property~(P)} on the measure $\mu$. This property essentially means that $\mu$ is valued in a subring of $k$ that maps to enough fields of positive characteristic. Since the property holds in all cases of interest to this paper, we do not bother to explain it in detail; see \cite[Definition~7.17]{HarmanSnowden} for details.

\begin{theorem}[{\cite[Theorem~13.2]{HarmanSnowden}}]
Assume that $k$ is a field, $\mu$ is quasi-regular, and Property~(P) holds.
\begin{enumerate}
\item Every object of $\uRep(G)$ is the union of its finite length subobjects.
\item The category $\uRep^{\rf}(G)$ of finite length objects is pre-Tannakian.
\item If $\mu$ is regular then $\uRep^{\rf}(G)$ is semi-simple.
\end{enumerate}
\end{theorem}

\begin{remark}
One can view part (c) above as an analogue of Maschke's theorem.
\end{remark}

As a corollary of the theorem, we find that in the regular case one can abstractly recover the category $\uRep^{\rf}(G)$ without going through the completed group algebra:

\begin{corollary}
If $\mu$ is regular and $k$ is a field then $\uRep^{\rf}(G)$ is the Karoubian envelope of $\uPerm(G)$.
\end{corollary}

\subsection{Induction and restriction} \label{ss:induction}

Let $G$ be an admissible group equipped with a quasi-regular measure $\mu$, and let $U$ be an open subgroup of $G$. The measure $\mu$ restricts to a quasi-regular measure on $U$ \cite[Proposition~3.23]{HarmanSnowden}, and there is a restriction functor \cite[\S 10.5]{HarmanSnowden}
\begin{displaymath}
\Res^G_U \colon \uRep(G) \to \uRep(U).
\end{displaymath}
This functor is simply restriction of scalars along the natural homomorphism $A(U) \to A(G)$.

We now discuss induction (which was not treated in \cite{HarmanSnowden}). Let $N$ be an object of $\uRep(U)$. We define $\Ind_U^G(N)$ to be the space of all functions $\phi \colon G \to N$ satisfying the following two conditions:
\begin{itemize}
\item $\phi$ is left $U$-equivariant, i.e., $\phi(ug)=u\phi(g)$ for all $u \in U$ and $g \in G$
\item $\phi$ is right $G$-smooth, i.e., there is an open subgroup $V$ of $G$ such that $\phi(gv)=\phi(g)$ for all $v \in V$ and $g \in G$.
\end{itemize}
Given $a \in A(G)$ and $\phi \in \Ind_U^G(N)$, we define a function $a\phi \colon G \to N$ by
\begin{displaymath}
(a\phi)(g) = \int_{G/V} a_V(h) \phi(gh) dh,
\end{displaymath}
where $V$ is an open subgroup of $G$ for which $\phi$ is right invariant. (Here $a_V$ is the component of $a$ in $\cC(G/V)$; recall that $A(G)$ is the inverse limit of these spaces.) The function $a \phi$ is well-defined (i.e., independent of the choice of $V$), and belongs to $\Ind_U^G(N)$; furthermore, this construction defines on $\Ind_U^G(N)$ the structure of a smooth $A(G)$-module. We omit the verification of these statements. We call $\Ind_U^G(N)$ the \defn{induction} of $N$ from $U$ to $G$. One easily sees that induction defines an additive functor
\begin{displaymath}
\Ind_U^G \colon \uRep(U) \to \uRep(G).
\end{displaymath}
The following proposition gives its main properties.

\begin{proposition} \label{prop:induction}
We have the following:
\begin{enumerate}
\item The induction functor is both left and right adjoint to the restriction functor.
\item The induction functor is continuous and co-continuous (and thus exact).
\item For a $U$-set $X$, we have $\Ind_U^G(\cC(X))=\cC(I_U^G(X))$, where $I_U^G(X)=G \times^U X$ is the induction of the $U$-set $X$ to $G$.
\end{enumerate}
\end{proposition}

\begin{proof}
(a) Let $M$ and $N$ be objects of $\uRep(G)$ and $\uRep(U)$ respectively. To show that induction is right adjoint to restriction, we must give a natural isomorphism
\begin{displaymath}
\Hom(M, \Ind_U^G(N)) = \Hom(\Res^G_U(M), N).
\end{displaymath}
Suppose $\alpha \colon M \to \Ind_U^G(N)$ is given. We define $\alpha' \colon \Res^G_U(M) \to N$ by $\alpha'(m)=\alpha(m)(1)$, where the right side indicates the value of the function $\alpha(m)$ at $1 \in G$. Now suppose that $\beta \colon \Res^G_U(M) \to N$ is given. We define $\beta' \colon M \to \Ind_U^G(N)$ by $\beta'(m)(g)=\beta(gm)$. We leave to the reader the verification that these constructions are well-defined and mutually inverse.

To show that induction is left adjoint to restriction, we must give a natural isomorphism
\begin{displaymath}
\Hom(\Ind_U^G(N), M) = \Hom(N, \Res^G_U(M)).
\end{displaymath}
Suppose $\alpha \colon \Ind_U^G(N) \to M$ is given. We define $\alpha' \colon N \to \Res^G_U(M)$ by $\alpha'(n)=\alpha(\gamma_n)$, where $\gamma_n \colon G \to N$ is the function defined by
\begin{displaymath}
\gamma_n(g) = \begin{cases}
gn & \text{if $g \in U$} \\
0 & \text{otherwise} \end{cases}
\end{displaymath}
Now suppose that $\beta \colon N \to \Res^G_U(M)$ is given. We define $\beta' \colon \Ind^G_U(N) \to M$ by $\beta'(\phi) = \int_{G/U} g \phi(g^{-1}) dg$, where here we are integrating a module-valued function (see \cite[\S 11.9]{HarmanSnowden}). We again leave the necessary verifications to the reader.

(b) This follows from (a) and general properties of adjoints.

(c) Let $M$ be an object of $\uRep(G)$. We have
\begin{align*}
\Hom_{A(G)}(\Ind_U^G(\cC(X)), M)
&= \Hom_{A(U)}(\cC(X), M)
= \Hom_U(X, M) \\
&= \Hom_G(I_U^G(X), M)
= \Hom_{A(G)}(\cC(I_U^G(X)), M).
\end{align*}
In the first step, we used (a); in the second, the mapping property for Schwartz space; in the third, adjunction for $I^G_U$ (at the $G$-set level); and in the fourth, the property for Schwartz space again. The result now follows from Yoneda's lemma.
\end{proof}

There is one more general fact about induction we require, namely, the projection formula:

\begin{proposition} \label{prop:projection}
Let $M$ and $N$ be objects of $\uRep(G)$ and $\uRep(U)$ respectively. Then there is a natural isomorphism
\begin{displaymath}
\Ind_U^G(\Res^G_U(M) \uotimes N) \cong M \uotimes \Ind_U^G(N)
\end{displaymath}
\end{proposition}

\begin{proof}
This can be deduced from general considerations (see, e.g., \cite[Proposition~4.6]{Harriger}), but we will give a direct argument. There is a natural map (from the left side to the right side) coming from adjunctions; we must show that it is an isomorphism. Since all functors involved are right exact, it suffices (by choosing presentations) to verify this when $M=\cC(X)$ and $N=\cC(Y)$, where $X$ and $Y$ are $G$- and $U$-sets. This follows from the projection formula at the set level, i.e., we have a natural isomorphism of $G$-sets
\begin{displaymath}
I_U^G(R^G_U(X) \times Y)=X \times I_U^G(Y),
\end{displaymath}
where $R^G_U$ is the restriction functor from $G$-sets to $U$-sets.
\end{proof}

\section{Automorphisms of the line} \label{s:AutR}

In this section, we introduce the main group of interest: the automorphism group of the line. We recall some results about this case from \cite{HarmanSnowden}, and develop some additional basic theory.

\subsection{The group} \label{ss:AutR}

Let $G=\Aut(\bR,<)$ be the group of all order-preserving bijections $\bR \to \bR$. The symbol $G$ will denote this group for the remainder of the paper. It is easy to see that $G$ is oligomorphic with respect to its action on $\bR$.

Let $\bR^{(n)}$ denote the set of $n$-element subsets of $\bR$. We identify $\bR$ with the subset of $\bR^n$ consisting of tuples $(x_1, \ldots, x_n)$ satisfying $x_1<x_2<\cdots<x_n$. One easily sees that $G$ acts transitively on $\bR^{(n)}$. For $a \in \bR^{(n)}$, we let $G(a)$ denote the subgroup of $G$ consisting of elements that fix each $a_i$ (which is equivalent to fixing the set $a$). Recall that a subgroup of $G$ is open if it contains some $G(a)$. In fact:

\begin{proposition}
Every open subgroup of $G$ is of the form $G(a)$ for some $a \in \bR^{(n)}$ and $n \in \bN$.
\end{proposition}

\begin{proof}
See \cite[Proposition~17.1]{HarmanSnowden}.
\end{proof}

\begin{corollary} \label{cor:open-class}
Every transitive $G$-set\footnote{Recall our convention that ``$G$-set'' means ``set with a \emph{smooth} $G$-action.''} is isomorphic to $\bR^{(n)}$, for some $n$.
\end{corollary}

For an open interval $I \subset \bR$, let $G_I$ be the group of order-preserving bijections $I \to I$. Choosing an order-preserving bijection $i \colon I \to \bR$ yields an isomorphism $f \colon G_I \to G$ of topological groups. Different choices of $i$ lead to different choices of $f$; however, the various $f$'s only differ by an inner automorphisms of $G$. Thus we have a well-defined identification $G_I \cong G$ up to inner automorphisms.

For $a \in \bR^{(n)}$, the set $\bR \setminus a$ is a disjoint union of $n+1$ open intervals $I_0, \ldots, I_n$. The group $G(a)$ preserves these intervals, and so there is an induced homomorphism
\begin{displaymath}
G(a) \to G_{I_0} \times \cdots \times G_{I_n}.
\end{displaymath}
One easily sees that this map is an isomorphism of topological groups. Thus $G(a)$ is isomorphic to $G^{n+1}$. This self-similarity present in the structure of $G$ will play an important role in its representation theory.

Since $G_I$ is isomorphic to $G$, Corollary~\ref{cor:open-class} implies that every transitive $G_I$-set is isomorphic to $I^{(n)}$ for some $n$. Combined with the above classification and description of open subgroups of $G$, we obtain the following description of $\hat{G}$ sets:

\begin{proposition} \label{prop:Ghat}
Every finitary $\hat{G}$-set is a finite union of sets of the form $I_1^{(n_1)} \times \cdots \times I_r^{(n_r)}$ where the $I_i$ are open intervals in $\bR$ and the $n_i$ are natural numbers.
\end{proposition}

\subsection{Measure and integration} \label{ss:AutR-meas}

Let $k$ be a commutative ring. By \cite[Theorem~17.7]{HarmanSnowden}, the group $G$ admits a unique $k$-valued measure $\mu$ satisfying
\begin{displaymath}
\mu(I_1^{(n_1)} \times \cdots \times I_r^{(n_r)}) = (-1)^{n_1+\cdots+n_r},
\end{displaymath}
where the $I_i$'s are disjoint open intervals of $\bR$ and the $n_i$'s are natural numbers. For a general $\hat{G}$-set $X$, we have $\mu(X)=\chi_c(X)$, where $\chi_c$ is compactly supported Euler characteristic; one can define $\chi_c(X)$ by appealing to Proposition~\ref{prop:Ghat}. Since $\mu(\bR^{(n)})=(-1)^n$, it follows that the measure $\mu$ is regular. We refer to $\mu$ as the \defn{principal measure} on $G$, and it is the only measure for $G$ that we will use. In what follows, we simply write $\vol(X)$ in place of $\mu(X)$, and we no longer use the symbol $\mu$ for a measure.

Integration with respect to the principal measure, as defined in \S \ref{ss:oligo-int}, is essentially the Euler calculus of Schapira and Viro (see, e.g., \cite{Viro}). Here are some important examples of integrals:
\begin{displaymath}
\int_{\bR} \delta_p(x) dx = 1, \qquad
\int_{\bR} 1_{(p,q)}(x) dx = -1, \qquad
\int_{\bR} 1_{(p,q]}(x) dx = 0.
\end{displaymath}
Here $p<q$ are real numbers, $\delta_p$ is the point mass at $p$, and $1_X$ is the characteristic function of the set $X$.

\begin{remark}
The group $G$ admits essentially four measures in total; more precisely, the universal ring $\Theta(G)$ (discussed after Definition~\ref{def:measure}) is isomorphic to $\bZ^4$ \cite[Theorem~17.7]{HarmanSnowden}. The three non-principal measures are not quasi-regular though.
\end{remark}

\subsection{Representation categories} \label{ss:repcat}

We now come to the main object of study of this paper: the category $\uRep(G)$, taken with respect to the principal measure on $G$ over the field $k$. One can consider this category over a general commutative ring, but we confine our attention to the field case. As discussed in \S \ref{ss:oligo-rep}, $\uRep(G)$ is a semi-simple Grothendieck abelian category equipped with a tensor product $\uotimes$, and the category $\uRep^{\rf}(G)$ of finite length objects is pre-Tannakian.

Associated to each $G$-set $X$, we have the Schwartz space $\cC(X)$ in $\uRep(G)$. By the classification of transitive $G$-sets (Corollary~\ref{cor:open-class}), $\cC(X)$ decomposes into a direct sum of $\cC(\bR^{(n)})$'s. Thus every simple object of $\uRep(G)$ is a summand of some $\cC(\bR^{(n)})$. The Schwartz spaces $\cC(\bR^{(n)})$ will play a central role in our analysis of $\uRep(G)$.

We introduce one additional piece of notation/terminology: we use the term ``$\ul{G}$-module'' or ``representation of $\ul{G}$'' (or similar) for objects of the category $\uRep(G)$, and ``$\ul{G}$-map'' (or similar) for morphisms in $\uRep(G)$. Thus $\ul{G}$ is essentially synonymous with the completed group algebra $A(G)$, with the convention that ``$\ul{G}$-module'' means ``smooth module over $A(G)$.'' This terminology makes the language closer to classsical representation theory, which we find helpful.

\subsection{Orbits on products}

We now arrive at the connection between $\uRep(G)$ and Delannoy paths. For the moment, we establish a simple bijection; in \S \ref{s:delannoy}, we will find deeper connections.

\begin{proposition} \label{prop:orbit}
The $G$-orbits on $\bR^{(n)} \times \bR^{(m)}$ are naturally in bijective correspondence with $(n,m)$-Delannoy paths. In particular, the number of orbits is the Delannoy number $D(n,m)$.
\end{proposition}

\begin{figure}
\begin{center}
\begin{tikzpicture}[scale=0.5,baseline=0]
\draw[color=gray!50] (0,2)--(14,2);
\draw[red,fill=red] (1,2) circle (1.1ex);
\draw[red,fill=red] (3,2) circle (1.1ex);
\draw[fill=yellow] (4.8,1.8) rectangle (5.2,2.2);
\draw[blue,line width=2pt] (7,2) circle (1ex);
\draw[red,fill=red] (9,2) circle (1ex);
\draw[fill=yellow] (10.8,1.8) rectangle (11.2,2.2);
\draw[blue,line width=2pt] (13,2) circle (1.1ex);
\node at (1,1.35) {\tiny 1};
\node at (3,1.35) {\tiny 2};
\node at (5,1.35) {\tiny 3};
\node at (7,1.35) {\tiny 4};
\node at (9,1.35) {\tiny 5};
\node at (11,1.35) {\tiny 6};
\node at (13,1.35) {\tiny 7};
\end{tikzpicture}
\qquad\qquad
\begin{tikzpicture}[scale=0.5]
\draw[step=1, color=gray!50] (0, 0) grid (5,4);
\draw[line width=2pt] (0,0)--(1,0)--(2,0)--(3,1)--(3,2)--(4,2)--(5,3)--(5,4);
\end{tikzpicture}
\end{center}
\caption{An illustration of the proof of Proposition~\ref{prop:orbit}, with $x=(1,2,3,5,6)$ and $y=(3,4,6,7)$. The solid red dots are placed at lone $x_i$'s, the open blue circles at lone $y_j$'s, and yellow squares at points that are both of the form $x_i$ and $y_j$.} \label{fig:orbit}
\end{figure}
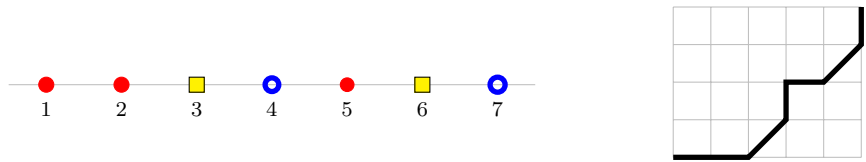

\begin{proof}
Let $(x,y)$ be an element of $\bR^{(n)} \times \bR^{(m)}$. We associate to $(x,y)$ an $(n,m)$-Delannoy path $p(x,y)$ as follows. Put a red dot on the real line at each $x_i$ and a blue dot at each $y_j$; a red and blue dot at the same point make a yellow dot. Now start in the plane at $(0,0)$ and read the real line from $-\infty$ to $\infty$. Whenever a red dot is encountered, take a $(1,0)$ step; for blue dots, take a $(0,1)$ step; and for yellow, a $(1,1)$ step. The resulting path is $p(x,y)$. See Figure~\ref{fig:orbit} for an example.

Let $(x',y')$ be a second element of $\bR^{(n)} \times \bR^{(m)}$. It is clear that if $(x,y)$ and $(x',y')$ belong to the same $G$-orbit then $p(x,y)=p(x',y')$. Conversely, suppose that $p(x,y)=p(x',y')$. Regarding $x$ and $y$ as subsets of $\bR$ size $n$ and $m$, let $z$ be their union; similarly, let $z'$ be the union of $x'$ and $y'$. The cardinality of $z$ is the length (=number of steps) of the path $p(x,y)$, and similarly for $z'$. We thus see that $z$ and $z'$ have the same cardinality, say $\ell$. Since $G$ acts transitively on $\bR^{(\ell)}$, there is $g \in G$ such that $gz=z'$, i.e., $gz_r=z'_r$ for all $1 \le r \le \ell$. Now, whether $z_r$ is an $x_i$, $y_j$, or both can be determined from $p(x,y)$. Since $p(x,y)=p(x',y')$, it follows that $gx=x'$ and $gy=y'$, and so $(x,y)$ and $(x',y')$ are in the same orbit.
\end{proof}

\begin{corollary} \label{cor:hom-dim}
We have
\begin{displaymath}
\dim \Hom_{\ul{G}}(\cC(\bR^{(n)}), \cC(\bR^{(m)})) = D(n,m).
\end{displaymath}
\end{corollary}

\begin{proof}
Recall that $\ul{G}$-maps $\cC(\bR^{(n)}) \to \cC(\bR^{(m)})$ are given by $G$-invariant $\bR^{(m)} \times \bR^{(n)}$ matrices. Thus the space of maps is isomorphic to the space of $k$-valued functions on the orbit space $G \backslash (\bR^{(m)} \times \bR^{(n)})$, and so the result follows.
\end{proof}

\begin{remark}
The $G$-orbits on $\bR^{(n_1)} \times \cdots \times \bR^{(n_d)}$ can similarly be parametrized by $d$-dimensional Delannoy paths, as defined in \S \ref{ss:path}.
\end{remark}

\subsection{The Grothendieck group} \label{ss:groth}

Let $\rK$ be the Grothendieck group of the category $\uRep^{\rf}(G)$. For a finite length $\ul{G}$-representation $V$, we let $[V]$ denote its class in $\rK$. Classes of the form $[V]$ are called \defn{effective}. The classes $[V]$, with $V$ a simple representation, form a basis of $\rK$; we will explicitly describe this basis in Corollary~\ref{cor:K-basis} below. We now explain a number of features of $\rK$.

\textit{(a) The unit class.} We let $1 \in \rK$ be the class of the trivial representation.

\textit{(b) The ring structure.} The tensor product $\uotimes$ on $\uRep(G)$ induces a product on $\rK$, which we call the \defn{standard product}. (We use this term to distinguish it from two other products discussed below: the induction and concatenation products.) Under the standard product, $\rK$ is an associative, unital, and commutative ring; the unit element is~1.

\textit{(c) The Hom pairing.} We let
\begin{displaymath}
\langle, \rangle \colon \rK \times \rK \to \bZ
\end{displaymath}
be the Hom pairing, defined on effective classes by
\begin{displaymath}
\langle [V], [W] \rangle = \dim_k \Hom_{\ul{G}}(V,W).
\end{displaymath}
It is symmetric since $\uRep(G)$ is semi-simple. The classes of simple objects form an orthogonal basis; in fact, it follows from Theorem~\ref{thm:simple} that the simples are absolutely simple and so the classes of simples form an orthonormal basis.

\textit{(d) Induction and restriction.} Let $G(0) \subset G$ be the stabilizer of the point $0 \in \bR$. By definition, $G(0)$ is an open subgroup of $G$. We thus have functors
\begin{displaymath}
\Ind \colon \uRep(G(0)) \to \uRep(G), \qquad
\Res \colon \uRep(G) \to \uRep(G(0))
\end{displaymath}
as in \S \ref{ss:induction}. The group $G(0)$ is isomorphic to $G \times G$, and so the Grothendieck group of $\uRep^{\rf}(G(0))$ is isomorphic to $\rK \otimes \rK$. The above functors therefore induce additive maps
\begin{displaymath}
\ind \colon \rK \otimes \rK \to \rK, \qquad
\res \colon \rK \to \rK \otimes \rK.
\end{displaymath}
These maps are associative and co-associative, respectively (see below). The product $\ind$ is not unital or commutative. We will see that $\res$ is co-unital, but it is not co-commutative. We note that $\res$ is a ring homomorphism (with respect to the standard product), since $\Res$ is a tensor functor.

[We briefly sketch a proof of associativity. Let $G(0,1) \subset G$ be the stabilizer of $\{0,1\} \in \bR^{(2)}$, which is isomorphic to $G \times G \times G$. Induction from $G(0,1)$ to $G$ defines a map $\rK \otimes \rK \otimes \rK \to \rK$. One can factor this induction in two ways: first induce from $G(0,1)$ to $G(0)$, and then to $G$; or first induce to $G(1)$, and then to $G$. These equivalent factorizations prove that $\ind$ is associative. We note that associativity also follows from our explicit computation of $\ind$ in Theorem~\ref{thm:ind}.]

\textit{(e) Frobenius reciprocity.} For $x,y,z \in \rK$, we have
\begin{displaymath}
\langle \ind(x \otimes y), z \rangle = \langle x \otimes y, \res(z) \rangle.
\end{displaymath}
On the right side, $\langle, \rangle$ is the Hom pairing on $\rK \otimes \rK$, which is simply the tensor product of the Hom pairing on $\rK$ with itself. This follows from Frobenius reciprocity at the representation level, as discussed in \S \ref{ss:induction}.

\textit{(f) The projection formula.} For $x,y,z \in \rK$, we have
\begin{displaymath}
\ind(\res(x) \cdot (y \otimes z))=x \cdot \ind(y \otimes z).
\end{displaymath}
This follows from Proposition~\ref{prop:projection}.

\textit{(g) Mackey theory.} For $x,y \in \rK$, we have
\begin{displaymath}
\res(\ind(x \otimes y)) = x \otimes y + \ind_{1,2}(x \otimes \res(y)) + \ind_{2,3}(\res(x) \otimes y),
\end{displaymath}
where $\ind_{1,2}$ denotes the induction product on the first two tensor factors. Since we will not use this formula, we omit the proof.

\textit{(h) The filtration.} We let $\rK_{\le n}$ be the $\bZ$-submodule of $\rK$ spanned by the simple classes appearing in $\cC(\bR^{(m)})$ for some $0 \le m \le n$; see Remark~\ref{rmk:K-filt-basis} for an explicit basis of $\rK_{\le n}$. The $\rK_{\le n}$ define an increasing exhaustive filtration of $\rK$. This filtration is compatible with multiplication in the sense that
\begin{displaymath}
\rK_{\le n} \cdot \rK_{\le m} \subset \rK_{\le n+m}.
\end{displaymath}
Indeed, if $N$ and $M$ are simple $\ul{G}$-submodules of $\cC(\bR^{(n)})$ and $\cC(\bR^{(m)})$ then $N \uotimes M$ is a $\ul{G}$-submodule of
\begin{displaymath}
\cC(\bR^{(n)}) \uotimes \cC(\bR^{(m)}) = \cC(\bR^{(n)} \times \bR^{(m)}),
\end{displaymath}
and $\bR^{(n)} \times \bR^{(m)}$ decomposes into orbits of the form $\bR^{(s)}$ with $s \le n+m$. The filtration is also compatible with the co-multiplication $\res$, in the sense that
\begin{displaymath}
\res(\rK_{\le n}) \subset \sum_{i+j=n} \rK_{\le i} \otimes \rK_{\le j}.
\end{displaymath}
One can prove this directly by examining the restriction of $\cC(\bR^{(n)})$, or deduce it from our explicit computation of $\res$ in Theorem~\ref{thm:res}.

\textit{(i) Other structure.} There is some other structure on $\rK$ that we are not yet ready to define in detail, but we mention here. In \S \ref{ss:concat}, we define a third product on $\rK$, called the concatentation product, and denoted $\odot$. We show that $\res$ admits a co-unit $\delta$ (Propositions~\ref{prop:delta} and~\ref{prop:counit}). We also show that $\rK$ forms a Hopf algebra under the standard product and the restriction co-product (Proposition~\ref{prop:hopf}). Finally, we show that $\rK$ has the structure of a $\lambda$-ring (see \S \ref{ss:lambda-bg}); in fact, it is a binomial ring (see Corollary~\ref{cor:binom}).

\section{Simple representations} \label{s:simple}

In this section, we construct and classify the simple objects of $\uRep(G)$. We also determine the simple decomposition of the Schwartz space $\cC(\bR^{(n)})$.

\subsection{Construction of simples} \label{ss:simples}

To begin, we construct certain simples, which will turn out to be all of them. To parametrize these objects, we introduce the following notion:

\begin{definition} \label{defn:weight}
A \emph{weight} is a word $\lambda$ in the alphabet $\{\wa,\wb\}$. We write $\ell(\lambda)$ for the length $\lambda$, and let $\Lambda$ denote the set of all weights.
\end{definition}

By a half-open interval in $\bR$, we mean a non-empty interval of the form $(b,a]$ or $[a,b)$, where $a \in \bR$ and $b \in \bR \cup \{\pm \infty\}$. We define the \defn{type} of a half-open interval to be $\wa$ if its right endpoint is included, and $\wb$ if its left endpoint is included. For two subsets $I$ and $J$ of $\bR$, we write $I<J$ to mean that every element of $I$ is less than every element of $J$. For example, we have $[a,b)<[b,c)$ if $a<b<c$.

Now consider a tuple $\bI=(I_1,\ldots,I_n)$ of half-open intervals. We define the \defn{type} of $\bI$ to be the weight $\lambda_1 \cdots \lambda_n$ where $\lambda_i$ is the type of $I_i$. We write $\bI$ still for the hypercube $I_1 \times \cdots I_n$ in $\bR^n$, and we let $\phi_{\bI} \in \cC(\bR^n)$ be its characteristic function. We say that $\bI$ is \defn{ordered} if $I_1<\cdots<I_n$. In this case, $\bI$ is contained in $\bR^{(n)}$, and we regard $\phi_{\bI}$ as an element of $\cC(\bR^{(n)})$. We can now introduce some extremely important representations:

\begin{definition} \label{defn:simple}
For a weight $\lambda$ of length $n$, we define $L_{\lambda}$ to be the $\ul{G}$-submodule of $\cC(\bR^{(n)})$ generated by the $\phi_{\bI}$, where $\bI$ is an ordered tuple of half-open intervals of type $\lambda$. We also let $a_{\lambda}$ be the class of $L_{\lambda}$ in the Grothendieck group $\rK$.
\end{definition}

The following is the main theorem of \S \ref{ss:simples}:

\begin{theorem} \label{thm:simple}
We have the following:
\begin{enumerate}
\item The module $L_{\lambda}$ is simple, and $\End_{\ul{G}}(L_{\lambda})=k$.
\item If $\lambda$ and $\mu$ are distinct weights then $L_{\lambda}$ and $L_{\mu}$ are non-isomorphic.
\end{enumerate}
\end{theorem}

\begin{figure}
\begin{center}
\begin{tikzpicture}[scale=0.5]
\draw[step=1, color=gray!50] (0, 0) grid (4,4);
\draw[line width=2pt] (0,0)--(1,0)--(1,1)--(2,2)--(2,3)--(3,3)--(4,4);
\end{tikzpicture}
\qquad\qquad\qquad
\begin{tikzpicture}[scale=0.5]
\draw[step=1, color=gray!50] (0, 0) grid (4,4);
\draw[line width=2pt] (0,0)--(1,1)--(2,1)--(3,2)--(4,3)--(4,4);
\end{tikzpicture}
\end{center}
\caption{Two $(4,4)$-Delannoy paths. The first is quasi-diagonal and turns left, diagonally, right, and diagonally. The second is not quasi-diagonal.} \label{fig:quasi-diag}
\end{figure}
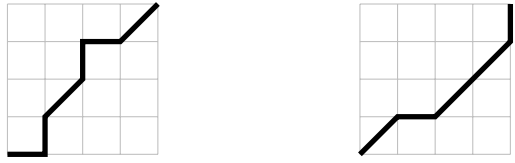

We require two lemmas before proving the theorem. For an $(n,n)$-Delannoy path $p$, let $A_p \in \Mat_{\bR^{(n)}}$ be the indicator function of the orbit on $\bR^{(n)} \times \bR^{(n)}$ corresponding to $p$ via Proposition~\ref{prop:orbit}. We say that $p$ is \emph{quasi-diagonal} if it passes through every vertex along the main diagonal. Suppose $p$ is quasi-diagonal. At the $i$th square along along the main diagonal, $p$ can behave in three possible ways:
\begin{itemize}
\item $p$ \emph{turns diagonally} if it goes from $(i-1,i-1)$ to $(i,i)$;
\item $p$ \emph{turns right} if it goes from $(i-1,i-1)$ to $(i-1,i)$ and then to $(i,i)$; or
\item $p$ \defn{turns left} if it goes from $(i-1,i-1)$ to $(i,i-1)$ and then to $(i,i)$.
\end{itemize}
See Figure~\ref{fig:quasi-diag} for an illustration.

\begin{lemma} \label{lem:simple-2}
Let $\lambda \in \Lambda$ have length $n$, let $p$ be a $(n,n)$-Delannoy path, and let $\phi_{\bI}$ be a generator of $L_{\lambda}$. If $p$ is not quasi-diagonal then $A_p \cdot \phi_{\bI} = 0$. Suppose $p$ is quasi-diagonal. For $1 \le i \le n$, define
\begin{displaymath}
\epsilon_i(\lambda,p) = \begin{cases}
1 & \text{if $p$ turns diagonally in the $(i,i)$ square} \\
-1 & \text{if $\lambda_i=\wa$ and $p$ turns right in the $(i,i)$ square} \\
-1 & \text{if $\lambda_i=\wb$ and $p$ turns left in the $(i,i)$ square} \\
0 & \text{otherwise} \end{cases}
\end{displaymath}
and define $\epsilon(\lambda,p) = \prod_{i=1}^n \epsilon_i(\lambda,p)$. Then
\begin{displaymath}
A_p \cdot \phi_{\bI}=\epsilon(\lambda,p) \cdot \phi_{\bI}.
\end{displaymath}
\end{lemma}

\begin{proof}
Let $X \subset \bR^{(n)} \times \bR^{(n)}$ be the orbit corresponding to $p$. We have
\begin{displaymath}
(A_p \phi_{\bI})(y) = \vol(\{x \in \bR^{(n)} \mid \text{$x \in \bI$ and $(y,x) \in X$}\}).
\end{displaymath}
The set on the right side above factors as $J_1 \times \cdots \times J_n$, where $J_i$ is a subset of $I_i$. To be a bit more precise, fix $y$ and $i$. Then the condition $(y,x) \in X$ puts us in one of two cases, as far as $x_i$ is concerned:
\begin{itemize}
\item $x_i=y_j$ for some $j$; in this case, $J_i$ is $\{y_j\} \cap I_i$.
\item $y_j<x_i<y_{j+1}$ for some $j$; in this case, $J_i$ is $(y_j,y_{j+1}) \cap I_i$. (We use the convention $y_0=-\infty$ and $y_{n+1}=+\infty$.)
\end{itemize}

Suppose $y \not \in \bI$. We claim that $(A_p \phi_{\bI})(y)=0$. Indeed, if any $J_i$ is empty or all of $I_i$ then the total volume is~0. To ensure $J_i$ is non-empty and not all of $I_i$, it must contain some $y_j$. If this holds for all $i$, then each $I_i$ contains exactly one $y_j$. Since the $I_i$'s and $y_j$'s are both in increasing order, we must have $y_i \in I_i$ for each $i$, that is, $y \in \bI$.

Now suppose $y \in \bI$. First we claim that if $y$ is not quasi-diagonal then $(A_p \phi_{\bI})(y)=0$. Indeed if $(A_p \phi_{\bI})(y) \neq 0$ then each $I_i$ contains both an $x_i$ and $y_i$, so before we enter the interval the path is at the vertex $(i-1,i-1)$ and as we exit it is at $(i,i)$.  Second, if $p$ is quasi-diagonal we want to show that $\epsilon_i(\lambda,p) = \vol(J_i)$.  We split into cases accordingly using the above description of $J_i$ (with $j=i$ or $j=i+1$ according to how the path turns in the $i$th square):

\begin{itemize}
\item If the path turns diagonally in the $i$th square, then $\epsilon_i(\lambda,p) = 1$ and $J_i = \{y_i\} \cap I_i = \{y_i\}$, which has volume $1$. 
\item If $\lambda_i=\wa$ and the path turns right in the $i$th square, then $\epsilon_i(\lambda,p) = -1$ and $J_i = (y_{i-1}, y_{i}) \cap (x_i,x_{i+1}] = (x_i, y_i)$ has volume $-1$.
\item If $\lambda_i=\wb$ and the path turns right in the $i$th square, then $\epsilon_i(\lambda,w) = 0$ and $J_i = (y_{i-1}, y_{i}) \cap [x_i,x_{i+1}) = [x_i, y_i)$ has volume $0$.
\item If $\lambda_i=\wa$ and the path turns left in the $i$th square, then $\epsilon_i(\lambda,w) = 0$ and $J_i = (y_{i}, y_{i+1}) \cap (x_i,x_{i+1}] = (y_i, x_{i+1}]$ has volume $0$.
\item If $\lambda_i=\wb$ and the path turns right in the $i$th square, then $\epsilon_i(\lambda,w) = -1$ and $J_i = (y_{i}, y_{i+1}) \cap [x_i,x_{i+1}) = (y_i, x_{i+1})$ has volume $-1$. \qedhere
\end{itemize}
\end{proof}

\begin{lemma} \label{lem:simple-1}
Let $p \colon \bR^{(n)} \to \bR^{(n-1)}$ be the projection away from the last coordinate. Then $p_*(L_{\lambda})=0$.
\end{lemma}

\begin{proof}
It suffices to show that $p_*(\phi_{\bI})=0$ for all $\phi_{\bI}$ in $L_{\lambda}$. Fix $(x_1, \ldots, x_{n-1})$ in $\bR^{(n-1)}$, and let $X \subset \bR$ be the set of $x_n$'s such that $(x_1, \ldots, x_n) \in \bI$. Then the value of $p_*(\phi_{\bI})$ at $(x_1, \ldots, x_{n-1})$ is the volume of $X$. If $x_i \not\in I_i$ for some $1 \le i \le n-1$ then $X$ is empty and the volume is~0. If $x_i \in I_i$ for all $i$ then $X=I_n$, which also has volume~0.
\end{proof}

\begin{proof}[Proof of Theorem~\ref{thm:simple}]
(a) Lemma~\ref{lem:simple-2} shows that $A_p$ acts by a scalar on $L_{\lambda}$. Since the $A_p$'s form a basis of the endomorphism algebra of $\cC(\bR^{(n)})$, it follows that every endomorphism acts as a scalar on $L_{\lambda}$. Since $\cC(\bR^{(n)})$ is semi-simple, every endomorphism of $L_{\lambda}$ is induced by an endomorphism of $\cC(\bR^{(n)})$. We thus find $\End_{\ul{G}}(L_{\lambda})=k$. Since $L_{\lambda}$ is semi-simple, it follows that $L_{\lambda}$ is simple. Note that since every endomorphism of $\cC(\bR^{(n)})$ preserves $L_{\lambda}$, it follows that $L_{\lambda}$ has multiplicity one in $\cC(\bR^{(n)})$. (We will prove a more general statement in Theorem~\ref{thm:simple}.)

(b) Let $\lambda$ and $\mu$ be distinct weights; we show that $L_{\lambda}$ and $L_{\mu}$ are non-isomorphic. Put $n=\ell(\lambda)$ and $m=\ell(\mu)$, and assume $m \le n$ without loss of generality. First suppose that $m=n$. By Lemma~\ref{lem:simple-2}, we see that $A_p$ acts by distinct scalars on $L_{\lambda}$ and $L_{\mu}$ for appropriate $p$: in fact, note that there is a unique quasi-diagonal path $p$ having no diagonal turns for which $A_p$ is non-zero on $L_{\lambda}$, and this path determines $\lambda$. Thus $L_{\lambda}$ and $L_{\mu}$ are distinct subspaces of $\cC(\bR^{(n)})$. Since each has multiplicity one in $\cC(\bR^{(n)})$, it follows that they are non-isomorphic.

Now suppose $m<n$. We first claim that if $L$ is a simple appearing in $\cC(\bR^{(n-1)})$ then $L$ is not isomorphic to $L_{\lambda}$. Indeed, let $p$ be the projection in Lemma~\ref{lem:simple-1}, and note that $p^* \colon \cC(\bR^{(n-1)}) \to \cC(\bR^{(n)})$ is injective. If $L$ were isomorphic to $L_{\lambda}$ then we would have $L_{\lambda}=p^*(L)$ since $L_{\lambda}$ has multiplicity one in $\cC(\bR^{(n)})$. But $p_*(L_{\lambda})=0$ by Lemma~\ref{lem:simple-1}, while $p_*p^*=-1$ (since the fibers of $p$ are open intervals), a contradiction. This proves the claim. Now simply note that $L_{\mu}$ appears in $\cC(\bR^{(n-1)})$, as the standard projection $\bR^{(n-1)} \to \bR^{(m)}$ gives an injection $\cC(\bR^{(m)}) \to \cC(\bR^{(n-1)})$.
\end{proof}

\subsection{Evaluation maps} \label{ss:eval}

For $a \in \bR^{(n)}$, define a map
\begin{displaymath}
\ev_a \colon \cC(\bR^{(n)}) \to k, \qquad \ev_a(\phi)=\phi(a).
\end{displaymath}
We call $\ev_a$ the \defn{evaluation at $a$} map. It is clear that $\ev_a$ is a map of $\ul{G}(a)$-representations. We will require the following simple observation:

\begin{proposition} \label{prop:eval-ind}
Let $a_1, \ldots, a_r \in \bR^{(n)}$ be distinct, and let $\lambda$ be a weight of length $n$. Then $\ev_{a_1}, \ldots, \ev_{a_r}$ restrict to linearly independent functionals on $L_{\lambda}$.
\end{proposition}

\begin{proof}
Fix $1 \le i \le r$. We can then find $\bI=(I_1<\cdots<I_n)$ of type $\lambda$ such that $a_i \in \bI$, but $a_j \not\in \bI$ for all $j \ne i$. It follows that $\ev_{a_i}(\phi_{\bI})=1$, but $\ev_{a_j}(\phi_{\bI})=0$ for $j \ne i$. This completes the proof.
\end{proof}

\subsection{Decomposition of Schwartz space}

We now determine the simple decomposition of Schwartz space:

\begin{theorem} \label{thm:decomp}
For $n \ge 0$, we have an isomorphism
\begin{displaymath}
\cC(\bR^{(n)}) = \bigoplus_{\ell(\lambda) \le n} L_{\lambda}^{\oplus m_{\lambda}}, \qquad m_{\lambda}=\binom{n}{\ell(\lambda)}.
\end{displaymath}
\end{theorem}

\begin{proof}
Define $m_{\lambda}$ be the multiplicity of $L_{\lambda}$ in $\cC(\bR^{(n)})$, for any weight $\lambda$. We first give a lower bound for $m_{\lambda}$. Let $a \in \bR^{(n)}$, and put $s=\ell(\lambda)$. We have
\begin{displaymath}
m_{\lambda} = \dim \Hom_{\ul{G}}(L_{\lambda}, \cC(\bR^{(n)})) = \dim \Hom_{\ul{G}(a)}(L_{\lambda}, k),
\end{displaymath}
where the second step is Frobenius reciprocity. The functionals $\ev_b$, with $b \subset a$ of cardinality $s$, are $\ul{G}(a)$-equivariant and linearly independent by Proposition~\ref{prop:eval-ind}. Thus $m_{\lambda} \ge \binom{n}{s}$.

We have
\begin{displaymath}
\cC(\bR^{(n)}) = \big( \bigoplus_{\lambda} L_{\lambda}^{\oplus m_{\lambda}} \big) \oplus X
\end{displaymath}
where $X$ contains the simples not of the form $L_{\lambda}$. We now compare the endomorphism rings of the two sides. The dimension of $\End_{\ul{G}}(\cC(\bR^{(n)}))$ is the central Delannoy number $D(n)$ (Corollary~\ref{cor:hom-dim}), while $\End_{\ul{G}}(L_{\lambda})$ is one-dimensional (Theorem~\ref{thm:simple}). We thus find
\begin{displaymath}
D(n) = \sum_{\lambda} m_{\lambda}^2 + d,
\end{displaymath}
where $d$ is the dimension of $\End_{\ul{G}}(X)$. We have
\begin{displaymath}
\sum_{\lambda} m_{\lambda}^2 = \sum_{s \ge 0} \sum_{\ell(\lambda)=s} m_{\lambda}^2 \ge \sum_{s \ge 0} 2^s \binom{n}{s}^2 = D(n),
\end{displaymath}
where in the final step we used a well-known formula for $D(n)$ \cite[Example~2]{Sulanke1}. It follows that we must have $m_{\lambda}=\binom{n}{s}$ where $s=\ell(\lambda)$, and $d=0$ (which implies $X=0$).
\end{proof}

\begin{example} \label{ex:decomp}
For $n=1,2$, the theorem gives the following:
\begin{align*}
\cC(\bR) &= L_{\wa} \oplus L_{\wb} \oplus \bbone \\
\cC(\bR^{(2)}) &= L_{\wa\wa} \oplus L_{\wa\wb} \oplus L_{\wb\wa} \oplus L_{\wb\wb} \oplus L_{\wa}^{\oplus 2} \oplus L_{\wb}^{\oplus 2} \oplus \bbone.
\end{align*}
Here $\bbone$ is the trivial representation. The first decomposition above was obtained ``by hand'' in \cite[\S 17.6]{HarmanSnowden}.
\end{example}

\begin{remark} \label{rmk:selberg}
The following two quantities are equal to each other, and to $D(n)$
\begin{displaymath}
\sum_{k=0}^n \frac{(2n-k)!}{(n-k)! (n-k)! k!} = \sum_{k=0}^n \binom{n}{k}^2 2^k,
\end{displaymath}
see \cite[Examples~1,3]{Sulanke1}. In light of Theorem~\ref{thm:decomp}, the above equality can be seen as analogous to the Selberg trace formula. The left side is the geometric side, as one obtains it by counting the orbits of $G$ on the geometric object $\bR^{(n)} \times \bR^{(n)}$. The right side is the spectral side, as one easily obtains it from the simple decomposition of $\cC(\bR^{(n)})$.
\end{remark}

The theorem has a number of corollaries.

\begin{corollary} \label{cor:schwartz-length}
The length of $\cC(\bR^{(n)})$ is $3^n$.
\end{corollary}

\begin{proof}
From the decomposition of $\cC(\bR^{(n)})$ we find that its length is
\begin{displaymath}
\sum_{\lambda} \binom{n}{\ell(\lambda)} = \sum_{s \ge 0} \sum_{\ell(\lambda)=s} \binom{n}{s} = \sum_{s \ge 0} \binom{n}{s} 2^s = 3^n,
\end{displaymath}
where in the third step we use that there are $2^s$ weights of length $s$.
\end{proof}

\begin{corollary} \label{cor:inv-dim}
Let $\lambda$ be a weight of length $s$ and let $a \in \bR^{(n)}$. Then the dimension of the $G(a)$-invariants of $L_{\lambda}$ is $\binom{n}{s}$.
\end{corollary}

\begin{proof}
The dimension of $L_{\lambda}^{G(a)}$ is the multiplicity of $L_{\lambda}$ in $\cC(\bR^{(n)})$ by Frobenius reciprocity, which is $\binom{n}{s}$. (We note that $G$-invariants and $\ul{G}$-invariants coincide by \cite[\S 11.7]{HarmanSnowden}.)
\end{proof}

\begin{corollary} \label{cor:all-simples}
Every simple object of $\uRep(G)$ is isomorphic to $L_{\lambda}$ for some $\lambda$.
\end{corollary}

\begin{proof}
Every object of $\uRep(G)$ is a quotient of a sum of $\cC(\bR^{(n)})$'s, so all simples must appear in some Schwartz space $\cC(\bR^{(n)})$. The theorem shows that the only simples in $\cC(\bR^{(n)})$ are those of the form $L_{\lambda}$.
\end{proof}

Recall that $a_{\lambda}$ is the class of the simple $L_{\lambda}$ in the Grothendieck group $\rK$. We note one final corollary, that we will use constantly (often without mention).

\begin{corollary} \label{cor:K-basis}
The elements $\{a_{\lambda}\}_{\lambda \in \Lambda}$ are a $\bZ$-basis of $\rK$.
\end{corollary}

\begin{remark} \label{rmk:K-filt-basis}
We also see that the elements $a_{\lambda}$ with $\ell(\lambda) \le n$ form a $\bZ$-basis of $\rK_{\le n}$.
\end{remark}

\begin{remark} \label{rmk:field-ind}
Write $\rK(k)$ for the Grothendieck group, indicating the dependence on the field $k$. If $k'$ is a second field, Corollary~\ref{cor:K-basis} gives an isomorphism $\iota \colon \rK(k) \to \rK(k')$ via $a_{\lambda} \mapsto a_{\lambda}$. In fact, $\iota$ is canonical and respects the extra structure on $\rK$, as we now explain.

There are two main cases to discuss. First, if $k$ is the prime subfield of $k'$ then it follows from the construction of $L_{\lambda}$ that extension of scalars induces $\iota$. Second, we have a natural map $\rK(\bQ) \to \rK(\bF_p)$ obtained by picking an integral lattice and reducing mod $p$. One can show that this induces $\iota$; here one must use the existence of a good integral form of $\uRep(G)$, and Proposition~\ref{prop:simple-dual-descr} is helpful too.
\end{remark}

\subsection{Duality} \label{ss:dual}

For a weight $\lambda$, let $\lambda^{\vee}$ be the weight obtained by interchanging $\wa$ and $\wb$.

\begin{proposition} \label{prop:simple-dual}
Let $\lambda$ be a weight of length $n$. The pairing
\begin{displaymath}
\langle, \rangle \colon L_{\lambda} \times L_{\lambda^{\vee}} \to k, \qquad
\langle \phi, \psi \rangle = \int_{\bR^{(n)}} \phi(x) \psi(x) dx
\end{displaymath}
identifies $L_{\lambda^{\vee}}$ with the dual of $L_{\lambda}$.
\end{proposition}

\begin{proof}
The pairing is clearly $\ul{G}$-equivariant, and thus induces a map of $\ul{G}$-modules $L_{\lambda^{\vee}} \to L_{\lambda}^{\vee}$. Since both objects are simple, this map is an isomorphism provided it is non-zero. It thus suffices to exhibit one non-zero pairing between functions in $L_{\lambda}$ and $L_{\lambda^{\vee}}$.

Choose real numbers
\begin{displaymath}
a_1<b_1<c_1<d_1<\cdots<a_n<b_n<c_n<d_n.
\end{displaymath}
Define intervals $I_i$ and $J_i$ as follows:
\begin{itemize}
\item If $\lambda_i=\wa$ then $I_i=(a_i,c_i]$ and $J_i=[b_i,d_i)$.
\item If $\lambda_i=\wb$ then $I_i=[b_i,d_i)$ and $J_i=(a_i,c_i]$.
\end{itemize}
Then $\bI=(I_1<\cdots<I_n)$ and $\bJ=(J_1<\cdots<J_n)$ have types $\lambda$ and $\lambda^{\vee}$, and so $\phi_{\bI}$ and $\phi_{\bJ}$ belong to $L_{\lambda}$ and $L_{\lambda^{\vee}}$. The pairing $\langle \phi_{\bI}, \phi_{\bJ} \rangle$ is equal to the volume of $\bI \cap \bJ$. By our construction, this set factors as $\prod_{i=1}^n (I_i \cap J_i)$, and $I_i \cap J_i=[b_i,c_i]$ is a closed interval. Thus the volume is 1, which completes the proof.
\end{proof}

\begin{remark}
The group $G$ admits an outer automorphism $\sigma$ obtained by reversing the real line. This induces an involution of $\rK$ that is different from the one induced by duality: we have $a_{\lambda}^{\sigma}=a_{\mu}$, where $\mu$ is the reverse of $\lambda^{\vee}$.
\end{remark}

\subsection{A dual description of the simples} \label{ss:simple-dual-descr}

We defined $L_{\lambda}$ as the subrepresentation of $\cC(\bR^{(n)})$ generated by certain elements. We now give a dual description of $L_{\lambda}$, that is, we describe it as the subrepresentation of $\cC(\bR^{(n)})$ consisting of elements satisfying certain linear equations. We fix $n \ge 0$ in what follows.

To begin, we explicitly describe the embeddings of simple objects into Schwartz space. For a subset $S$ of $[n]$ of cardinality $m$, let $p_S \colon \bR^{(n)} \to \bR^{(m)}$ be the projection onto the $S$ coordinates. This map induces a map of representations $p_S^* \colon \cC(\bR^{(m)}) \to \cC(\bR^{(n)})$.

\begin{proposition} \label{prop:simple-embeddings}
Let $\lambda$ be a weight of length $\le n$. For $S \subset [n]$ of cardinality $m$, let $i_S \colon L_{\lambda} \to \cC(\bR^{(n)})$ be the restriction of $p_S^*$ to $L_{\lambda}$. Then the $i_S$ form a basis of the multiplicity space $\Hom_{\ul{G}}(L_{\lambda}, \cC(\bR^{(n)}))$.
\end{proposition}

\begin{proof}
By Theorem~\ref{thm:decomp}, the multiplicity space in question has dimension $\binom{n}{m}$. It thus suffices to prove that the elements $i_S$ are linearly independent. To see this, pick $a \in \bR^{(n)}$. Then the composition
\begin{displaymath}
\xymatrix@C=4em{
L_{\lambda} \ar[r]^-{i_S} & \cC(\bR^{(n)}) \ar[r]^-{\ev_a} & k }
\end{displaymath}
is equal to $\ev_{p_S(a)}$. Since these functionals are linearly independent on $L_{\lambda}$ (Proposition~ \ref{prop:eval-ind}), it follows that the $i_S$ are linearly independent.
\end{proof}

For $1 \le i \le n$, let $p_i \colon \bR^{(n)} \to \bR^{(n-1)}$ be the projection away from the $i$th coordinate; this is just $p_S$ where $S=[n] \setminus \{i\}$.

\begin{proposition} \label{prop:XY}
Define subrepresentations $X$ and $Y$ of $\cC(\bR^{(n)})$ by
\begin{displaymath}
X = \bigcap_{i=1}^n \ker((p_i)_*), \qquad Y = \sum_{i=1}^n \im(p_i^*).
\end{displaymath}
Then $X$ (resp.\ $Y$) is the sum of the $L_{\lambda}$-isotypic spaces of $\cC(\bR^{(n)})$ where $\ell(\lambda)=n$ (resp.\ $\ell(\lambda)<n$). In particular, $\cC(\bR^{(n)})=X \oplus Y$.
\end{proposition}

\begin{proof}
It is clear that $Y$ only contains simples $L_{\lambda}$ with $\ell(\lambda)<n$. Now, if $S$ is proper subset of $[n]$ then it does not contain some element $i$, and so $p_S$ factors through $p_i$; it follows that $\im(p_S^*) \subset \im(p_i^*)$. It thus follows from Proposition~\ref{prop:simple-embeddings} that $Y$ contains every copy of $L_{\lambda}$ in $\cC(\bR^{(n)})$ if $\ell(\lambda)<n$. We have thus verified the statement about $Y$.

The maps $(p_i)_*$ and $p_i^*$ are adjoint with respect to the standard pairings on $\cC(\bR^{(n)})$ and $\cC(\bR^{(n-1)})$. It follows that $X$ is exactly the orthogonal complement of $Y$. This yields the statement about $X$.
\end{proof}

We now come to the desired description of $L_{\lambda}$.

\begin{proposition} \label{prop:simple-dual-descr}
Let $\lambda$ be a weight of length $n$. Then $L_{\lambda}$ consists of those functions $\phi$ in $\cC(\bR^{(n)})$ satisfying the following two conditions:
\begin{enumerate}
\item We have $(p_i)_*(\phi)=0$ for all $1 \le i \le n$.
\item For every weight $\mu \ne \lambda^{\vee}$ of length $n$ and every $\psi \in L_{\mu}$ we have $\langle \phi, \psi \rangle=0$.
\end{enumerate}
Moreover, in (b) it is enough to consider functions $\psi$ of the form $\phi_{\bI}$.
\end{proposition}

\begin{proof}
The first statement follows immediately from Propositions~\ref{prop:simple-dual} and~\ref{prop:XY}. We now explain that the second statement (about the $\phi_{\bI}$ being enough). Suppose $\phi \in \cC(\bR^{(n)})$ is orthogonal to all $\phi_{\bI}$ in some $L_{\mu}$. Let $U$ be an open subgroup of $G$ such that $\phi$ is $U$ invariant. The $\phi_{\bI}$ generate $L_{\mu}$ as a representation of $U$ (see \cite[Corollary~11.23]{HarmanSnowden} and the preceeding discussion). It follows that $\phi$ is orthogonal to all functions in $L_{\mu}$.
\end{proof}

\subsection{The concatenation product} \label{ss:concat}

We define the \defn{concatenation product} to be the unique bilinear map
\begin{displaymath}
\odot \colon \rK \times \rK \to \rK
\end{displaymath}
satisfying
\begin{displaymath}
a_{\lambda} \odot a_{\mu} = a_{\lambda \mu}
\end{displaymath}
for all weights $\lambda$ and $\mu$; here $\lambda \mu$ simply denotes the concatenation of the words $\lambda$ and $\mu$. This definition is justified since the $a_{\lambda}$'s form a $\bZ$-basis of $\rK$ (Corollary~\ref{cor:K-basis}). The concatenation product gives $\rK$ the structure of an associative unital ring; it fact, it is simply the non-commutative polynomial ring on the variables $a_{\wa}$ and $a_{\wb}$.

\section{Invariants} \label{s:inv}

In the previous section, we computed the dimension of the space of $G(a)$-invariants on the simple $L_{\lambda}$ (Corollary~\ref{cor:inv-dim}). In this section, we explicitly identify a basis of this space. Using this, we write down the projection operator for $L_{\lambda}$ in $\cC(\bR^{(n)})$, which allows us to compute the categorical dimension of $L_{\lambda}$ (Corollary~\ref{cor:cat-dim}). We close this section by showing how one can determine the structure an arbitrary representation by analyzing its invariant spaces.

\subsection{The key set}

Let $\lambda$ be a weight of length $n$. We define
\begin{displaymath}
\Psi_{\lambda} \subset \bR^{(n)} \times \bR^{(n)}
\end{displaymath}
to be the subset consisting of all pairs $(x,y)$ satisfying the following conditions, for $1 \le i \le n$:
\begin{itemize}
\item If $\lambda_i=\wa$ then $x_i \le y_i$.
\item If $\lambda_i=\wb$ then $y_i \le x_i$.
\item If $i<n$ then $x_i<y_{i+1}$ and $y_i<x_{i+1}$.
\end{itemize}
Here are some examples:
\begin{itemize}
\item For $\lambda=\wa\wa$ we have $x_1 \le y_1 < x_2 \le y_2$.
\item For $\lambda=\wa\wb$ we have $x_1 \le y_1 < y_2 \le x_2$.
\item For $\lambda=\wb\wa$ we have $y_1 \le x_1 < x_2 \le y_2$.
\item For $\lambda=\wb\wb$ we have $y_1 \le x_1 < y_2 \le x_2$.
\end{itemize}
It is clear that $\Psi_{\lambda}$ is a $G$-stable set; however, it is not a single orbit. The set $\Psi_{\lambda}$ will play a key role in the rest of this section. Note that if $\tau$ is the map on $\bR^{(n)} \times \bR^{(n)}$ given by $\tau(x,y)=(y,x)$ then $\tau(\Psi_{\lambda})=\Psi_{\lambda^{\vee}}$.

\begin{remark} \label{rmk:key-set}
Let $p$ be the quasi-diagonal $(n,n)$-Delannoy path where the $i$th turn is rigjt (resp.\ left) if $\lambda_i=\wa$ (resp.\ $\lambda_i=\wb$). The the $G$-orbit on $\bR^{(n)} \times \bR^{(n)}$ corresponding to $p$ in Proposition~\ref{prop:orbit} is defined just like $\Psi_{\lambda}$, except with all $\le$ changed to $<$.
\end{remark}

\subsection{Invariant functions}

For $a \in \bR^{(n)}$, define $\Psi_{\lambda}^a \subset \bR^{(n)}$ to be the set of points $x$ such that $(x,a)$ belongs to $\Psi_{\lambda}$, and let $\psi_{\lambda}^a \in \cC(\bR^{(n)})$ be the characteristic function of $\Psi_{\lambda}^a$. The main interest in these functions lies in the following result:

\begin{proposition} \label{prop:inv}
The space of $G(a)$-invariants in $L_{\lambda}$ is spanned by $\psi^a_{\lambda}$.
\end{proposition}

\begin{proof}
Since $\Psi^a_{\lambda}$ is stable by $G(a)$, it follows that $\psi^a_{\lambda}$ is a $G(a)$-invariant element of $\cC(\bR^{(n)})$. Since we already know that that $G(a)$-invariants of $L_{\lambda}$ are one-dimensional (Corollary~\ref{cor:inv-dim}), it is enough to show that $\psi^a_{\lambda}$ belongs to $L_{\lambda}$. For $1 \le i \le n$, let $p_i \colon \bR^{(n)} \to \bR^{(n-1)}$ be the projection away from the $i$th coordinate. Recall (Proposition~\ref{prop:simple-dual-descr}) that an element $\theta$ of $\cC(\bR^{(n)})$ belongs to $L_{\lambda}$ if and only if the following two conditions hold:
\begin{enumerate}
\item $(p_i)_*(\theta)=0$ for all $1 \le i \le n$
\item $\theta$ is orthogonal to all $L_{\mu}$'s with $\ell(\mu)=n$, except for $\mu=\lambda^{\vee}$.
\end{enumerate}
We verify these two conditions for $\theta=\psi^a_{\lambda}$.

(a) Fix $(x_1, \ldots, \hat{x}_i, \ldots, x_n)$ in $\bR^{(n-1)}$. Let $X \subset \bR$ be the set of $x_i$'s such that $(x_1, \ldots, x_i, \ldots, x_n)$ belongs to $\Psi^a_{\lambda}$. Then the value of $(p_i)_* \psi^a_{\lambda}$ at $(x_1, \ldots, \hat{x}_i, \ldots, x_n)$ is the volume of $X$. We claim that $X$ is a half-open interval (or empty); this will prove the statement as such a set has volume~0. First suppose that $\lambda_i=\wa$. There are three cases:
\begin{itemize}
\item if $i=1$ then $X=(-\infty, a_1]$
\item if $i>1$ and $\lambda_{i-1}=\wa$ then $X=(a_{i-1},a_i]$
\item if $i>1$ and $\lambda_{i-1}=\wb$  then $X=(x_{i-1}, a_i]$.
\end{itemize}
The case $\lambda_i=\wb$ is similar. This proves the claim.

(b) Suppose $\mu$ is a weight of length $n$ that is not equal to $\lambda^{\vee}$; this precisely means that there is some $1 \le i \le n$ such that $\mu_i=\lambda_i$. To show that $\psi^a_{\lambda}$ is orthogonal to $L_{\mu}$, it suffices to show that it is orthogonal to the generators $\phi_{\bI}$. Thus let $\phi_{\bI}$ be given, where $\bI=(I_1<\cdots<I_n)$ are intervals of type $\mu$. By Fubini's theorem, we have
\begin{displaymath}
\langle \psi^a_{\lambda}, \phi_{\bI} \rangle = \int_{\bR^{(n)}} \psi^a_{\lambda}(x) \phi_{\bI}(x) dx
=\int_{I_1} \cdots \int_{I_n} \psi^a_{\lambda}(x_1, \ldots, x_n) dx_1 \cdots dx_n.
\end{displaymath}
In fact, we can put the integrals on the right in whatever order we prefer. We do the $i$th integral first. Thus fix $x_j \in I_j$ for all $j \ne i$, and consider the integral
\begin{displaymath}
\int_{I_i} \psi_{\lambda}^a(x_1, \ldots, x_i, \ldots, x_n) dx_i
\end{displaymath}
It suffices to show that this integral vanishes. Let $X \subset \bR$ be defined as in (a). Then the above integral is simply the volume of $X \cap I_i$. Suppose $\lambda_i=\wa$. Then, as we saw in (a), $X$ is a half-open interval containing its right endpoint (or empty). Since $I_i$ has type $\mu_i=\lambda_i$, it is a half-open interval containing its right endpoint. Thus $X \cap I_i$ is also a half-open interval containing its right endpoint (or empty), and therefore has volume~0. The case $\lambda_i=\wb$ is similar.
\end{proof}

Since the $\ul{G}(a)$-invariant space of $L_{\lambda}$ is one-dimensional, it follows from semi-simplicity that the $\ul{G}(a)$-co-invariant space is also one-dimensional, i.e., up to scalar multiples, there is a unique $\ul{G}(a)$-equivariant functional $L_{\lambda} \to k$. We know two such functionals, namely, the evaluation map $\ev_a$ (see \S \ref{ss:eval}), and the pairing $\langle -, \psi^a_{\lambda^{\vee}} \rangle$. These two functionals are necessarily linearly dependent; in fact, they are equal:

\begin{proposition} \label{prop:psi-ev}
For any $\phi \in L_{\lambda}$, we have
\begin{displaymath}
\phi(a) = \int_{\bR^{(n)}} \phi(x) \psi^a_{\lambda^{\vee}}(x) dx.
\end{displaymath}
In other words, $\ev_a=\langle -, \psi^a_{\lambda^{\vee}} \rangle$.
\end{proposition}

\begin{proof}
By the above discussion, we have $\ev_a = \beta \cdot \langle -, \psi^a_{\lambda^{\vee}} \rangle$, for some scalar $\beta$. Let $\bI=(I_1<\cdots<I_n)$ be intervals of type $\lambda$ such that $a_i$ is the closed endpoint of $I_i$ for all $i$. We have $\ev_a(\phi_{\bI})=1$. To show $\beta=1$, it suffices to show that $\langle \phi_{\bI}, \psi^a_{\lambda^{\vee}} \rangle=1$. The value of this pairing is the volume of $X=\Psi^a_{\lambda^{\vee}} \cap \bI$. We claim that $X=\{a\}$, which will complete the proof. It is clear that $a \in X$. Conversely, suppose $b \in X$. Say $\lambda_i=\wa$. Since $b_i \in I_i$, we have $b_i \le a_i$, and since $b \in \Psi^a_{\lambda^{\vee}}$, we have $a_i \le b_i$. Thus $b_i=a_i$. The case $\lambda_i=\wb$ is similar. Thus $a=b$, as required.
\end{proof}

\begin{corollary} \label{cor:psi-ind}
Let $a_1, \ldots, a_r$ be distinct elements of $\bR^{(n)}$. Then $\psi^{a_1}_{\lambda}, \ldots, \psi^{a_r}_{\lambda}$ are linearly independent functions.
\end{corollary}

\begin{proof}
This follows from the above proposition and the fact that $\ev_{a_1}, \ldots, \ev_{a_r}$ are linearly independent functionals on $L_{\lambda^{\vee}}$ (Proposition~\ref{prop:eval-ind}).
\end{proof}

\begin{corollary}
Let $m \ge n$, and let $b \in \bR^{(m)}$. Then the functions $\psi_{\lambda}^a$, with $a$ an $n$-element subset of $b$, form a basis of $L_{\lambda}^{G(b)}$.
\end{corollary}

\begin{proof}
The $\psi_{\lambda}^a$ belong to $L_{\lambda}^{G(b)}$ by Proposition~\ref{prop:inv}. By Corollary~\ref{cor:psi-ind}, these functions are linearly independent. Since $L_{\lambda}^{G(b)}$ has dimension $\binom{m}{n}$ (Corollary~\ref{cor:inv-dim}), the result follows.
\end{proof}

\subsection{The projection map}

Define
\begin{displaymath}
\pi_{\lambda} \colon \bR^{(n)} \times \bR^{(n)} \to k
\end{displaymath}
to be the characteristic function of $\Psi_{\lambda}$. We regard $\pi_{\lambda}$ as an $\bR^{(n)} \times \bR^{(n)}$ matrix; as such it defines an endomorphism of $\cC(\bR^{(n)})$.

\begin{proposition} \label{prop:projector}
The matrix $\pi_{\lambda}$ defines the projection map $\cC(\bR^{(n)}) \to L_{\lambda}$, i.e., $\pi_{\lambda}$ is idempotent and its image is exactly $L_{\lambda}$.
\end{proposition}

\begin{proof}
The ``columns'' of the matrix $\pi_{\lambda}$ are the functions $\psi^a_{\lambda}$ considered above. Since these belong to $L_{\lambda}$ (Proposition~\ref{prop:inv}), it follows that the image of $\pi_{\lambda}$ is contained in $L_{\lambda}$. To complete the proof, it suffices to show that $\pi_{\lambda}$ acts by the identity on elements of $L_{\lambda}$. This is exactly Proposition~\ref{prop:psi-ev}; note that $\psi^a_{\lambda^{\vee}}(x)=\pi_{\lambda}(a,x)$.
\end{proof}

Recall that an object in a pre-Tannakian category has a ``categorical dimension'' \cite[\S 4.7]{Etingof}.

\begin{corollary} \label{cor:cat-dim}
The categorical dimension of $L_{\lambda}$ is $(-1)^{\ell(\lambda)}$.
\end{corollary}

\begin{proof}
The categorical dimension $d$ of $L_{\lambda}$ is equal to the categorical trace of $\pi_{\lambda}$, which is just the integral of the diagonal of $\pi_{\lambda}$ \cite[Proposition~8.11]{HarmanSnowden}. Since $\Psi_{\lambda}$ contains the diagonal, we find $d=\vol(\bR^{(n)})$, which is $(-1)^n$, as claimed.
\end{proof}

\begin{remark}
Suppose $k$ has positive characteristic $p$. The paper \cite{EtingofHarmanOstrik} defines a notion of ``$p$-adic dimension'' for objects in (certain) pre-Tannakian categories over $k$. The $p$-adic dimension of $L_{\lambda}$ is $(-1)^{\ell(\lambda)}$. One can see this by using the integral version of $\uRep(G)$, and the fact that $L_{\lambda}$ exists there.
\end{remark}

\begin{proposition} \label{prop:delta}
There exists a unique ring homomorphism $\delta \colon \rK \to \bZ$ satisfying $\delta(a_{\lambda})=(-1)^{\ell(\lambda)}$.
\end{proposition}

\begin{proof}
Suppose $k$ has characteristic~0. Then the map $\rK \to k$ assigning to each class its categorical dimension is a ring homomorphism. Since this takes values in $\bZ \subset k$ by Corollary~\ref{cor:cat-dim}, the result follows. In positive characteristic, one can either use $p$-adic dimension, or simply appeal to the fact that $\rK$ is independent of the field $k$ (see Remark~\ref{rmk:field-ind}).
\end{proof}

\subsection{Operators on invariants} \label{ss:inv-ops}

Let $V$ be a simple $\ul{G}$-module. We know that $V$ is isomorphic to $L_{\lambda}$ for some weight $\lambda$. By Corollary~\ref{cor:inv-dim}, the length of $\lambda$ is the minimal $n$ for which $V^{G(a)}$ is non-zero for $a \in \bR^{(n)}$. We now show how one can recover $\lambda$ itself via certain operators on this invariant space.

Fix a $\ul{G}$-module $V$ and let $a \in \bR^{(n)}$. We define an operator
\begin{displaymath}
e_{\lambda} \colon V^{G(a)} \to V^{G(a)}
\end{displaymath}
as follows. Let $x \in V^{G(a)}$ be given. Then
\begin{displaymath}
e_{\lambda}(x) = \int_{G/G(a)} \psi_{\lambda}^a(ga)\,gx\ dg.
\end{displaymath}
The element $e_{\lambda}(x)$ is $G(a)$-invariant since the function $\psi_{\lambda}^a$ is $G(a)$-invariant. The following is our main result:

\begin{proposition} \label{prop:inv-op}
The operator $e_{\lambda}$ is idempotent. Its rank is the multiplicity of $L_{\lambda}$ in $V$.
\end{proposition}

\begin{proof}
Since $e_{\lambda}$ is natural in $V$, it suffices to treat the case where $V$ is simple. We show that $e_{\lambda}=0$ if $V$ is not isomorphic to $L_{\lambda}$, and $e_{\lambda}$ is the identity if $V=L_{\lambda}$.

First suppose that $V$ is not isomorphic to $L_{\lambda}$. Let $x \in V^{G(a)}$. There is a unique map of $\ul{G}$-modules $f \colon \cC(\bR^{(n)}) \to V$ satisfying $f(\delta_a)=x$ \cite[Corollary~11.19]{HarmanSnowden}. For $\phi \in \cC(\bR^{(n)})$, we have
\begin{displaymath}
f(\phi)=\int_{G/G(a)} \phi(ga)\,gx\ dg.
\end{displaymath}
Thus $e_{\lambda}(x)=f(\psi^a_{\lambda})$. Since $V$ is not isomorphic to $L_{\lambda}$, we see that $f$ restricts to the zero map $L_{\lambda} \to V$. Since $\psi^a_{\lambda}$ belongs to $L_{\lambda}$, it follows that $e_{\lambda}(x)=0$.

Now suppose that $V=L_{\lambda}$. The space $V^{G(a)}$ is one-dimensional and spanned by $\psi^a_{\lambda}$. It follows that $e_{\lambda}(\psi^a_{\lambda})=\alpha \cdot \psi^a_{\lambda}$ for some scalar $\alpha$, and we must show $\alpha=1$. Since $\psi^a_{\lambda}(a)=1$, we have $\alpha=(e_{\lambda} \psi^a_{\lambda})(a)$. For $g \in G$, we have $g \psi^a_{\lambda}=\psi^{ga}_{\lambda}$. Thus
\begin{displaymath}
\alpha = \int_{G/G(a)} \psi^a_{\lambda}(ga) \psi^{ga}_{\lambda}(a) dg
=\int_{\bR^{(n)}} \psi^a_{\lambda}(b) \psi^b_{\lambda}(a) db
=\int_{\bR^{(n)}} \psi^a_{\lambda}(b) \psi^a_{\lambda^{\vee}}(b) db
=\psi^a_{\lambda}(a)=1,
\end{displaymath}
where in the fourth step we used Proposition~\ref{prop:psi-ev}.
\end{proof}

\section{Induction and restriction} \label{s:ind}

In this section, we study induction and restriction between $G$ and its open subgroups. Our main results give explicit formulas for these operations. As applications, we show that $\rK$ is a Hopf algebra and we classify its primitive elements.

\subsection{Induction}

Recall (\S \ref{ss:groth}) the induction product
\begin{displaymath}
\ind \colon \rK \otimes \rK \to \rK,
\end{displaymath}
which is induced by the induction functor from $G(0) \cong G \times G$ to $G$. We sometimes write $\ind(x,y)$ in place of $\ind(x \otimes y)$ for readability. We now compute this operation explicitly, in terms of the concatenation product $\odot$ (see \S \ref{ss:concat}):

\begin{theorem} \label{thm:ind}
For $x,y \in \rK$, we have
\begin{displaymath}
\ind(x,y) =
x \odot [\cC(\bR)] \odot y,
\end{displaymath}
where $[\cC(\bR)]=a_{\wa}+a_{\wb}+1$.
\end{theorem}

See Example~\ref{ex:decomp} for the formula for $[\cC(\bR)]$. We note two corollaries of the theorem:

\begin{corollary} \label{cor:ind}
For weights $\lambda$ and $\mu$, we have
\begin{displaymath}
\ind(a_{\lambda}, a_{\mu}) = a_{\lambda\wa\mu} + a_{\lambda\wb\mu} + a_{\lambda\mu}
\end{displaymath}
\end{corollary}

\begin{corollary} \label{cor:ind-cat}
For $x,y,z \in \rK$, we have
\begin{displaymath}
\ind(x \odot y, z)=x \odot \ind(y,z), \qquad
\ind(x, y \odot z) = \ind(x,y) \odot z
\end{displaymath}
\end{corollary}

\begin{remark}
The Schwartz space $\cC(\bR^{(n)})$ is the induction to $G$ of the trivial representation of $G(a)$, where $a \in \bR^{(n)}$. The theorem therefore gives $[\cC(\bR^{(n)})]=[\cC(\bR)]^{\odot n}$. This agrees with the decomposition of $\cC(\bR^{(n)})$ given in Theorem~\ref{thm:decomp}.
\end{remark}

\begin{remark} \label{rmk:concat}
We have defined the concatenation product on the Grothendieck group, but not at the level of representations. Corollary~\ref{cor:ind} suggests a possible definition of $L_{\lambda} \odot L_{\mu}$, namely, it can be realized as the $\ul{G}$-submodule of $\Ind(L_{\lambda} \boxtimes L_{\mu})$ generated by the $G(a)$-invariant, where $a \in \bR^{(n)}$ with $n=\ell(\lambda)+\ell(\mu)$.
\end{remark}

We now prove Theorem~\ref{thm:ind}. The proof has essentially two steps. First, we show that $\Ind(L_{\lambda} \boxtimes L_{\mu})$ has length three by computing the dimensions of invariant spaces. Second, we identify the three simple constituents using the operators studied in \S \ref{ss:inv-ops}.

We begin by introducing a convenient device for counting invariants. Let $V$ be a finite length object of $\uRep(G)$. We define the \emph{Hilbert function} $h_V$ of $V$ by letting $h_V(n)$ be the dimension of the $G(a)$-invariants of $V$, for $a \in \bR^{(n)}$; this is well-defined since the various $G(a)$ subgroups are all conjugate. We also define the \emph{Hilbert series} of $V$ by
\begin{displaymath}
\rH_V(t) = \sum_{n \ge 0} h_V(t) t^n.
\end{displaymath}
The following result tells us what these invariants look like.

\begin{lemma} \label{lem:hilbert}
Let $V$ be a finite length object of $\uRep(G)$. Let $m_{\lambda}$ be the multiplicity of $L_{\lambda}$ in $V$, and let $m_n=\sum_{\ell(\lambda)=n} m_{\lambda}$. Then
\begin{displaymath}
\rH_V(t) = \sum_{n \ge 0} m_n \cdot \frac{t^n}{(1-t)^{n+1}}.
\end{displaymath}
\end{lemma}

\begin{proof}
It suffices to treat the case where $V$ is simple; this follows from Corollary~\ref{cor:inv-dim}.
\end{proof}

Write  $\bR=\bR_- \sqcup \{0\} \sqcup \bR_+$ where $\bR_-=(-\infty,0)$ and $\bR_+=(0,\infty)$, and let $G_{\pm}$ be the group of order-preserving automorphisms of $\bR_{\pm}$. These groups are isomorphic to $G$, and we have $G(0)=G_- \times G_+$. We the isomorphism $G_{\pm} \cong G$, we can apply concepts about $G$-representations to $G_{\pm}$-representations, such as Hilbert series.

\begin{lemma} \label{lem:ind-1}
Let $W_1$ and $W_2$ be representations of $\ul{G}_-$ and $\ul{G}_+$, and let $V$ be the induction of $W_1 \boxtimes W_2$ from $G(0)$ to $G$. Then
\begin{displaymath}
\rH_V(t) = (1+t) \rH_{W_1}(t) \rH_{W_2}(t)
\end{displaymath}
\end{lemma}

\begin{proof}
Recall that $h_V(n)$ is the dimension of the invariant space $V^{G(a)}$, where $a \in \bR^{(n)}$. By Frobenius reciprocity, this is equal to the dimension of $\Hom_{\ul{G}}(\cC(\bR^{(n)}), V)$. Applying Frobenius reciprocity in the other argument, we find
\begin{displaymath}
h_V(n) = \dim \Hom_{\ul{G}(0)}(\cC(\bR^{(n)}), W_1 \boxtimes W_2).
\end{displaymath}
Now, we have
\begin{displaymath}
\bR^{(n)} = \bigsqcup_{r+s+t=n, s=0,1} \bR_-^{(r)} \times \bR_+^{(t)},
\end{displaymath}
and so
\begin{displaymath}
\cC(\bR^{(n)}) = \bigoplus_{r+s+t=n, s=0,1} \cC(\bR_-^{(r)}) \boxtimes \cC(\bR_+^{(t)}).
\end{displaymath}
Thus we have
\begin{displaymath}
h_V(n) = \sum_{r+t=n} h_{W_1}(r) h_{W_2}(t) + \sum_{r+t=n-1} h_{W_1}(r) h_{W_2}(t),
\end{displaymath}
where the second term is omitted if $n=0$. This yields the stated formula.
\end{proof}

\begin{lemma} \label{lem:ind-2}
Let $W_1$, $W_2$, and $V$ be as in the previous lemma, and suppose that $W_1$ and $W_2$ are simple. Then $V$ has length three.
\end{lemma}

\begin{proof}
Suppose $W_1$ and $W_2$ are associated to weights of length $n$ and $m$. By Lemmas~\ref{lem:hilbert} and~\ref{lem:ind-1}, we have
\begin{align*}
\rH_V(t) = (1+t) \rH_{W_1}(t) \rH_{W_2}(t)
&= (1+t) \cdot \frac{t^n}{(1-t)^{n+1}} \cdot \frac{t^m}{(1-t)^{m+1}} \\
&= ((1-t)+2t) \cdot \frac{t^{n+m}}{(1-t)^{n+m+2}} \\
&= \frac{t^{n+m}}{(1-t)^{n+m+1}} + 2 \frac{t^{n+m+1}}{(1-t)^{n+m+2}}.
\end{align*}
By Lemma~\ref{lem:hilbert} again, we thus see that $V$ has three simples, with one of length $n+m$ and two of length $n+m+1$.
\end{proof}

\begin{proof}[Proof of Theorem~\ref{thm:ind}]
Write $L^+_{\lambda}$ for the $G_+$ version of $L_{\lambda}$. This is the $\ul{G}_+$-submodule of $\cC(\bR_+^{(n)})$ generated by the $\phi_{\bI}$'s, where $\bI=(I_1<\cdots<I_n)$ are intervals of type $\lambda$ contains in $\bR_+$ (here $n$ is the length of $\lambda$). We similarly have $L^-_{\lambda}$

Let $W_1=L^-_{\lambda}$ and $W_2=L^+_{\mu}$, and let $V$ be the induction of $W_1 \boxtimes W_2$. We show that
\begin{displaymath}
V = L_{\lambda \mu} \oplus L_{\lambda\wa\mu} \oplus L_{\lambda\wb\mu},
\end{displaymath}
and this will complete the proof. Since we know $V$ has length three by Lemma~\ref{lem:ind-2}, it is enough to show that each simple on the right appears in $V$. By Frobenius reciprocity, it is enough to show that $L_{\lambda}^- \boxtimes L_{\lambda}^+$ appears in the restriction to $G(0)$ of each of these simples.

Let $n=\ell(\lambda)$ and $m=\ell(\mu)$. We have natural maps
\begin{displaymath}
L_{\lambda}^- \boxtimes L_{\mu}^+ \to \cC(\bR_-^{(n)}) \boxtimes \cC(\bR_+^{(m)}) \to \Res^G_{G(0)} \cC(\bR^{(n+m)}),
\end{displaymath}
where the first comes from the inclusions, and the second takes the product of the pullbacks of the two functions. This map carries $\phi_{\bI} \boxtimes \phi_{\bJ}$ to $\phi_{(\bI, \bJ)}$, where $(\bI, \bJ)$ means $(I_1, \ldots, I_n, J_1, \ldots, J_n)$. Since $\phi_{(\bI, \bJ)}$ belongs to $L_{\lambda \mu}$, we thus have a non-zero $G(0)$-map $L^-_{\lambda} \boxtimes L^+_{\mu} \to L_{\lambda \mu}$, as desired.

Let $\nu=\lambda\wa\mu$, and fix $a \in \bR_-^{(n)}$ and $b \in \bR_+^{(m)}$. Let $e^-_{\lambda}$ be the operator on the $G(a)$-invariants of a $\ul{G}_-$-representation defined in \S \ref{ss:inv-ops}, and similarly define $e_{\mu}^+$. To show that $L^-_{\lambda} \boxtimes L^+_{\mu}$ appears in $\Res^G_{G(0)}(L_{\nu})$, it is enough (by Proposition~\ref{prop:inv-op}) to show that $e_{\lambda}^- e_{\mu}^+$ is non-zero on the $G_-(a) \times G_+(b)$ invariants. Note that $G_-(a) \times G_+(b)=G(c)$, where $c=(a,0,b) \in \bR^{(n+m+1)}$. We must therefore show that $e_{\lambda}^- e_{\mu}^+$ is non-zero on $L_{\nu}^{G(c)}$. This space is one-dimensional, and spanned by $\psi^c_{\nu}$ (Proposition~\ref{prop:inv}). Thus it is enough to show that $e_{\lambda}^-$ and $e_{\mu}^+$ are separately non-zero on $\psi^c_{\nu}$.

The argument here is similar to the final paragraph of the proof of Proposition~\ref{prop:inv-op}. We have $e_{\lambda}^-(\psi^c_{\nu}) = \alpha \cdot \psi^c_{\nu}$ for some scalar $\alpha$, and we must show $\alpha \ne 0$. We have
\begin{displaymath}
\alpha = \int_{\bR_-^{(n)}} \psi^a_{\lambda}(p) \cdot \psi_{\nu}^{(p,0,b)}(c) dp,
\end{displaymath}
and so $\alpha=\vol(X)$, where
\begin{displaymath}
X=\{p \in \bR_-^{(n)} \mid \text{$p \in \Psi^a_{\lambda}$ and $(a,0,b) \in \Psi^{(p,0,b)}_{\nu}$} \}.
\end{displaymath}
One finds that $X=\{a\}$, and so $\alpha=1$. The case of $e_{\mu}^+$ is similar.

We have thus shown that $L^-_{\lambda} \boxtimes L^+_{\mu}$ appears in the restriction of $L_{\lambda\wa\nu}$. The case of $L_{\lambda\wb\nu}$ is similar.
\end{proof}

\subsection{Restriction} \label{ss:res}

Recall (\S \ref{ss:groth}) that we also have a map
\begin{displaymath}
\res \colon \rK \to \rK \otimes \rK
\end{displaymath}
given by restricting from $G$ to $G \times G$. We now determine this map explicitly. We first introduce a piece of notation: for a weight $\lambda=\lambda_1 \cdots \lambda_n$ of length $n$, we let $\lambda[i,j]$ denote the substring of $\lambda$ between indices $i$ and $j$ (inclusively), i.e., $\lambda_i \cdots \lambda_j$. As with intervals, we use parentheses to exclude the edge values, e.g., $\lambda[i,j)=\lambda_i \cdots \lambda_{j-1}$.

\begin{theorem} \label{thm:res}
Let $\lambda$ be a word of length $n$. Then
\begin{displaymath}
\res(a_{\lambda}) = \sum_{i=0}^n a_{\lambda[1,i]} \otimes a_{\lambda(i,n]} + \sum_{i=1}^n a_{\lambda[1,i)} \otimes a_{\lambda(i,n]}.
\end{displaymath}
In words, the first sum is over all ways of breaking $\lambda$ into two between letters (or at the ends), while the second sum is over all ways of breaking $\lambda$ into two by deleting letters.
\end{theorem}

\begin{proof}
Let $\alpha$ and $\beta$ be weights, and let $i=\ell(\alpha)$. We have
\begin{align*}
\langle \res(a_{\lambda}), a_{\alpha} \otimes a_{\beta} \rangle
&= \langle a_{\lambda}, \ind(a_{\alpha} \otimes a_{\beta}) \rangle \\
&= \langle a_{\lambda}, a_{\alpha\wa\beta} \rangle + \langle a_{\lambda}, a_{\alpha\wb\beta} \rangle + \langle a_{\lambda}, a_{\alpha\beta} \rangle
\end{align*}
where in the first step we used Frobenius reciprocity (see \S \ref{ss:groth}), and in the second Corollary~\ref{cor:ind}. The first two terms above, counted together, give~1 if $\alpha=\lambda[1,i]$ and $\beta=\lambda[i+2,n]$ and $0 \le i \le n-1$, and~0 otherwise. The third term above gives~1 if $\alpha=\lambda[1,i]$ and $\beta=\lambda[i+1,n]$, and~0 otherwise. This yields the stated formula.
\end{proof}

The theorem gives us the following very important corollary for how restriction interacts with the filtration on $\rK$.

\begin{corollary} \label{cor:res}
Let $x \in \rK_{\le n}$. Then
\begin{displaymath}
\res(x)=x \otimes 1 + 1 \otimes x + y
\end{displaymath}
where $y$ belongs to $\sum_{i+j=n, i,j>0} \rK_{\le i} \otimes \rK_{\le j}$.
\end{corollary}

The point here is that $y$ is built out of elements in smaller filtration degree; thus $\res(x)$ has a very simple form up to an error that can be controlled. This observation is used in several proofs in the remainder of the paper; the following example illustrates the basic mechanism used in these arguments.

\begin{example} \label{ex:res-dim}
Theorem~\ref{thm:res} gives an alternative approach to computing the categorical dimension of simple objects. Let $\lambda$ of length $n$ be given, and suppose we know $a_{\mu}$ has dimension $(-1)^{\ell(\mu)}$ for $\ell(\mu)<n$. From the theorem, we have
\begin{displaymath}
\res(a_{\lambda}) = (a_{\lambda} \otimes 1 + 1 \otimes a_{\lambda}) + \sum_{i=1}^{n-1} a_{\lambda[1,i]} \otimes a_{\lambda(i,n]} + \sum_{i=1}^n a_{\lambda[1,i)} \otimes a_{\lambda(i,n]}.
\end{displaymath}
Categorical dimension is invariant under restriction, and multiplicative for the external tensor product. Thus if $d$ is the dimension of $a_{\lambda}$ then the above equation gives
\begin{displaymath}
d=2d+(n-1) (-1)^n + n (-1)^{n-1},
\end{displaymath}
and so $d=(-1)^n$.
\end{example}

\subsection{The Hopf algebra structure on $\rK$} \label{ss:hopf}

Recall (Proposition~\ref{prop:delta}) that we have an algebra homomorphism $\delta \colon \rK \to \bZ$ defined by $\delta(a_{\lambda})=(-1)^{\ell(\lambda)}$.

\begin{proposition} \label{prop:counit}
The map $\delta$ is a co-unit for $\res$.
\end{proposition}

\begin{proof}
We must show that the composition
\begin{displaymath}
\xymatrix@C=3em{
\rK \ar[r]^-{\res} & \rK \otimes \rK \ar[r]^-{\delta \otimes 1} & \rK }
\end{displaymath}
is the identity, and similarly for $1 \otimes \delta$. Let $\lambda$ be a weight of length $n$. Applying Theorem~\ref{thm:res}, we have
\begin{align*}
(\delta \otimes 1)(\res(a_{\lambda}))
&= \sum_{i=0}^n \delta(a_{\lambda[1,i]}) a_{\lambda(i,n]} + \sum_{i=1}^n \delta(a_{\lambda[1,i)}) a_{\lambda(i,n]} \\
&= a_{\lambda} + \sum_{i=1}^n \big( \delta(a_{\lambda[1,i]}) + \delta(a_{\lambda[1,i)}) \big) a_{\lambda(i,n]}.
\end{align*}
The $a_{\lambda}$ here is the $i=0$ term of the first sum. Since the lengths of $\lambda[1,i]$ and $\lambda[1,i)$ differ by one, the $\delta$'s in the second line cancel, and so the final result is simply $a_{\lambda}$. This shows that $(\delta \otimes 1) \circ \res$ is the identity. The $1 \otimes \delta$ case is similar.
\end{proof}

We now know that $\rK$ is a bi-algebra, under the standard multiplication $m$ and the co-multiplication given by $\res$; the unit is~1 and the co-unit is $\delta$. In fact:

\begin{proposition} \label{prop:hopf}
The bi-algebra $\rK$ is a Hopf algebra.
\end{proposition}

\begin{proof}
Recall that if $f,g \colon \rK \to \rK$ are additive maps, their \defn{convolution} $f \ast g$ is the composition
\begin{displaymath}
\xymatrix@C=3em{
\rK \ar[r]^-{\res} & \rK \otimes \rK \ar[r]^-{f \otimes g} & \rK \otimes \rK \ar[r]^-m & \rK }
\end{displaymath}
Convolution gives the set $\Hom(\bK,\bK)$ of all additive maps the structure of a unital and associative algebra; see \cite[\S 4.2]{hopf}. The unit is the composition $\eta \delta$, where $\eta \colon \bZ \to \rK$ is the unit map. An \defn{antipode} is an additive map $S \colon \rK \to \rK$ that is a two-sided inverse to the identity map under convolution, that is, we have $S \ast \id = \id \ast S = \eta \delta$. To prove the proposition, we must show that $\rK$ admits an antipode.

Consider the equation $\id \ast S=\eta \delta$ on a map $S$. Let $\lambda$ be a word of length $n$. Evaluating at $a_{\lambda}$, this takes the form 
\begin{displaymath}
m((\id \otimes S)(\res(a_{\lambda})))=(-1)^n.
\end{displaymath}
Appealing to our computation of $\res$ (Theorem~\ref{thm:res}), this becomes
\begin{displaymath}
\sum_{i=0}^n a_{\lambda[1,i]} S(a_{\lambda(i,n]}) + \sum_{i=1}^n a_{[1,i)} S(a_{\lambda(i,n]}) = (-1)^n.
\end{displaymath}
The $i=0$ term in the first sum is $S(a_{\lambda})$. Every other input to $S$ in the above equation has the form $a_{\mu}$ where $\ell(\mu)<n$. Thus the above equation allows us to recursively solve for $S(a_{\lambda})$.

The above reasoning shows that there is a unique right inverse to $\id$ under convolution. A similar argument shows that there is a unique left inverse. Sinve convolution is associative, these two inverses are equal and are thus the antipode.
\end{proof}

\begin{example}
Proposition~\ref{prop:hopf} shows how to recursively compute the anitpode $S$. Here are a few examples:
\begin{align*}
S(a_{\wa}) &= -a_{\wa}-2 \\
S(a_{\wa\wa}) &= a_{\wa\wa} + 3 a_{\wa} +3 \\
S(a_{\wa\wb}) &= a_{\wb\wa}+2a_{\wa}+2a_{\wb}+4.
\end{align*}
We do not know the general formula for $S(a_{\lambda})$. However, the leading term is always $(-1)^{\ell(\lambda)} a_{\mu}$, where $\mu$ is the reverse of $\lambda$; see Remark~\ref{rmk:gr-hopf}.
\end{example}

\subsection{Classification of primitive elements} \label{ss:prim}

An element $x \in \rK$ is \defn{primitive} if
\begin{displaymath}
\res(x)=(x \otimes 1)+(1 \otimes x).
\end{displaymath}
The collection of primitive elements of $\rK$ forms an additive subgroup, which we denote by $\rK^{\rm prim}$. We now determine it.

\begin{proposition} \label{prop:prim}
The group $\rK^{\rm prim}$ has a $\bZ$-basis consisting of $a_{\wa}+1$ and $a_{\wb}+1$.
\end{proposition}

\begin{proof}
From Theorem~\ref{thm:res}, we have
\begin{displaymath}
\res(a_{\wa}+1) = (a_{\wa}\otimes 1)+(1 \otimes a_{\wa})+2;
\end{displaymath}
note that $\res(1)=1$. We also have
\begin{displaymath}
(a_{\wa}+1) \otimes 1 + 1 \otimes (a_{\wa}+1) = (a_{\wa}\otimes 1)+(1 \otimes a_{\wa})+2.
\end{displaymath}
Thus $a_{\wa}+1$ is primitive, and similarly for $a_{\wb}+1$.

Let $\gamma \colon \rK \otimes \rK \to \rK$ be the concatentation product, i.e., $\gamma(a_{\lambda} \otimes a_{\mu})=a_{\lambda \mu}$. From Theorem~\ref{thm:res}, we have $\gamma(\res(a_{\lambda})) = (n+1)a_{\lambda}+y$ where $y$ belongs to $\rK_{\le n-1}$. It follows that for any $x \in \rK_{\le n}$, we have $\gamma(\res(x)) =(n+1) x+y$ for some $y \in \rK_{\le n-1}$. On the other hand, if $x$ is primitive then clearly $\gamma(\res(x))=2x$. We thus see that $x$ can be primitive only if $x$ belongs to $\rK_{\le 1}$; in other words, $\rK^{\rm prim} \subset \rK_{\le 1}$.

Now, the $\bZ$-module $\rK_{\le 1}$ is free of rank three. Since no multiple of~1 is primitive, it follows that $\rK^{\rm prim}$ has rank $\le 2$. But we have already exhibited two independent elements, so $\rK^{\rm prim}$ has rank~2. Since these two elements generate a saturated $\bZ$-submodule of $\rK_{\le 1}$, they must generate all of $\rK^{\rm prim}$.
\end{proof}

Let $\Delta \colon \rK \to \rK \otimes \rK$ be the map given by
\begin{displaymath}
\Delta(x) = \res(x) - x \otimes 1 - 1 \otimes x.
\end{displaymath}
Thus $x$ is primitive if and only if $\Delta(x)=0$. Since there are very few primitive elements, $\Delta$ is nearly injective. On the other hand, $\Delta$ maps $\rK_{\le n}$ into $\rK_{\le n-1} \otimes \rK_{\le n-1}$ by Corollary~\ref{cor:res}. This provides a very useful set-up for inductive arguments, which we employ a few times.

\section{Tensor products} \label{s:tensor}

In this section, we give a combinatorial rule for the decomposition of a tensor product of two simple representations. As an application, we show that $\rK \otimes \bQ$ is a polynomial ring with variables indexed by Lydon words (Corollary~\ref{cor:K-poly}).

\subsection{The tensor product formula}

Let $\lambda$ and $\mu$ be two weights of lengths $n$ and $m$. A \defn{ruffle}, or \defn{$(n,m)$-ruffle}, is a surjective function let $r \colon [n] \amalg [m] \to [\ell]$ that is injective and order-preserving on $[n]$ and on $[m]$. We think of a ruffle like a shuffle except that there can be collisions. Fix a ruffle $r$, and let $r_1$ and $r_2$ be the restrictions of $r$ to $[n]$ and $[m]$. We say there is a \defn{collision} at $i \in [\ell]$ if $i$ belongs to the images of $r_1$ and $r_2$. Suppose $i=r_1(j)=r_2(\ell)$ is a collision. We say that the collision is \defn{positive} if $\lambda_j=\mu_{\ell}=\wa$, \defn{negative} if $\lambda_j=\mu_{\ell}=\wb$, and \defn{neutral} otherwise. We define the \defn{pre-product} of $\lambda$ and $\mu$ with respect to $r$ to be the length $\ell$ word defined as follows. At the $i$th spot, we put:
\begin{itemize}
\item the letter $\lambda_j$ if $i=r_1(j)$ is not a collision.
\item the letter $\mu_{\ell}$ if $i=r_2(\ell)$ is not a collision.
\item the letter $\wa$ if $i$ is a positive collision
\item the letter $\wb$ if $i$ is a negative collision.
\item the letter $?$ if $i$ is a neutral collision.
\end{itemize}
We now define a multiset $P_r(\lambda, \mu)$ by replacing each $?$ in the pre-product with either $\wa$, $\wb$, or the empty string (in which case the ? character is simply deleted), in all possible ways. Thus if there are $q$ question marks then $P_r(\lambda,\mu)$ has cardinality $3^q$ (as a multiset). Finally, we define
\begin{displaymath}
m'(a_{\lambda}, a_{\mu}) = \sum_r \sum_{\nu \in P_r(\lambda,\mu)} a_{\nu},
\end{displaymath}
where the outer sum is over all $(n,m)$-ruffles. We extend $m'$ bilinearly to a product on $\rK$.

\begin{example} \label{ex:tensor2}
Consider the two weights $\lambda=\wa$ and $\mu=\wa\wb$. There are five ruffles. The corresponding pre-products are
\begin{displaymath}
\ul{\wa}\wa\wb, \quad \wa\ul{\wa}\wb, \quad \wa\wb\ul{\wa}, \quad \ul{\wa}\wb, \quad \wa \ul{?}
\end{displaymath}
We have underlined the position that $\lambda$ occupies in the above words. In the first three cases above there are no collisions (i.e., they are shuffles); in the fourth case there is a positive collision; and in the fifth case there is a neutral collision. In the first four cases, $P_r(\lambda,\mu)$ contains just the indicated word; in the final case, $P_r(\lambda,\mu)=\{\wa,\wa\wa,\wa\wb\}$. Thus
\begin{displaymath}
m'(a_{\wa}, a_{\wa\wb})=2 a_{\wa\wa\wb}+a_{\wa\wb\wa}+2a_{\wa\wb}+a_{\wa\wa}+a_{\wa}. \qedhere
\end{displaymath}
\end{example}

Throughout this section, we write $m$ for the standard product on $\rK$, i.e., the one induced by tensor products. The following is our main result on tensor products:

\begin{theorem} \label{thm:tensor}
We have $m=m'$.
\end{theorem}

The theorem gives the following combinatorial rule for computing tensor products.

\begin{corollary} \label{cor:tensor}
Let $\lambda$ and $\mu$ be weights of lengths $n$ and $m$. Then
\begin{displaymath}
L_{\lambda} \uotimes L_{\mu} \cong \bigoplus_r \bigoplus_{\nu \in P_r(\lambda,\mu)} L_{\nu},
\end{displaymath}
where the outer sum is over all $(n,m)$-ruffles $r$.
\end{corollary}

\begin{example} \label{ex:tensor3}
Combining Example~\ref{ex:tensor2} and Corollary~\ref{cor:tensor}, we find
\begin{displaymath}
L_{\wa} \uotimes L_{\wa\wb} \cong L_{\wa\wa\wb}^{\oplus 2} \oplus L_{\wa\wb\wa} \oplus L_{\wa\wb}^{\oplus 2} \oplus L_{\wa\wa} \oplus L_{\wa}. \qedhere
\end{displaymath}
\end{example}

To prove the theorem, we compute one special type of product for $m$ (Proposition~\ref{prop:prod-special}) and establish the projection formula for $m'$ (Proposition~\ref{prop:comb-proj}); this turns out to be enough to deduce that $m=m'$ (see Proposition~\ref{prop:prod-char}).

\subsection{Some simple tensor products} \label{ss:simple-tensor}

We now compute a family of particularly simple products. Let $\pi(n)$ be a string of length $n$ consisting of all $\wa$'s.

\begin{proposition} \label{prop:prod-special}
For $n \ge 0$, we have
\begin{displaymath}
a_{\pi(1)} \cdot a_{\pi(n)}=(n+1) a_{\pi(n+1)} + n a_{\pi(n)}.
\end{displaymath}
\end{proposition}

We first give a non-rigorous but illuminative argument. The tensor product $L_{\pi(1)} \uotimes L_{\pi(n)}$ is the subrepresentation of $\cC(\bR \times \bR^{(n)})$ generated by the functions $\phi_{(J,I_1,\ldots,I_n)}$ where $J$ and each $I_i$ has type $\wa$, and $I_1<\cdots<I_n$. The functions with $J=I_i$ generate a copy of $L_{\pi(n)}$, while the functions with $I_i<J<I_{i+1}$ (or $J<I_1$ or $I_n<J$) generate a copy of $L_{\pi(n+1)}$. This gives the stated decomposition.

One can make the above argument rigorous; however, there are many details to handle. Instead, we give a different proof simply by computing in the Grothendieck group, which avoids any subtleties. This proof is perhaps less insightful, but demonstrates how rigid the various structures on $\rK$ are.

Put $x_n=a_{\pi(1)} a_{\pi(n)}$ and $y_n=(n+1) a_{\pi(n+1)}+na_{\pi(n)}$. Thus we must show $x_n=y_n$. We will eventually do this by induction. We prove two lemmas first. Recall from \S \ref{ss:prim} that $\Delta(x)=\res(x)-x \otimes 1 - 1 \otimes x$.

\begin{lemma} \label{lem:special-1}
Let $n>0$ be given, and suppose $x_i=y_i$ for $0 \le i < n$. Then $\Delta(x_n)=\Delta(y_n)$.
\end{lemma}

\begin{proof}
Put
\begin{displaymath}
\Sigma_n = \sum_{i+j=n} a_{\pi(i)} \otimes a_{\pi(j)},
\end{displaymath}
so that $\res(a_{\pi(n)}) = \Sigma_n + \Sigma_{n-1}$ by Theorem~\ref{thm:res}. Since $\res$ is a ring homomorphism, we have
\begin{displaymath}
\res(x_n)
= (\Sigma_1 + 1) \cdot (\Sigma_n+\Sigma_{n-1})
\end{displaymath}
We have
\begin{align*}
(a_{\pi(1)} \otimes 1) \cdot \Sigma_n
&= \sum_{i+j=n} x_i \otimes a_{\pi(j)} \\
&= (x_n-y_n) \otimes 1 + \sum_{i+j=n} ((i+1) a_{\pi(i+1)}+i a_{\pi(i)}) \otimes a_{\pi(j)} \\
&= (x_n-y_n) \otimes 1 + \sum_{i+j=n+1} ia_{\pi(i)} \otimes a_{\pi(j)}
+ \sum_{i+j=n} ia_{\pi(i)} \otimes a_{\pi(j)}
\end{align*}
There is a similar formula for $(1 \otimes a_{\pi(1)}) \cdot \Sigma_n$. Summing these, we find
\begin{displaymath}
\Sigma_1 \cdot \Sigma_n \\
= (x_n-y_n) \otimes 1 + 1 \otimes (x_n-y_n) + (n+1) \Sigma_{n+1} + n \Sigma_n
\end{displaymath}
A similar computation gives
\begin{displaymath}
\Sigma_1 \cdot \Sigma_{n-1} = n \Sigma_n + (n-1) \Sigma_{n-1}
\end{displaymath}
Combining all the above, we find
\begin{displaymath}
\res(x_n) = (x_n-y_n) \otimes 1 + 1 \otimes (x_n-y_n) + (n+1) \Sigma_{n+1} + (2n+1)  \Sigma_n + n \Sigma_{n-1}.
\end{displaymath}
We also have
\begin{displaymath}
\res(y_n) = (n+1) (\Sigma_{n+1}+\Sigma_n) + n(\Sigma_n+\Sigma_{n-1}) = (n+1) \Sigma_{n+1} + (2n+1) \Sigma_n + n \Sigma_{n-1}.
\end{displaymath}
Thus
\begin{displaymath}
\res(x_n)-\res(y_n) =  (x_n-y_n) \otimes 1 + 1 \otimes (x_n-y_n),
\end{displaymath}
and so $\Delta(x_n)=\Delta(y_n)$, as required.
\end{proof}

\begin{lemma} \label{lem:special-2}
We have
\begin{displaymath}
a_{\wa} a_{\wa}=2a_{\wa\wa}+a_{\wa}, \qquad
a_{\wa} a_{\wb} = a_{\wa\wb}+a_{\wb\wa}+a_{\wa}+a_{\wb}+1, \qquad
a_{\wb} a_{\wb} = 2a_{\wb\wb} + a_{\wb}
\end{displaymath}
\end{lemma}

\begin{proof}
From Lemma~\ref{lem:special-1}, we have $\Delta(x_1)=\Delta(y_1)$. It follows that $x_1-y_1$ belongs to the kernel of $\Delta$, and is therefore primitive. Thus by Proposition~\ref{prop:prim}, we have
\begin{displaymath}
a_{\wa}a_{\wa}=2a_{\wa\wa}+a_{\wa}+p(a_{\wa}+1)+q(a_{\wb}+1)
\end{displaymath}
for integers $p$ and $q$. Since this is an effective class, we have $p \ge -1$ and $q \ge 0$. We also have
\begin{displaymath}
p+q=\langle a_{\wa} a_{\wa}, 1 \rangle = \langle a_{\wa}, a_{\wb} \rangle=0,
\end{displaymath}
where in the second step we used adjunction. We thus have
\begin{displaymath}
a_{\wa} a_{\wa}=2a_{\wa\wa}+(1-q) a_{\wa} + q a_{\wb}
\end{displaymath}
and $q \in \{0,1\}$. Next,
\begin{align*}
a_{\wa} \cdot (a_{\wa}+a_{\wb}+1)
&= a_{\wa} \cdot \ind(1)= \ind(\res(a_{\wa})) \\
&= \ind(a_{\wa} \otimes 1+1 \otimes a_{\wa}+1) \\
&= 2a_{\wa\wa}+a_{\wa\wb} + a_{\wb\wa}+3a_{\wa}+a_{\wb}+1
\end{align*}
where we have used Theorems~\ref{thm:ind} and~\ref{thm:res}, and the projection formula. Thus
\begin{displaymath}
a_{\wa} a_{\wb} = a_{\wa\wb} + a_{\wb\wa} + (1+q) a_{\wa} + (1-q) a_{\wb} + 1
\end{displaymath}
Since the above element is obviously self-dual, we must have $q=0$. This yields the first two formulas, while the third follows from the first by duality.
\end{proof}

\begin{proof}[Proof of Proposition~\ref{prop:prod-special}]
We must show that $x_n=y_n$ for all $n$. The $n=0$ case is trivial, while the $n=1$ case follows from Lemma~\ref{lem:special-2}. Now let $n \ge 2$ be given, and suppose $x_i=y_i$ for $0 \le i <n$. By Lemma~\ref{lem:special-1}, we have $\Delta(x_n)=\Delta(y_n)$. Thus $x_n-y_n$ is primitive, and so by Proposition~\ref{prop:prim} we have
\begin{displaymath}
x_n=y_n+p(a_{\wa}+1)+q(a_{\wb}+1)
\end{displaymath}
for integers $p$ and $q$. We have
\begin{displaymath}
p+q=\langle x_n, 1 \rangle = \langle a_{\pi(n)}, a_{\wb} \rangle = 0
\end{displaymath}
where in the second step we used adjunction. We also have
\begin{displaymath}
q=\langle x_n, a_{\wb} \rangle = \langle a_{\pi(n)}, a_{\wb} a_{\wb} \rangle
\end{displaymath}
where again we have used adjunction. From the computation of $a_{\wb} a_{\wb}$ in Lemma~\ref{lem:special-2}, the above pairing vanishes. Thus $p=q=0$, and so $x_n=y_n$ as required.
\end{proof}

Let $\bZ\langle x \rangle$ be the ring of integer-valued polynomials in the variable $x$. Put $b_n=\binom{x}{n}$. The elements $\{b_n\}_{n \ge 0}$ form a $\bZ$-basis of $\bZ\langle x \rangle$.

\begin{corollary} \label{cor:integer-valued}
The additive map $\phi \colon \bZ\langle x \rangle \to \rK$ defined by $\phi(b_n)=a_{\pi(n)}$ is a ring homomorphism.
\end{corollary}

\begin{proof}
We have $b_1 b_n=(n+1)b_{n+1}+nb_n$. Since $b_1$ generates $\bZ\langle x \rangle \otimes \bQ$ as a $\bQ$-algebra, it follows that $\phi \otimes \bQ$ is a ring homomorphism. Since $\rK$ is torsion-free, it follows that $\phi$ is a ring homomorphism.
\end{proof}

\subsection{Characterization of $m$}

Let
\begin{displaymath}
q \colon \rK \times \rK \to \rK
\end{displaymath}
be a bilinear map, which we view as a product on $\rK$. We write $q$ still for the induced product on $\rK \otimes \rK$. We consider the following conditions on $\mu$:
\begin{enumerate}
\item The product $q$ is commutative and unital with unit $a_{\emptyset}$.
\item The projection formula holds: for $x,y,z \in \rK$, we have
\begin{displaymath}
\ind(q(\res(x), y \otimes z)) = q(x, \ind(y \otimes z)).
\end{displaymath}
\item For any $n \ge 0$, we have
\begin{displaymath}
q(a_{\pi(1)}, a_{\pi(n)}) = (n+1) a_{\pi(n+1)} + n a_{\pi(n)}.
\end{displaymath}
\end{enumerate}
We characterize the standard product $m$ via these conditions:

\begin{proposition} \label{prop:prod-char}
The standard product is the unique product satisfying (a)--(c) above.
\end{proposition}

We will first require a lemma. Put
\begin{displaymath}
\rI_{\le n} = \sum_{i+j=n-1} \ind(\rK_{\le i} \otimes \rK_{\le j}).
\end{displaymath}
Note that $\rI_{\le n}$ is contained in $\rK_{\le n}$ by Theorem~\ref{thm:ind}. We define the \defn{Hamming distance} between two weights of the same length to be the number of positions at which they differ.

\begin{lemma} \label{lem:hamming}
Let $\lambda$ and $\mu$ be two weights of length $n$ with Hamming distance $d$. Then in the quotient
\begin{displaymath}
\rK_{\le n}/(\rK_{\le n-1}+\rI_{\le n})
\end{displaymath}
we have $a_{\lambda}=(-1)^d a_{\mu}$. In particular, the above group is cyclic and generated by $a_{\pi(n)}$.
\end{lemma}

\begin{proof}
It suffices to treat the case $d=1$. Thus, after swapping $\lambda$ and $\mu$ if necessary, we have $\lambda=\alpha\wa\beta$ and $\mu=\alpha\wb\beta$ for two weights $\alpha$ and $\beta$ whose lengths sum to $n-1$. In the quotient group we have
\begin{displaymath}
\ind(a_{\alpha} \otimes a_{\beta}) = a_{\lambda}+a_{\mu}+a_{\alpha \beta} = 0,
\end{displaymath}
and also $a_{\alpha \beta}=0$, and so the result follows.
\end{proof}

\begin{proof}[Proof of Proposition~\ref{prop:prod-char}]
We have seen in \S \ref{ss:groth} that $m$ satisfies (a) and (b), while (c) follows is given by Proposition~\ref{prop:prod-special}. Suppose now that $q$ is an arbitrary product satisfying (a)--(c). We must show $q=m$. Consider the following statement:
\begin{itemize}[leftmargin=4em]
\item[$(S_n)$:] If $x \in \rK_{\le r}$ and $y \in \rK_{\le s}$ with $r+s \le n$ then $q(x, y)=m(x, y)$.
\end{itemize}
We prove $(S_n)$ by induction on $n$. The statement $(S_0)$ follows from (a), since $\rK_{\le 0}$ is spanned by~1. Suppose now that $(S_{n-1})$ holds, and let us prove $(S_n)$.

Let $x$ and $y$ as in $(S_n)$ be given. We must show
\begin{equation} \label{eq:goal}
q(x, y)=m(x, y)
\end{equation}
If $y$ belongs to $\rK_{\le s-1}$ then \eqref{eq:goal} holds by $(S_{n-1})$. Suppose now that $y=\ind(a_{\alpha} \otimes a_{\beta})$ where $\ell(\alpha)+\ell(\beta) \le s-1$. We have
\begin{displaymath}
q(x, y)=\ind(q(\res(x), a_{\alpha} \otimes a_{\beta})),
\end{displaymath}
and similarly for $m$; this is where (b) is used. Note that $\res(x)$ is contained in $(\rK \otimes \rK)_{\le r}$ by Theorem~\ref{thm:res}. By $(S_{n-1})$, we thus have
\begin{displaymath}
q(\res(x), a_{\alpha} \otimes a_{\beta}) = m(\res(x), a_{\alpha} \otimes a_{\beta}).
\end{displaymath}
Applying $\ind$, we find \eqref{eq:goal} holds. We thus see that \eqref{eq:goal} holds if $y$ belongs to $\rK_{\le s-1} + \rI_{\le s}$; of course, there is an analogous statement for $x$ (since $q$ is commutative).

Applying Lemma~\ref{lem:hamming}, write
\begin{displaymath}
x=u a_{\pi(r)}+x', \qquad y=va_{\pi(s)}+y'
\end{displaymath}
where $u,v \in \bZ$ and $x' \in \rK_{\le r-1}+\rI_{\le r}$ and $y' \in \rK_{\le s-1}+\rI_{\le s}$. We have
\begin{displaymath}
q(x,y) = uv q(a_{\pi(i)}, a_{\pi(j)}) + u q(a_{\pi(i)}, y') + v q(x', a_{\pi(j)}) + q(x',y').
\end{displaymath}
Call the four terms above $q_1, \ldots, q_4$. We have a similar expression for $m$, yielding terms $m_1, \ldots, m_4$. We have $q_1=m_1$ by an easy inductive argument, or, more directly, by the reasoning in Corollary~\ref{cor:integer-valued}; in any case, this is where assumption (c) is used. In each of the remaining terms, we have $q_i=m_i$ by the previous paragraph. Thus \eqref{eq:goal} holds, which completes the proof.
\end{proof}

\begin{example} \label{ex:tensor}
The proof of Proposition~\ref{prop:prod-char} explains how to actually compute arbitrary products using the projection formula and Proposition~\ref{prop:prod-special}. We illustrate the simplest case by computing $a_{\wa} \cdot a_{\wb}$. First we express $a_{\wb}$ in terms of inductions and $a_{\pi(n)}$'s. By Theorem~\ref{thm:ind}, we have
\begin{displaymath}
a_{\wb}=\ind(1)-a_{\wa}-1.
\end{displaymath}
Multiplying by $a_{\wa}$, we find
\begin{displaymath}
a_{\wa} \cdot a_{\wb} = \ind(\res(a_{\wa}))-a_{\wa} \cdot a_{\wa}- a_{\wa}.
\end{displaymath}
We have $a_{\wa} \cdot a_{\wa} = 2 a_{\wa\wa} + a_{\wa}$ by Proposition~\ref{prop:prod-special}. We have
\begin{displaymath}
\res(a_{\wa})=(a_{\wa} \otimes 1) + (1 \otimes a_{\wa}) + 1,
\end{displaymath}
by Theorem~\ref{thm:res}, and so
\begin{displaymath}
\ind(\res(a_{\wa})) = 2a_{\wa\wa} + a_{\wa\wb} + a_{\wb\wa} + 3a_{\wa} + a_{\wb} + 1
\end{displaymath}
by Theorem~\ref{thm:ind}. Combining the above, we thus find
\begin{displaymath}
a_{\wa} \cdot a_{\wb} = a_{\wa\wb} + a_{\wb\wa} + a_{\wa} + a_{\wb} + 1.
\end{displaymath}
Note that we computed the above product in Lemma~\ref{lem:special-2}, which was used in the proof of Proposition~\ref{prop:prod-special}, so the above computation is somewhat circular; it is just intended to illustrate the procedure provided by the proof of Proposition~\ref{prop:prod-char}.
\end{example}

\begin{remark}
Lemma~\ref{lem:hamming} implies that the elements $a_{\pi(n)}$ generate $\rK$ under the induction product (which is a non-unital ring).
\end{remark}

\subsection{The projection formula for $m'$}

We would like to apply Proposition~\ref{prop:prod-char} to show that $m'=m$. For this, we need to know that $m'$ satisfies the projection formula. We now prove this.

\begin{proposition} \label{prop:comb-proj}
For $x,y,z \in \rK$, we have
\begin{displaymath}
m'(x, \ind(y \otimes z)) = \ind(m'(\res(x), y \otimes z)).
\end{displaymath}
\end{proposition}

Before giving the proof, we require two lemmas. We begin by verifying the projection formula in a simple case.

\begin{lemma} \label{lem:comb-proj-2}
For $x \in \rK$ we have $m'(x, \ind(1))=\ind(\res(x))$.
\end{lemma}

\begin{proof}
It suffices to treat the case $x=a_{\lambda}$. Let $n=\ell(\lambda)$. Define elements $x_2,x_3,x_4,x_5$ as follows:
\begin{itemize}
\item $x_2$ is the sum of all $a_{\mu}$'s where $\mu$ is obtained by inserting a $\wa$ in between two letters of $\lambda$, or at the start or end of $\lambda$; this sum has $n+1$ terms.
\item $x_3$ is defined like $x_2$, except we insert $\wb$.
\item $x_4$ is the sum of all $a_{\mu}$'s where $\mu$ is obtained by toggling a letter of $\lambda$; this sum has $n$ terms.
\item $x_5$ is the sum of all $a_{\mu}$'s where $\mu$ is obtained by deleting a letter from $\lambda$; this sum has $n$ terms.
\end{itemize}
We will show that each side of the equation in the statement of the lemma is equal to
\begin{equation} \label{eq:comb-proj}
(2n+1) x + x_2 + x_3 + x_4 + x_5,
\end{equation}
which will complete the proof.

We begin with the left side. From the computation of $\ind(1)$, we have
\begin{displaymath}
m'(x, \ind(1)) = m'(x,a_{\wa})+m'(x,a_{\wb})+x.
\end{displaymath}
We claim that this is equal to \eqref{eq:comb-proj}. The terms in $m'(x, a_{\wa})$ without collisions produces $x_2$; similarly, the terms in $m'(x, a_{\wb})$ without collisions produces $x_3$. We now consider collisions at the $i$th letter of $\lambda$, for both products at once. Together, these produce
\begin{displaymath}
a_{\lambda[1,i)} \odot (2a_{\lambda_i}+a_{\lambda_i}^{\vee}+1) \odot a_{\lambda(i,n]},
\end{displaymath}
where $a_{\lambda_i}^{\vee}$ just means that we toggle $\lambda_i$. Indeed, when we collide with $\lambda_i$, we just get $\lambda_i$, while in the other case we get $a_{\wa}+a_{\wb}+1$, which is $a_{\lambda_i}+a_{\lambda_i}^{\vee}+1$. Summing over $i$, we obtain $2nx+x_4+x_5$. This proves the claim.

We next claim that the right side of is also equal to \eqref{eq:comb-proj}. By Theroem~\ref{thm:res}, we have $\res(x)=y_1+y_2$, where
\begin{displaymath}
y_1 = \sum_{i=0}^n a_{\lambda[1,i]} \otimes a_{\lambda(i,n]}, \qquad
y_2 = \sum_{i=1}^n a_{\lambda[1,i)} \otimes a_{\lambda(i,n]}.
\end{displaymath}
Thus $y_1$ is the sum of ways of breaking $\lambda$ into two pieces between letters (or at the ends), while $y_2$ is the sum of ways of breaking $\lambda$ into two pieces by deleting letters. By Theorem~\ref{thm:ind}, we have
\begin{displaymath}
\ind(y_1) = \sum_{i=0}^n \big( a_{\lambda[1,i]} \odot (1+a_{\wa}+a_{\wb}) \odot a_{\lambda(i,n]} \big) = (n+1)x+x_2+x_3.
\end{displaymath}
Similarly,
\begin{displaymath}
\ind(y_2) = \sum_{i=1}^n \big( a_{\lambda[1,i)} \odot (1+a_{\wa}+a_{\wb}) \odot a_{\lambda(i,n]} \big) = x_5+nx+x_4.
\end{displaymath}
In each term in the above sum, one choice of $a_{\wa}$ or $a_{\wb}$ will replace the deleted letter, and these terms amount to $nx$; the terms corresponding to the other choice amount to $x_4$. This proves the claim, and completes the proof.
\end{proof}

For $\alpha, \beta \in \{\wa,\wb\}$, we define a quantity $c(\alpha, \beta) \in \rK$ as follows:
\begin{displaymath}
c(\alpha,\beta) = \begin{cases}
a_{\wa} & \text{if $\alpha=\beta=\wa$} \\
a_{\wa}+a_{\wb}+1 & \text{if $\alpha \ne \beta$} \\
a_{\wb} & \text{if $\alpha=\beta=\wb$} \end{cases}
\end{displaymath}
Note that $c(\alpha,\beta)$ is exactly the quantity used in the definition of $m'$ when a collision is encountered. The following lemma gives a recursive characterization of $m'$:

\begin{lemma} \label{lem:comb-proj-1}
Let $\lambda$ and $\rho$ be weights of lengths $r$ and~1, and let $x \in \rK$. We then have
\begin{displaymath}
m'(a_{\lambda}, a_{\rho} \odot x) = \sum_{i=0}^r a_{\lambda[1,i]} \odot a_{\rho} \odot m'(a_{\lambda(i,r]}, x) + \sum_{i=1}^r a_{\lambda[1,i)} \odot c(\lambda_i, \rho) \odot m'(a_{\lambda(i,r]}, x)
\end{displaymath}
There is a similar formula when the order of $a_{\rho} \odot x$ is reversed.
\end{lemma}

\begin{proof}
The first sum is the case where the ruffle has no collision and $i$ is taken to be maximal such that $[1,i]$ in the first variable is mapped to $[1,i]$ by the ruffle. The second sum is the case where the ruffle has a collision at $i$.
\end{proof}

\begin{proof}[Proof of Proposition~\ref{prop:comb-proj}]
Consider the following statement:
\begin{itemize}
\item[$(S_n)$] Given $x \in \rK$ and $y \otimes z \in (\rK \otimes \rK)_{\le n}$, we have
\begin{displaymath}
m'(x, \ind(y \otimes z)) = \ind(m'(\res(x), y \otimes z)).
\end{displaymath}
\end{itemize}
We prove $(S_n)$ by induction on $n$. The $n=0$ case was established in Lemma~\ref{lem:comb-proj-2}. Suppose now that $n \ge 1$ and $(S_{n-1})$ holds; we will prove $(S_n)$.

It suffices to treat the case where $x$, $y$, and $z$ are basis vectors. Moreover, we can assume that $y$ and $z$ are not both~1 since this is covered by $(S_0)$; by symmetry, we can assume that $y \ne 1$. We can thus write $x=a_{\lambda}$, $y=a_{\rho \alpha}$, and $z=a_{\beta}$, where $\rho$ is a single letter, and $\ell(\alpha)+\ell(\beta) \le n-1$. Put $w=\ind(a_{\alpha} \otimes a_{\beta})$ Note that $y=a_{\rho} \odot a_{\alpha}$, and $\ind(y \otimes z)=a_{\rho} \odot w$ by Corollary~\ref{cor:ind-cat}. Put
\begin{displaymath}
A = m'(a_{\lambda}, a_{\rho} \odot w), \qquad
B = m'(\res(a_{\lambda}), a_{\rho \alpha} \otimes a_{\beta}).
\end{displaymath}
We must prove $A=\ind(B)$.

Computing $\res(a_{\lambda})$ via Theorem~\ref{thm:res}, we find $B=X+Y$ where
\begin{align*}
X &= \sum_{j=0}^r m'(a_{\lambda[1,j]}, a_{\rho} \odot a_{\alpha}) \otimes m'(a_{\lambda(j,r]}, a_{\beta}) \\
Y &= \sum_{j=1}^r m'(a_{\lambda[1,j)}, a_{\rho} \odot a_{\alpha}) \otimes m'(a_{\lambda(j,r]}, a_{\beta})
\end{align*}
Applying Lemma~\ref{lem:comb-proj-1} to the first $m'$ in $X$, we obtain $X = X_1 + X_2$ where
\begin{align*}
X_1 &=\sum_{j=0}^r \sum_{i=0}^j (a_{\lambda[1,i]} \odot a_{\rho} \odot m'(a_{\lambda(i,j]}, a_{\alpha})) \otimes m'(a_{\lambda(j,r]}, a_{\beta}) \\
X_2 &= \sum_{j=0}^r \sum_{i=1}^j (a_{\lambda[1,i)} \odot c(\lambda_i, \rho) \odot m'(a_{\lambda(i,j]}, a_{\alpha})) \otimes m'(a_{\lambda(j,r]}, a_{\beta}).
\end{align*}
Similarly, we obtain $Y = Y_1 + Y_2$ where
\begin{align*}
Y_1 &= \sum_{j=1}^r \sum_{i=0}^{j-1} (a_{\lambda[1,i]} \odot a_{\rho} \odot m'(a_{\lambda(i,j)}, a_{\alpha})) \otimes m'(a_{\lambda(j,r]}, a_{\beta}) \\
Y_2 &= \sum_{j=1}^r \sum_{i=1}^{j-1} (a_{\lambda[1,i)} \odot c(\lambda_i, \rho) \odot m'(a_{\lambda(i,j)}, a_{\alpha})) \otimes m'(a_{\lambda(j,r]}, a_{\beta}).
\end{align*}
Reversing the order of the sums in $X_1+Y_1$ and using Theorem~\ref{thm:res}, we obtain
\begin{center}
\resizebox{0.95\hsize}{!}{%
\begin{math}
\begin{aligned}
& \sum_{i=0}^r ((a_{\lambda[1,i]} \odot a_{\rho}) \otimes 1) \odot (\sum_{j=i}^r  m'(a_{\lambda(i,j]}, a_{\alpha})) \otimes m'(a_{\lambda(j,r]}, a_{\beta})
+\sum_{j=i+1}^r m'(a_{\lambda(i,j)}, a_{\alpha})) \otimes m'(a_{\lambda(j,r]}, a_{\beta})) \\
=& \sum_{i=0}^r ((a_{\lambda[1,i]} \odot a_{\rho}) \otimes 1) \odot m'(\res(a_{\lambda(i,r]}), a_{\alpha} \otimes a_{\beta})
\end{aligned}
\end{math}
}
\end{center}
Applying $\ind$ to this, and appealing to Corollary~\ref{cor:ind-cat} and $(S_{n-1})$, we obtain
\begin{displaymath}
\ind(X_1+Y_1) = \sum_{i=0}^r a_{\lambda[1,i]} \odot a_{\rho} \odot m'(a_{\lambda(i,r]}, w)
\end{displaymath}
Similarly, we find
\begin{displaymath}
X_2+Y_2 = \sum_{i=1}^r (a_{\lambda[1,i)} \odot c(\lambda_i,\rho)) \otimes 1) \odot m'(\res(\lambda_{(i,r]}), a_{\alpha} \otimes a_{\beta}),
\end{displaymath}
and so
\begin{displaymath}
\ind(X_2+Y_2) = \sum_{i=1}^r a_{\lambda[1,i)} \odot c(\lambda_i,\rho) \odot m'(\lambda_{(i,r]}, w),
\end{displaymath}
Putting this all together, we obtain
\begin{displaymath}
\ind(B) = 
\sum_{i=0}^r a_{\lambda[1,i]} \odot a_{\rho} \odot m'(a_{\lambda(i,r]}, w)
+\sum_{i=1}^r a_{\lambda[1,i)} \odot c(\lambda_i,\rho) \odot m'(\lambda_{(i,r]}, w).
\end{displaymath}
This is equal to $A$ by Lemma~\ref{lem:comb-proj-1}.
\end{proof}

\subsection{Proof of Theorem~\ref{thm:tensor}}

We now prove $m'=m$ by showing that $m'$ satisfies conditions (a)--(c) of Proposition~\ref{prop:prod-char}. It is clear that $a_{\emptyset}$ is the unit for $m'$, and that $m'$ is commutative (since the notion of ruffle is symmetric); thus condition (a) holds. We have verified that $m'$ satisfies the projection formula (Proposition~\ref{prop:comb-proj}), which gives (b). Finally, (c) is a simple combinatorially exercise, carried out in the following lemma:

\begin{lemma}
We have $m'(a_{\pi(1)}, a_{\pi(n)})=(n+1) a_{\pi(n+1)} + na_{\pi(n)}$.
\end{lemma}

\begin{proof}
Let $\lambda=\pi(1)$ and $\mu=\pi(n)$. There are $n+1$ shuffles (i.e., ruffles without collisions), and for each one $r$ the set $P_r(\lambda,\mu)$ is the singleton containing $\pi(n+1)$. There are $n$ ruffles where the letter from $\lambda$ collides with one of the letters with $\mu$; for such ruffles $r$, the set $P_r(\lambda,\mu)$ is the singleton containing $\pi(n)$. We thus see that
\begin{displaymath}
m(a_{\lambda}, a_{\mu}) = (n+1) a_{\pi(n+1)}+n a_{\pi(n)},
\end{displaymath}
which verifies (c).
\end{proof}

\subsection{The structure of $\rK$ as a ring} \label{ss:K-ring}

Now that we have an explicit formula for the product on $\rK$, we can describe $\rK$ as ring. Recall that $\rK$ is filtered (as a ring) by the $\rK_{\le n}$. We let $\gr(K)$ be the associated graded ring. For a weight $\lambda$ of length $n$, we let $\ol{a}_{\lambda}$ be the degree $n$ element of $\gr(K)$ defined by $a_{\lambda}$; these elements form a $\bZ$-basis of $\gr(K)$.

Recall that an \defn{$(n,m)$-shuffle} is an $(n,m)$-ruffle $s \colon [n] \amalg [m] \to [\ell]$ with $\ell=n+m$, i.e., $s$ is bijective. If $s$ is such a shuffle and $\lambda$ and $\mu$ are words of lengths $n$ and $m$, then $P_s(\lambda,\mu)$ contains a unique weight, which we denote by $s(\lambda,\mu)$; it has length $n+m$.

\begin{proposition}
For weights $\lambda$ and $\mu$ of lengths $n$ and $m$, we have
\begin{displaymath}
\ol{a}_{\lambda} \cdot \ol{a}_{\mu} = \sum_s \ol{a}_{s(\lambda,\mu)},
\end{displaymath}
where the sum is over all $(n,m)$-shuffles $s$.
\end{proposition}

\begin{proof}
According to Theorem~\ref{thm:tensor}, we have
\begin{displaymath}
a_{\lambda} \cdot a_{\mu} = \sum_r \sum_{\nu \in P_r(\lambda,\mu)} a_{\nu},
\end{displaymath}
where the sum is over all $(n,m)$-ruffles $r$. We thus find
\begin{displaymath}
\ol{a}_{\lambda} \cdot \ol{a}_{\mu} = \sum_r \sum_{\substack{\nu \in P_r(\lambda,\mu), \\ \ell(\nu)=n+m}} \ol{a}_{\nu},
\end{displaymath}
as all other terms in the first formula become~0 in $\gr(K)$. If $r$ is a ruffle that is not a shuffle, then there is necessarily a collision, and so $\ell(\nu)<n+m$ for any $\nu \in P_r(\lambda,\mu)$. Thus in the above formula, it suffices to sum over shuffles. If $r$ is a shuffle then $P_r(\lambda,\mu)$ contains a unique element $r(\lambda,\mu)$, and it has length $n+m$. This gives the stated result.
\end{proof}

Recall that the \defn{shuffle algebra} on an alphabet has a $\bZ$-basis consisting of all words, and the product of two words is the sum of all shuffles; see \cite[\S 1.4]{Reutenauer}. The above proposition thus yields:

\begin{corollary}
The ring $\gr(K)$ is the shuffle algebra on the alphabet $\{\wa,\wb\}$.
\end{corollary}

Now, define an order on $\{\wa,\wb\}$ by $\wa < \wb$, and use this to lexicographically order the set of weights $\Lambda$. We say that $\lambda \in \Lambda$ is a \defn{Lyndon weight} if it is non-empty and lexicographically smaller than its proper suffixes, that is, whenever $\lambda=\mu \nu$ with $\mu$ non-empty we have $\lambda<\nu$.

\begin{corollary} \label{cor:K-poly}
The ring $\rK \otimes \bQ$ is the polynomial algebra in the elements $a_{\lambda}$ with $\lambda$ a Lyndon weight.
\end{corollary}

\begin{proof}
It is a well-known property of the shuffle algebra that $\gr(K) \otimes \bQ$ is a polynomial ring in the elements $\ol{a}_{\lambda}$ with $\lambda$ a Lyndon weight \cite[Theorem~6.1]{Reutenauer}. The stated result follows from this by a standard argument.
\end{proof}

\begin{remark} \label{rmk:gr-hopf}
The restriction map is also compatible with filrations. By Theorem~\ref{thm:res}, the induced co-multiplication on $\gr(\rK)$ agrees with the usual one on the shuffle algebra \cite[\S 1.5]{Reutenauer}. Thus $\gr(\rK)$ is isomorphic to the shuffle algebra as a bi-algebra, and hence as a Hopf algebra. In particular, writing $S$ for the antipode on $\rK$, the leading term of $S(a_{\lambda})$ is $(-1)^{\ell(\lambda)} a_{\mu}$, where $\mu$ is the reverse of $\lambda$, as this is the formula for the antipode on the shuffle algebra.
\end{remark}

\section{Adams operations} \label{s:adams}

In this section, we show that the Adams operations are trivial on the Grothendieck group $\rK$. As an application, we compute the action of Schur functors on $\rK$ (Proposition~\ref{prop:schur}).

\subsection{The $\lambda$- and Adams operations} \label{ss:lambda-bg}

Recall a \defn{$\lambda$-ring} is a commutative ring $R$ equipped with operations $\lambda^i \colon R \to R$, for $i \in \bN$, satisfying some conditions, and a \defn{special $\lambda$-ring} is one satisfying some further conditions. These definitions are reviewed in \cite[\S 5.2]{HarmanSnowden}, and general background can be found in \cite[\S 3]{Dieck} or \cite[\S I]{Knutson}. (Note: \cite{Knutson} uses the terminology ``pre-$\lambda$-ring'' and ``$\lambda$-ring'' in place of our ``$\lambda$-ring'' and ``special $\lambda$-ring.'') Given a $\lambda$-ring $R$, one defines the \defn{Adams operations} $\psi^i \colon R \to R$ in terms of the $\lambda^j$'s; see \cite[\S 3.4]{Dieck} or \cite[\S I.4]{Knutson}.

Suppose now that $\cC$ is a $k$-linear pre-Tannakian category with $k$ of characteristic~0. Then the Grothendieck group $\rK(\cC)$ admits the structure of a special $\lambda$-ring, with $\lambda^i([X])=[\lw^i(X)]$; see \cite[\S 3.3]{EtingofHarmanOstrik}. In the case where $\cC$ is the category of complex representations of a finite group, the Adams operations are given by $\chi_{\psi^i([V])}(g)=\chi_{[V]}(g^i)$, where $\chi$ denotes virtual character. In positive characteristic, $\rK(\cC)$ does not carry a natural $\lambda$-ring structure in general (see \cite[Example~3.3]{EtingofHarmanOstrik}); it does if $\cC$ is semi-simple, though even then the $\lambda$-ring structure may not be special (as the Verlinde category shows; see \cite[\S 3.3]{EtingofHarmanOstrik}).

By the above, we see that our Grothendieck group $\rK$ is a special $\lambda$-ring via $\lambda^i([V])=[\lw^i{V}]$ if $k$ has characteristic~0. In fact, this is true in positive characteristic as well; one can prove this directly, or by appealing to the fact that $\rK$ is independent of $k$ (Remark~\ref{rmk:field-ind}). The tensor product $\rK \otimes \rK$ carries a natural special $\lambda$-ring structure, as it is the Grothendieck group of $\uRep(G \times G)$.

\begin{proposition} \label{prop:res-lambda}
The map $\res \colon \rK \to \rK \otimes \rK$ is a map of $\lambda$-rings.
\end{proposition}

\begin{proof}
The map $\res$ is induced by the restriction functor
\begin{displaymath}
\uRep(G) \to \uRep(G(0)) \cong \uRep(G \times G),
\end{displaymath}
which is a tensor functor; thus the result follows.
\end{proof}

\subsection{The main theorem}

The following is the main result of this section:

\begin{theorem} \label{thm:adams}
The Adams operations are trivial, that is, $\psi^i$ is the identity for all $i$.
\end{theorem}

We again emphasize that this is an extremely special property of the Delannoy category: the only semi-simple categories we know with this property are $\uRep(G^n)$ for $n \in \bN$. In fact, the theorem is a formal consequence of results we have already proven, as the proof will show. The proof will take most of the remainder of the section.

 We observe one corollary of Theorem~\ref{thm:adams} here. Recall that a \defn{binomial ring} is a commutative ring that is $\bZ$-torsionfree and closed under the operations $x \mapsto \binom{x}{n}$ for all $n \ge 0$.

\begin{corollary} \label{cor:binom}
$\rK$ is a binomial ring, and $\lambda^i(x)=\binom{x}{i}$ for all $x \in \rK$.
\end{corollary}

\begin{proof}
This follows from results of Elliott \cite[Proposition~8.3]{Elliott} and Wilkerson \cite{Wilkerson} relating binomial rings and $\lambda$-rings.
\end{proof}

In particular, we see that one can compute the simple decomposition of exterior powers purely in terms of the ring structure on $\rK$. In fact, in \S \ref{ss:schur} we will see that this is true for arbitrary Schur functors as well.

\begin{remark} \label{rmk:adams}
Many elements of $G$ are conjugate to their powers, and this perhaps provides some intuition for Theorem~\ref{thm:adams}. To be a little more precise, suppose we had a character theory for $\ul{G}$-modules that behaved like that for finite groups. Let $V$ be a finite length $\ul{G}$-module. If $g \in G$ is a piecewise linear map then $g$ is conjugate to any power $g^i$ with $i \ge 1$. Thus $\chi_{\psi^i([V])}(g)=\chi_{[V]}(g^i)=\chi_{[V]}(g)$. Since the piecewise linear elements of $G$ are dense, we find $\chi_{\psi^i{[V]}}(g)=\chi_{[V]}(g)$ for all $g$, and so $\psi^i([V])=[V]$. We have not attempted to turn this argument into a rigorous proof.
\end{remark}

\subsection{A preliminary computation}

We begin by verifying that the Adams operations are trivial on the elements $a_{\wa}$ and $a_{\wb}$. We deduce this from the following result.

\begin{proposition} \label{prop:wedge-basic}
We have $\lw^n(L_{\wa}) \cong L_{\pi(n)}$ for all $n \ge 0$, and so $\lambda^n(a_{\wa})=a_{\pi(n)}$.
\end{proposition}

\begin{proof}
As with Proposition~\ref{prop:prod-special}, there is an illuminative proof at the representation level (see Remark~\ref{rmk:wedge-basic}), but we opt for a proof in the Grothendieck group which involves fewer subtleties. Let $z_n=\lambda^n(a_{\wa})$. We must show that $z_n=a_{\pi(n)}$. This is clear for $n \le 1$. We proceed by induction on $n$. Thus let $n \ge 2$ be given, and suppose that $z_i=a_{\pi(i)}$ holds for $0 \le i \le n-1$. We prove $z_n=a_{\pi(n)}$ by combining two constraints on $z_n$.

We have
\begin{displaymath}
\lw^{n-1}(L_{\wa}) \otimes L_{\wa} \cong L_{\pi(n-1)} \otimes L_{\pi(1)} \cong L_{\pi(n)}^{\oplus n} \oplus L_{\pi(n-1)}^{\oplus (n-1)},
\end{displaymath}
where in the first step we used the inductive hypothesis, and in the second step we used the tensor product rule (Theorem~\ref{thm:tensor}). Since $\lw^n(L_{\wa})$ is a summand of the above, we find
\begin{displaymath}
z_n=r a_{\pi(n)} + s a_{\pi(n-1)}
\end{displaymath}
for integers $0 \le r \le n$ and $0 \le s \le n-1$. This is our first constraint.

Since $\res$ is a map of $\lambda$-rings (Proposition~\ref{prop:res-lambda}), we have
\begin{displaymath}
\res(z_n) = \lambda^n(\res(a_{\wa})) = \lambda^n(a_{\wa}\otimes 1+1 \otimes a_{\wa}+1) = \sum_{i+j=n} z_i \otimes z_j + \sum_{i+j=n-1} z_i \otimes z_j
\end{displaymath}
Recall that $\Delta(x)=\res(x)-x \otimes 1 - 1 \otimes x$. Applying the inductive hypothesis to the above formula, we find
\begin{displaymath}
\Delta(z_n) = \sum_{i+j=n,\ i,j \ne 0} a_{\pi(i)} \otimes a_{\pi(j)} + \sum_{i+j=n-1} a_{\pi(i)} \otimes a_{\pi(j)}
\end{displaymath}
The right side above is equal to $\Delta(a_{\pi(n)})$ by the restriction formula (Theorem~\ref{thm:res}). Thus $z_n-a_{\pi(n)}$ is in the kernel of $\Delta$, and therefore a primitive element. It thus follows from the classification of primitive elements (Proposition~\ref{prop:prim}) that
\begin{displaymath}
z_n=a_{\pi(n)} + p(a_{\wa}+1) + q(a_{\wb}+1)
\end{displaymath}
for integers $p$ and $q$. This is our second constraint.

The two constraints combined yield $z_n=a_{\pi(n)}$. Indeed, there is no $a_{\wb}$ in the first constraint, and so $q=0$. Since $n \ge 2$, there is no~1 in the first constraint, and so $p=0$. (In fact, even for $n=1$ there is no~1 in the first constraint as $s \le n-1$.)
\end{proof}

\begin{corollary} \label{cor:adams-length-1}
We have $\psi^i(a_{\wa})=a_{\wa}$ and $\psi^i(a_{\wb})=a_{\wb}$ for all $i \ge 1$.
\end{corollary}

\begin{proof}
Recall that $\bZ\langle x \rangle$ is the ring of integer-valued polynomials, which has a $\bZ$-basis $b_n$ given by $b_n=\binom{x}{n}$. This ring is a special $\lambda$-ring with $\lambda^i(b_1)=b_i$. We have already seen that the map $\phi \colon \bZ\langle x \rangle \to \rK$ given by $\phi(b_n)=a_{\pi(n)}$ is a ring homomorphism (Corollary~\ref{cor:integer-valued}). By the previous corollary, we see that $\phi(\lambda^i(b_1))=\lambda^i(\phi(b_1))$. This implies that $\phi$ is a $\lambda$-ring homomorphism (see \cite[Lemma~5.26,~5.27]{HarmanSnowden}). Since the Adams operations are trivial on $\bZ\langle x \rangle$, it follows that they are trivial on $a_{\wa}$ as well. The case $a_{\wb}$ follows from symmetry.
\end{proof}

\begin{remark} \label{rmk:wedge-basic}
Here is the basic idea of how to prove Proposition~\ref{prop:wedge-basic} at the representation level (in characteristic $\ne 2$). Identify $L_{\wa}^{\uotimes n}$ with the submodule of $\cC(\bR^n)$ generated by functions $\phi_{\bI}$ where $\bI$ is an (unordered) tuple of type $\pi(n)$. We thus see that $\lw^n(L_{\wa})$ is generated by functions $\phi_{\bI}$, subject to relations $\phi_{\sigma \bI}=\sgn(\sigma) \cdot \phi_{\bI}$ for $\sigma \in \fS_n$. When the intervals in $\bI$ are disjoint, we may as well put them in order, and we thus obtain a copy of $L_{\pi(n)}$. If two intervals in $\bI$ are equal, then $\phi_{\bI}=0$ in $\lw^n(L_{\wa})$. One can show that any $\phi_{\bI}$ can be decomposed into $\phi_{\bJ}$'s with $\bJ$ disjoint or two components equal, which yields the result.
\end{remark}

\subsection{Proof of Theorem~\ref{thm:adams}}

We now prove the theorem. We begin with a few lemmas. For letters $a,b \in \{\wa,\wb\}$, let $S_{a,b}$ be the set of weights $\lambda$ of length $\ge 2$ that have first letter $a$ and final letter $b$. Let $\rK_{a,b}$ be the $\bZ$-submodule of $\rK$ spanned by elements $a_{\lambda}$ with $\lambda \in S_{a,b}$.

\begin{lemma} \label{lem:adams-1}
$\rK_{a,b}$ is closed under multiplication.
\end{lemma}

\begin{proof}
Consider a product $a_{\lambda} a_{\mu}$ where $\lambda$ and $\mu$ belong to $S_{a,b}$. We compute this product according to the combinatorial rule in Theorem~\ref{thm:tensor}. Consider a ruffle $r$ and the pre-product of $\lambda$ and $\mu$ associated to $r$. The first letter of the pre-product is neccessarily $a$: indeed, it must be the first letter of $\lambda$ (which is $a$), the first letter of $\mu$ (which is $a$), or the result of a collision of the first letters of $\lambda$ and $\mu$ (which can only yield $a$). Similarly, the final letter of the pre-product is $b$. The length of the pre-product is always at least the minimum length of $\lambda$ and $\mu$, and thus at least 2 in this case. Thus the pre-product starts with $a$, ends with $b$, and has length at least 2. It follows that all elements of $P_r(\lambda,\mu)$ belong to $S_{a,b}$. Thus $a_{\lambda} a_{\mu}$ is a sum of $a_{\nu}$'s with $\nu \in S_{a,b}$, and therefore belongs to $\rK_{a,b}$.
\end{proof}

Note that $\rK_{a,b}$ is not a subring of $\rK$ since it does not contain~1.

\begin{lemma} \label{lem:adams-2}
$\rK_{a,b}$ is closed under the $\lambda^i$ and $\psi^i$ operators for $i \ge 1$.
\end{lemma}

\begin{proof}
Let $\lambda$ be a non-empty weight in $S_{a,b}$. By Lemma~\ref{lem:adams-1}, $L_{\lambda}^{\uotimes i}$ decomposes into a sum of simples $L_{\nu}$ with $\nu \in S_{a,b}$ (here $i \ge 1$). Since $\lw^i(L_{\lambda})$ is a quotient of $L_{\lambda}^{\uotimes i}$, the same is true for it. Thus $\lambda^i(a_{\lambda})$ is contained in $\rK_{a,b}$. It now follows that $\rK_{a,b}$ is closed under $\lambda^i$ from Lemma~\ref{lem:adams-1} and the addition rule for $\lambda^i$. Finally, $\psi^i$ is a homogeneous degree $i$ polynomial in $\lambda^1, \ldots, \lambda^i$, and thus it too preserves $\rK_{a,b}$.
\end{proof}

\begin{proof}[Proof of Theorem~\ref{thm:adams}]
Fix $i \ge 1$. We show that $\psi^i$ is trivial on $\rK$ by inductively showing it is trivial on $\rK_{\le n}$ for all $n$. Triviality on $\rK_{\le 1}$ follows from Corollary~\ref{cor:adams-length-1}.

Suppose now that $n \ge 2$ and $\psi^i$ is trivial on $\rK_{\le n-1}$. Let $\lambda$ be a word of length $\ge 2$. By Theorem~\ref{thm:res}, we have
\begin{displaymath}
\res(a_{\lambda})=a_{\lambda} \otimes 1 + 1 \otimes a_{\lambda} + y
\end{displaymath}
where $y$ belongs to $\rK_{\le n-1} \otimes \rK_{\le n-1}$. We now applying $\psi^i$. This commutes with $\res$ by Proposition~\ref{prop:res-lambda} and fixes $y$ by the inductive hypothesis. We thus find
\begin{displaymath}
\res(\psi^i(a_{\lambda})) = \psi^i(a_{\lambda}) \otimes 1 + 1 \otimes \psi^i(a_{\lambda}) + y.
\end{displaymath}
It therefore follows that $\psi^i(a_{\lambda})-a_{\lambda}$ is a primitive element. It thus belongs to $\rK_{\le 1}$ by Proposition~\ref{prop:prim}. On the other hand, every basis vector appearing in $\psi^i(a_{\lambda})-a_{\lambda}$ is of the form $a_{\mu}$ where $\ell(\mu) \ge 2$ by Lemma~\ref{lem:adams-2}. Thus $\psi^i(a_{\lambda})-a_{\lambda}$ must vanish. It follows that $\psi^i$ is trivial on $\rK_{\le n}$, which completes the proof.
\end{proof}

\subsection{Application to Schur functors} \label{ss:schur}

We now assume $k$ is a field of characteristic~0. For a partition $\lambda$, let
\begin{displaymath}
s_{\lambda} \colon \rK \to \rK
\end{displaymath}
be the action of the Schur functor $\bS_{\lambda}$ on $\rK$. Thus if $x=[V]$ is an effective class then $s_{\lambda}(x)$ is the class of $\bS_{\lambda}(V)$. We note that $s_{\lambda}$ is not linear, but is some polynomial expression in the $\lambda^i$ operators. We now determine the action of $s_{\lambda}$ explicitly. To this end, let $p_{\lambda}$ be the integer-valued polynomial such that
\begin{displaymath}
p_{\lambda}(n) = \dim \bS_{\lambda}(k^n).
\end{displaymath}
Our main result is:

\begin{proposition} \label{prop:schur}
We have $s_{\lambda}(x) = p_{\lambda}(x)$ for all $x \in \rK$.
\end{proposition}

We note that $p_{\lambda}(x)$ is simply the result of applying the polynomial $p_{\lambda}$ to the ring element $x$. In particular, the proposition implies that the action of $s_{\lambda}$ can be obtained from just the ring structure on $\rK$, which is not true for Grothendieck groups in general.

\begin{proof}
Fix $x \in \rK$ and let $S$ be the ring of symmetric functions. We have a homomorphism of $\Lambda$-rings
\begin{displaymath}
\phi \colon S \to \rK, \qquad s_{\lambda} \mapsto s_{\lambda}(x).
\end{displaymath}
Let $P_i \in S$ denote the $i$th power sum symmetric function. We have $\phi(P_i)=\psi^i(x)$, essentially by definition of $\psi^i$. Since the Adams operations are trivial on $\rK$ (Theorem~\ref{thm:adams}), we thus see that $\phi(P_i)=x=\phi(P_1)$. Let $\ol{S}$ be the quotient of $S$ by the ideal generated by the $P_i-P_1$ for $i \ge 1$. Thus $\phi$ induces a ring homomorphism $\ol{\phi} \colon \ol{S} \to \rK$.

We have a ring homomorphism $\psi \colon S \otimes \bQ \to \bQ[t]$ induced by $\psi(s_{\lambda})=p_{\lambda}(t)$. Explicitly, for $s \in S$ and $n \in \bN$, we have
\begin{displaymath}
(\psi{s})(n)=s(1, \ldots, 1, 0, 0, \ldots),
\end{displaymath}
where there are $n$ 1's. In particular, $\psi(P_i)=t$ for all $i \ge 1$. Thus $\psi$ induces a ring homomorphism $\ol{\psi} \colon \ol{S} \otimes \bQ \to \bQ[t]$. Since $S \otimes \bQ$ is a polynomial ring in the $P_i$'s, it follows that $\ol{S} \otimes \bQ$ is a polynomial ring in $P_1$, and so $\ol{\psi}$ is a ring isomorphism. We have
\begin{displaymath}
\ol{\psi}(p_{\lambda}(P_1))=p_{\lambda}(\ol{\psi}(P_1))=p_{\lambda}(x) = \ol{\psi}(s_{\lambda}),
\end{displaymath}
and so $s_{\lambda}=p_{\lambda}(P_1)$ in $\ol{S} \otimes \bQ$.

Applying $\ol{\phi}$ to the identity just obtained, we see that $s_{\lambda}(x)=p_{\lambda}(x)$ holds in $\rK \otimes \bQ$. Since $\rK$ is torsion-free, this equality already holds in $\rK$.
\end{proof}

We can use the proposition to decompose the action of a Schur functor on $\uRep(G)$, as we now explain. Since $p_{\lambda}(t)$ is an integer-valued polynomial, we have an expression
\begin{displaymath}
p_{\lambda}(t) = \sum_{i \ge 0} c(\lambda,i) \binom{t}{i}
\end{displaymath}
for some integers $c(\lambda,i)$, almost all of which vanish. Recall that $\pi(n) \in \Lambda$ is the weight consisting of $n$ $\wa$'s.

\begin{corollary}
We have an irreducible decomposition
\begin{displaymath}
\bS_{\lambda}(L_{\wa}) = \bigoplus_{i \ge 0} L_{\pi(i)}^{\oplus c(\lambda,i)}.
\end{displaymath}
In particular, the $c(\lambda,i)$ are non-negative and the length of $\bS_{\lambda}(L_{\wa})$ is $\sum_{i \ge 0} c(\lambda,i)$.
\end{corollary}

\begin{proof}
We have $\lambda^i(x)=\binom{x}{i}$ for $x \in \rK$ (Corollary~\ref{cor:binom}), and $\lambda^i(a_{\wa})=a_{\pi(i)}$ (Proposition~\ref{prop:wedge-basic}). By Proposition~\ref{prop:schur}, we thus have
\begin{displaymath}
s_{\lambda}(a_{\wa}) = \sum_{i \ge 0} c(\lambda,i) \cdot a_{\pi(i)}.
\end{displaymath}
Since the $a_{\pi(i)}$'s are distinct basis vectors of $\rK$ and $s_{\lambda}(a_{\wa})$ is the class of the representation $\bS_{\lambda}(L_{\wa})$, it follows that $c(\lambda,i)$ is non-negative. The result follows.
\end{proof}

\begin{corollary}
 For any $V$ in $\uRep(G)$, we have
\begin{displaymath}
\bS_{\lambda}(V) \cong \bigoplus_{i \ge 0} (\lw^i{V})^{\oplus c(\lambda,i)}.
\end{displaymath}
\end{corollary}

\begin{proof}
It suffices to treat the finite length case. Again, we have $\lambda^i(x)=\binom{x}{i}$ for $x \in \rK$ (Corollary~\ref{cor:binom}). Applying the proposition to $[V] \in \rK$, we find
\begin{displaymath}
[\bS_{\lambda}(V)] = \sum_{i \ge 0} c(\lambda,i) \cdot [\lw^i(V)].
\end{displaymath}
Since $\uRep(G)$ is semi-simple, this equality in $\rK$ yields the stated isomorphism in $\uRep(G)$.
\end{proof}

\section{The path model} \label{s:delannoy}

\subsection{Delannoy paths} \label{ss:path}

Fix a vector $\ul{a} \in \bN^s$. A \defn{0-1 vector} in $\bR^s$ is a vector whose coordinates are all either~0 or~1. An \defn{$\ul{a}$-Delannoy path} is a tuple $p=(p_1, \ldots, p_{\ell})$ where each $p_i$ is a non-zero 0-1 vector and $a=p_1+\cdots+p_{\ell}$. We picture $p$ as a path in $\bR^s$ from 0 to $a$, composed of steps $p_1, \ldots, p_{\ell}$. The \defn{length} of the path $p$, denoted $\ell(p)$, is the number of steps, i.e., its length as a tuple. We write $\Gamma(\ul{a})$ for the set of $\ul{a}$-Delannoy paths.

\begin{example}
When $n=2$ there are three non-zero 0-1 vectors, namely, $(1,0)$, $(0,1)$, and $(1,1)$. Thus the above definition recovers the usual notion of planar Delannoy paths.
\end{example}

\begin{remark}
Combinatorially properties of higher dimensional Delannoy paths are studied in \cite{CDNS,Tarnauceanu}.
\end{remark}

Let $i \colon [t] \to [s]$ be an injection of finite sets, and let $i^* \colon \bR^s \to \bR^t$ be the corresponding projection. Given an $\ul{a}$-Delannoy path $p=(p_1, \ldots, p_{\ell})$, we defines its \defn{projection}, denoted $i^*(p)$, by taking the tuple $(i^*(p_1), \ldots, i^*(p_{\ell}))$ and deleting any entries that are zero. This is an $i^*(\ul{a})$-Delannoy path, possibly of shorter length. We thus have a function $i^* \colon \Gamma(\ul{a}) \to \Gamma(i^*(\ul{a}))$.

\subsection{The path category}

We now define a $k$-linear category based on Delannoy paths. Let $\cD(n,m)$ be the vector space with basis indexed by the set of paths $\Gamma(n,m)$. We write $[p]$ for the basis vector of $\cD(n,m)$ corresponding to $p \in \Gamma(n,m)$. We will define a composition law on the $\cD$'s. Before doing so, we prove a proposition.

\begin{proposition}
Let $\ul{a}=(a_1,a_2,a_3)$ be a vector in $\bN^3$, and suppose given
\begin{displaymath}
p_{12} \in \Gamma(a_1,a_2), \quad p_{23} \in \Gamma(a_2,a_3), \quad p_{13} \in \Gamma(a_1,a_3).
\end{displaymath}
Then there is at most one $\ul{a}$-Delannoy path $q$ such that
\begin{displaymath}
\pi_{1,2}(q)=p_{1,2}, \quad \pi_{2,3}(q)=p_{2,3}, \quad \pi_{1,3}(q)=p_{1,3},
\end{displaymath}
where $\pi_{i,j}$ is the usual projection map $\bR^3 \to \bR^2$.
\end{proposition}

\begin{proof}
Let $q$ and $q'$ be paths in $\Gamma(\ul{a})$ such that $\pi_{i,j}(q)=\pi_{i,j}(q')$ for all $(i,j)$ as in the proposition statement. We show that $q=q'$. It suffices to show that $q_1=q'_1$, for then we can remove the first step of each of $q$ and $q'$ and continue by induction.

First suppose that $q_1=(1,0,0)$. We then know that the first steps in $\pi_{1,2}(q)$ and $\pi_{1,3}(q)$ are $(1,0)$, and we have no information about the first step of $\pi_{2,3}(q)$. At most one of $\pi_{1,2}(q'_1)$ and $\pi_{1,3}(q'_1)$ vanishes; suppose the first one does not. Then $\pi_{1,2}(q'_1)$ is the first step in $\pi_{1,2}(q')=\pi_{1,2}(q)$, which we know to be $(1,0)$; it follows that $q'=(1,0,x)$. We thus find that $\pi_{1,3}(q'_1)=(1,x)$ does not vanish, and so it is the first step in $\pi_{1,3}(q')=\pi_{1,3}(q)$, which we know to be $(1,0)$; thus $x=0$. We have therefore shown $q'_1=(1,0,0)$ as required.

The other possibilities for $q_1$ are similar (or easier). Thus the result follows.
\end{proof}

Now, let $p_1 \in \cD(n,m)$, $p_2 \in \cD(m,\ell)$, and $p_3 \in \cD(n,\ell)$ be given. By the above proposition, there is at most one element $q \in \cD(n,m,\ell)$ with $\pi_{1,2}(q)=p_1$, $\pi_{2,3}(q)=p_2$ and $\pi_{1,3}(q)=p_3$. Put
\begin{displaymath}
\epsilon(p_1,p_2,p_3) = \begin{cases}
(-1)^{\ell(q)+\ell(p_3)} & \text{if $q$ exists} \\
0 & \text{otherwise} \end{cases}
\end{displaymath}
We define a $k$-bilinear composition law
\begin{displaymath}
\cD(n,m) \times \cD(m,\ell) \to \cD(n,\ell)
\end{displaymath}
by
\begin{displaymath}
[p_1] \circ [p_2] = \sum_{p_3 \in \Gamma(n,\ell)} \epsilon(p_1,p_2,p_3) [p_3].
\end{displaymath}
We now come to the main definition of \S \ref{s:delannoy}.

\begin{definition}
We define a $k$-linear category $\cD$ as follows. For each non-negative integer $n$ there is an object $X_n$. The space of morphisms $X_m \to X_n$ is $\cD(n,m)$. The composition law is as defined above.
\end{definition}

\begin{definition}
The \defn{Delannoy algebra} $\cD(n)$ is the algebra $\End_{\cD}(X_n)$.
\end{definition}

Note that we have not yet shown that composition is associative or admits identity elements, so the above definition is somewhat incomplete. It is possible to establish these properties directly; however, we will deduce them in the course of proving Theorem~\ref{thm:equiv}.

\subsection{The equivalence theorem}

Given an $(n,m)$-Delannoy path $p$, let $O_p \subset \bR^{(n)} \times \bR^{(m)}$ be the corresponding orbit, as defined in Proposition~\ref{prop:orbit}, and let $A_p$ be the characteristic function of $O_p$, thought of as an $\bR^{(n)} \times \bR^{(m)}$ matrix. We thus have a natural $k$-linear map
\begin{displaymath}
\cD(n,m) \to \Hom_{\ul{G}}(\cC(\bR^{(m)}), \cC(\bR^{(n)})), \qquad [p] \mapsto A_p.
\end{displaymath}
As this map is a bijection on bases, it is an isomorphism of $k$-vector spaces. We now verify that it is compatible with composition.

\begin{lemma} \label{lem:equiv}
Given $p_1 \in \cD(n,m)$ and $p_2 \in \cD(m,\ell)$, we have
\begin{displaymath}
A_{[p_1] \circ [p_2]}=A_{p_1} A_{p_2},
\end{displaymath}
where on the right side we use matrix multiplication.
\end{lemma}

\begin{proof}
Fix an $(n,\ell)$-Delannoy path $p_3$, and let us compute the coefficient $c$ of $A_{p_3}$ in the matrix product $A_{p_1} A_{p_2}$. By definition, this is the value of $A_{p_1} A_{p_2}$ at a point $(z,x)$ in $O_{p_3}$. Thus for such a point $(z,x)$ we have
\begin{displaymath}
c = \int_{\bR^{(m)}} A_{p_1}(z,y) A_{p_2}(y,x) dy.
\end{displaymath}
Integrating both sides over $(z,x) \in O_{p_3}$, we find
\begin{displaymath}
c \vol(O_{p_3}) = \int_{\bR^{(n)} \times \bR^{(m)} \times \bR^{(\ell)}} A_{p_1}(z,y) A_{p_2}(y,z) A_{p_3}(z,x) dx dy dz.
\end{displaymath}
The right side above is $\vol(X)$, where $X$ is the set of points $(x,y,z)$ such that $(z,y) \in O_{p_1}$, $(y,z) \in O_{p_2}$, and $(z,x) \in O_{p_3}$. In other words, $X$ is the union of those orbits that project to $O_{p_1}$, $O_{p_2}$, and $O_{p_3}$ under the three projections.

Now, the correspondence between orbits and Delannoy paths in Propsoition~\ref{prop:orbit} extends naturally to orbits on a threefold product. Moreover, this classification is compatible with projections, in the sense that a projeciton of an orbit is corresponds to the projection of the Delannoy path. Furthermore, the volume of an orbit is easily seen to be $(-1)^{\ell}$, where $\ell$ is the length of the Delannoy path.

It follows that if $q$ as in the definition of $\epsilon(p_1,p_2,p_3)$ exists then $X=O_q$, and otherwise $X$ is empty. In the first case, $c=\vol(O_q)/\vol(O_{p_3})=(-1)^{\ell(q)+\ell(p_3)}$, and in the second case $c=0$. In all cases, $c=\epsilon(p_1,p_2,p_3)$, which completes the proof.
\end{proof}

We can now prove our main result about $\cD$:

\begin{theorem} \label{thm:equiv}
We have the following:
\begin{enumerate}
\item The composition law in $\cD$ is associative and has identity elements.
\item There is a fully faithful functor $\Phi \colon \cD \to \uRep(G)$ defined on objects by $\Phi(X_n)=\cC(\bR^{(n)})$ and on morphisms by $\Phi([p])=A_p$.
\item $\Phi$ identifies $\cD$ with the full subcategory of $\uRep(G)$ spanned by the $\cC(\bR^{(n)})$'s.
\item $\Phi$ identifies the additive envelope of $\cD$ with $\uPerm(G)$.
\item $\Phi$ identifies the additive--Karoubian envelope of $\cD$ with $\uRep^{\rf}(G)$.
\end{enumerate}
\end{theorem}

\begin{proof}
(a) Lemma~\ref{lem:equiv} shows that composition in $\cD$ matches with composition of matrices in $\uPerm(G)$, and is therefore associative and has identity elements.

(b) Lemma~\ref{lem:equiv} shows this is a functor, and we already know it induces isomorphism on $\Hom$ spaces.

(c) Let $\cC_0$ be the full subcategory of $\uRep(G)$ spanned by the objects $\cC(\bR^{(n)})$. As $\Phi \colon \cD \to \cC_0$ is fully faithful and essentially surjective, it is an equivalence.

(d) This follows since $\uPerm(G)$ is the additive envelope of $\cC_0$.

(e) This follows since $\uRep^{\rf}(G)$ is the additive--Karoubian envelope of $\cC_0$.
\end{proof}

\begin{remark}
The equivalence theorem, combined with results on $\uRep(G)$ from \cite{HarmanSnowden}, implies that the Delannoy algebra $\cD(n)$ is semi-simple, and, more generally, that the additive--Karoubian envelope of $\cD$ is semi-simple. It is possible to prove these statements directly, without relying on \cite{HarmanSnowden}.
\end{remark}

\subsection{The tensor product}

Let $\cD^+$ be the additive envelope of $\cD$. We have seen in Theorem~\ref{thm:equiv} that $\cD^+$ is equivalent to $\uPerm(G)$. We can thus transfer the tensor product from $\uPerm(G)$ to $\cD^+$. We now say a few words about this.

First of all, in $\uPerm(G)$ we have
\begin{displaymath}
\cC(\bR^{(n)}) \uotimes \cC(\bR^{(m)}) = \cC(\bR^{(n)} \times \bR^{(m)}) = \bigoplus_{p \in \Gamma(n,m)} \cC(O_p),
\end{displaymath}
where here $O_p$ is the orbit corresponding to the Delannoy path $p$. It follows from the proof of Proposition~\ref{prop:orbit} that $O_p$ is isomorphic to $\bR^{(\ell)}$, where $\ell=\ell(p)$. We thus obtain the following description of the tensor product on $\cD^+$ at the level of objects:
\begin{displaymath}
X_n \uotimes X_m = \bigoplus_{p \in \Gamma(n,m)} X_{\ell(p)}.
\end{displaymath}
One can describe $\uotimes$ on morphisms, as well as the associator for $\uotimes$, directly in terms of Delannoy path but we will not do that here. We note that the associator is non-trivial.

There is an alternate approach that is a bit more clean. We define a tensor category $\cD'$ as follows. The objects are symbols $X_{\ul{n}}$ where $\ul{n}$ is a tuple of non-negative integers; intuitively, $X_{\ul{n}}$ is the tensor product of the $X_{n_i}$'s. Morphisms $X_{\ul{n}} \to X_{\ul{m}}$ are given by linear combinations of $\ul{n} \ul{m}$ Delannoy paths, where $\ul{n} \ul{m}$ indicates concatenation of tuples. Composition is defined analogously to $\cD$. The tensor product is defined on objects by $X_{\ul{n}} \uotimes X_{\ul{m}}=X_{\ul{n} \ul{m}}$, and by a simple formula on morphisms. The associator is trivial. One can then show that the additive envelopes of $\cD$ and $\cD'$ agree, by decomposing $X_{\ul{n}}$ with explicit idempotents. This allows one to transfer the tensor product from $\cD'$ to $\cD^+$.

\end{document}